\documentclass[preprint, 12pt]{elsarticle}

\usepackage{amsmath}
\usepackage{amsfonts}
\usepackage{amssymb}
\usepackage{mathtools}
\usepackage{graphicx}
\usepackage{color}
\usepackage{float}
\usepackage{fullpage}
\usepackage{dsfont}
\usepackage{epsfig}
\usepackage{setspace}
\usepackage{multicolumn}
\usepackage{algorithm2e}
\usepackage{longtable}
\usepackage{caption}
\usepackage{algcompatible}
\usepackage{booktabs}
\usepackage{url}
\usepackage{natbib}
 \bibpunct[, ]{(}{)}{,}{a}{}{,}%
\newtheorem{thm}{Theorem}[section]

\newtheorem{example}[thm]{Example}

\newtheorem{remark}[thm]{Remark}

\newcommand{\beq}{\begin{equation}}
\newcommand{\eeq}{\end{equation}}
\newcommand{\beqa}{\begin{eqnarray}}
\newcommand{\eeqa}{\end{eqnarray}}
\newcommand{\beqas}{\begin{eqnarray*}}
\newcommand{\eeqas}{\end{eqnarray*}}
\newcommand{\bi}{\begin{itemize}}
\newcommand{\ei}{\end{itemize}}

\newcommand{\lam}{{\lambda}}

\newcommand{\ignore}[1]{}
\def\defi{\vcentcolon=}

\journal{Computers and Operations Research}

\begin{document}

\begin{frontmatter}

\title{New Heuristics for the Operation of an Ambulance Fleet under Uncertainty}

\author[inst1]{Vincent Guigues}

\affiliation[inst1]{organization={School of Applied Mathematics, FGV},
            addressline={Praia de Botafogo 190}, 
            city={Rio de Janeiro},
            country={Brazil}}

\author[inst2]{Anton J. Kleywegt}
\author[inst3]{Victor Hugo Nascimento}

\affiliation[inst2]{organization={Georgia Institute of Technology},
            city={Atlanta},
            state={Georgia},
            postcode={30332-0205}, 
            country={USA}}

\affiliation[inst3]{organization={School of Applied Mathematics, FGV},
            addressline={Praia de Botafogo 190},
            city={Rio de Janeiro},
            country={Brazil}}

\begin{abstract}
The operation of an ambulance fleet involves ambulance selection decisions about which ambulance to dispatch to each emergency, and ambulance reassignment decisions about what each ambulance should do after it has finished the service associated with an emergency.
For ambulance selection decisions, we propose four new heuristics: the Best Myopic (BM) heuristic, a NonMyopic (NM) heuristic, and two greedy heuristics (GHP1 and GHP2).
For ambulance reassignment decisions, we propose several strategies to choose which emergency in queue to send an ambulance to or which ambulance station to send an ambulance to when it finishes service.
These heuristics are also used in a rollout approach: each time a new decision has to be made (when a call arrives or when an ambulance finishes service), a two-stage stochastic program is solved.
The proposed heuristics are used to efficiently compute the second stage cost of these problems.
We apply the rollout approach with our heuristics to data of the Emergency Medical Service (EMS) of a large city, and show that these methods outperform other heuristics that have been proposed for ambulance dispatch decisions.
We also show that better response times can be obtained using the rollout approach instead of using the heuristics without rollout.
Moreover, each decision is computed in a few seconds, which allows these methods to be used for the real-time management of a fleet of ambulances.
\end{abstract}

\begin{keyword}
 Stochastic programming \sep heuristics \sep ambulance dispatch
\sep emergency medical services
\MSC 90C15 \sep 90C90
\end{keyword}

\end{frontmatter}






\section{Introduction}
\label{sec:intro}

The efficient management of an ambulance fleet is of great importance for Emergency Medical Services (EMSs) and the customers whom they serve.
One part of such management is concerned with the location of ambulances and the dynamic allocation of ambulances to emergency calls.
The two main types of ambulance dispatch decisions are ambulance selection decisions and ambulance reassignment decisions.
When an emergency call is received by an EMS call center, then an ambulance selection decision is made, which determines which ambulance is dispatched to the new emergency, or whether the emergency is placed in a queue to be served later.
When an ambulance becomes available after completing service, then an ambulance reassignment decision is made, which determines what to do with that ambulance, including what emergency in queue to send the ambulance to or what station to send the ambulance to for staging until it is dispatched to an emergency.

\if{
for a setting with one priority level in \cite{jagt:17a}, two priority levels in \cite{band:12,mcla:13a,mcla:13b,mayo:13,band:14} and three priority levels in \cite{ande:07,lisay:16}.
Optimization models for ambulance selection have also been proposed in \cite{schm:12]
For ambulance reassignment, heuristics have been proposed in \cite{lees:12}.
}\fi


\paragraph{Ambulance operations}
We give a brief overview of ambulance operations.
When an emergency call arrives, an ambulance is sent to the location of the emergency immediately, or the emergency is placed in a queue, to be served later.
When an ambulance is sent to the location of an emergency, the ambulance goes through some or all of the following steps, called \emph{trips} in this paper, during the service of the emergency:
\begin{itemize}
\item[a)] ambulance travels to the scene of the emergency;
\item[b)] ambulance performs service at the scene of the emergency;
\item[c)] ambulance transports the patient to a hospital;
\item[d)] ambulance stays at the hospital to transfer the patient to the hospital;
\item[e)] ambulance goes to a cleaning station to clean the ambulance;
\item[f)] ambulance stays at the cleaning station for cleaning;
\item[g)] ambulance goes to an ambulance station for staging.
\end{itemize} 

An ambulance attending an emergency may or may not need to transport the patient(s) to a hospital.
Also, after service an ambulance is cleaned.
Cleaning on the spot (at the scene of the emergency or at the hospital) may be sufficient, or the ambulance may need to go to a cleaning station for a more thorough cleaning.
Therefore, we consider the following four sequences of ambulance trips during the service of an emergency (see Figure \ref{figgroupcalls} for a graphical representation of these trips):
\begin{itemize}
\item[($C_{1}$)] a), b), c), d), e), f), g) if patient(s) are transported to a hospital and the ambulance travels to a cleaning station;
\item[($C_{2}$)] a), b), c), d), g) if patient(s) are transported to a hospital but the ambulance does not travel to a cleaning station;
\item[($C_{3}$)] a), b), e), f), g) if patient(s) are not transported to a hospital but the ambulance travels to a cleaning station;
\item[($C_{4}$)] a), b), g) if patient(s) are not transported to a hospital and the ambulance does not travel to a cleaning station.
\end{itemize}
Note that in the case of $C_{1}$ and $C_{3}$ the ambulance reassignment decision is made after cleaning at the cleaning station, in the case of $C_{2}$ the ambulance reassignment decision is made after the patient(s) have been transferred to the hospital, and in the case of $C_{4}$ the ambulance reassignment decision is made after completion of service at the scene of the emergency.
Also note that an ambulance is available for dispatch to an emergency while traveling to an ambulance station, that is, during trip g).

\if{
The main objective of this paper is to propose
new heuristics for the operation of an ambulance
fleet under uncertainty with the problem assumptions mentioned before and to integrate them in a rolling horizon approach. The heuristics will satisfy
the following requirements:
\begin{itemize}
\item[1)] the decisions of the four groups described above
are taken based on the information available when these
decisions are taken and may only use statistical assumptions on 
the process of future calls (the exact locations, types, and time
instants
of future calls are not known);
\item[2)] the heuristics should be computed quickly (due to the nature
of the problem, new ambulance rides must be decided quickly);
\item[3)] the heuristics should provide quick response times when applied
and tailored to given emergency health services.
\end{itemize}
}\fi

\if{
The management of ambulance operations requires decisions to be made quickly, with important consequences for the mortality and morbidity of people.
Examples of such decisions are the following:
\begin{enumerate}
\item
When a call requesting emergency aid arrives, it has to be decided whether to dispatch an ambulance to the emergency location immediately, or whether to place the call in a queue of waiting calls.
If it is decided to dispatch an ambulance to the emergency location, then it also has to be decided which ambulance to dispatch to the scene.
\item
When an ambulance finishes its task with an emergency (either at the emergency location or at a hospital after dropping off patients), it has to be decided what the ambulance should do next.
If there are calls waiting in queue, and it is decided to send the ambulance to a call waiting in queue, then it also has to be decided to which call to send the ambulance.
If it is decided not to send the ambulance to a call waiting in queue, then it also has to be decided to which ambulance station, where ambulances wait for their next assignments, to send the ambulance.
\end{enumerate}
}\fi

These ambulance dispatch decisions involve trade-offs between current and future consequences.
For example, an ambulance may be dispatched to a current, less urgent, emergency, and in the process it may not be available for a future, more urgent, emergency.
These trade-offs are challenging for various reasons.
First, current emergencies are known (at least the location and something about the nature of the emergency is usually known), whereas typically future emergencies are not known.
However, it is known that certain types of emergencies tend to occur with greater frequencies in specific parts of the city and during specific times of the week, and therefore it is prudent to send more available ambulances to these parts during these times.
For example, penetrating trauma and traffic incidents tend to occur with greater frequency in some parts of the city on Friday evenings.
Second, different ambulances and crews have different capabilities to improve the outcomes for patients.
For example, many EMSs have basic life support (BLS) and advanced life support (ALS) ambulances.
Some also have other types of ambulances, such as intermediate life support ambulances, stroke units, and motorcycles.
Different crew members also have different capabilities, including Emergency Medical Technicians (EMTs), Advanced Emergency Medical Technicians (AEMTs), paramedics, and physicians with various specialties.
Third, the consequences of response time and ambulance/crew capabilities are different for different emergencies, and for many types of emergencies are not yet well known.
For example, it is well known that for cardiac arrest the CPR response time is crucial and is more important than the advanced capabilities of the ambulance and crew.
For some emergencies, it is known that advanced capabilities, such as the ability to administer intravenous treatment or specific pharmaceuticals, are more important.
And for many emergencies there is a trade-off between response time and ambulance/crew capabilities, so that given a choice of dispatching a BLS ambulance that is 10 minutes from the emergency or an ALS ambulance that is 20 minutes from the emergency, either may be better than the other depending on the type of emergency.

\paragraph{Related literature} Many ambulance-related optimization problems have been proposed in the literature.
Most of these problems address the location of ambulance stations or the assignment of ambulances to stations.
For example, the location set covering problem (LSCP) determines the locations of the minimum number of facilities that covers a given set of demand points, and was proposed by \cite{tore:71}.
Variations of the LSCP have been considered by \cite{berl:74,schi:79}, and \cite{dask:81}.
Another example is the maximal covering location problem (MCLP), that determines the locations of a given number of facilities to maximize the weighted set of demand points that is covered by the facilities, and was proposed by \cite{chur:74}.
Variations of the MCLP have been considered by \cite{dask:83,hoga:86,reve:89,repe:94,gend:97,ingo:08,erku:09,schm:10}, and \cite{sore:10}.
Various stochastic models, including queuing models and simulations, have been proposed by \cite{volz:71,swov:73a,swov:73b,fitz:73,lars:74,lars:75,hill:84,jarv:85,gold:90b,gold:91a,gold:91b,burw:93} and \cite{rest:09} for evaluating location decisions for stations and ambulances.
The solutions to the location problems mentioned above are sometimes used for the assignment of ambulances to stations. That is, each ambulance is assigned to a home station. Moreover, when an ambulance becomes available and is not dispatched to an emergency waiting in queue, then the ambulance is sent to its home station, as proposed in \cite{gold:90b,hend:04,rest:09,band:12,knig:12,maso:13,mayo:13}, and \cite{band:14}. A shortcoming of the home station approach is that even if coverage is optimal when all ambulances are at their home stations, when some ambulances are busy, the coverage can be far from optimal given the available ambulances.
To improve coverage when some ambulances are busy, it has been proposed to relocate ambulances among stations.
Various ambulance relocation problems have been considered by \cite{gend:01,brot:03,gend:06,nair:09,maxw:10,schm:12,alan:13,maxw:13,maxw:14,dege:15,jagt:15,bela:16,vanbarn:16}, and \cite{vanbarn:17}.

There are many important ambulance optimization problems besides the location of ambulance stations and the assignment of ambulances to stations.
One of the most important decisions in ambulance operations is to select the ambulance and crew to dispatch to an emergency.
The most popular dispatch policy in the literature is the simple closest-available-ambulance rule, used by \cite{hend:99,hend:04,maxw:09,maxw:10,maxw:13}, and \cite{alan:13}.
As the name indicates, the closest-available-ambulance rule dispatches the available ambulance that is closest (in terms of time or distance or some other metric) to the emergency.

A few papers have proposed alternatives to the closest-available-ambulance rule.
\cite{ande:07} proposed a measure of ``preparedness'' for each zone that measures how well available ambulances can respond to expected future emergencies in the zone.
They considered the following heuristic to select the ambulance to dispatch to an emergency for a system with three priority levels.
For priority~1 emergencies, the closest available ambulance is dispatched.
For priority~2 and~3 emergencies, the ambulance with expected travel time less than a specified threshold that will result in the least decrease in the minimum preparedness measure over all zones is dispatched.
\cite{lees:11} showed that the preparedness measure proposed by \cite{ande:07} resulted in worse performance than the closest-available-ambulance rule.
Then \cite{lees:11} proposed two modifications of the preparedness-based dispatching rule.
The first modification dispatches the available ambulance that maximizes the minimum preparedness measure over all zones divided by the travel time from the ambulance to the emergency location.
The second modification replaces the minimum preparedness measure over all zones in the calculations with other aggregates of the preparedness measures of different zones.

\citet{lees:17} pointed out a shortcoming of the alternative aggregates above, and proposed an ambulance selection method that chooses the ambulance that minimizes a weighted average response time.
\citet{lees:14} proposed a rule that takes both available and busy ambulances into account when making ambulance dispatch decisions.

\cite{mayo:13} proposed a method to partition the service region into districts and assign a number of ambulances to each district.
They used simulation to compare the performance of four dispatch policies for a setting with two priority levels; two types of policies specifying ambulance selection decisions if there is an ambulance available in the same district as the emergency, combined with two types of policies specifying ambulance selection decisions if there is no ambulance available in the same district.
If there is an ambulance available in the same district as the emergency, then the first type of policy dispatches the closest available ambulance within the same district, and the second type of policy applies a heuristic ambulance selection rule to each district.
If there is no ambulance available in the same district as the emergency, then the first type of policy assumes that an alternate emergency response, for example provided by the fire department, is automatically dispatched within the same district, and the second type of policy dispatches an ambulance from another district using a preference list of ambulances.
Similarly, \cite{band:14} proposed an ambulance dispatching heuristic, and used simulation to compare the performance of the heuristic and the closest-available-ambulance dispatching rule for a setting with two priority levels.
\cite{lisay:16} used simulation to compare the closest-available-ambulance selection rule with a policy that dispatches the closest available ambulance to priority~1 emergencies, and the ambulance within a specified response time radius which has the least utilization to priority~2 and~3 emergencies.
\cite{jagt:17a} and \cite{jagt:17b} compared two ambulance selection policies with the closest-available-ambulance policy.
\cite{schm:12} used approximate dynamic programming to control ambulance selection and ambulance reassignment for an EMS in Vienna. 
In \cite{band:12}, the ambulance selection problem was also modeled as a continuous-time Markov decision process, which can be solved if the number of zones and number of ambulances are sufficiently small.
In \cite{guiklevhn2022}, a two-stage stochastic optimization model, that incorporated ambulance selection and ambulance reassignment decisions, was proposed.
An advantage of the heuristic-based solution methods proposed in our paper is that the decisions are computed much faster than using the methods proposed in \cite{schm:12,band:12} and \cite{guiklevhn2022}.

Note that \cite{ande:07,mayo:13} and \cite{band:14} made provision for different emergency types in the form of priority levels, whereas \cite{lees:11} made provision for only one emergency type.
Also, none of these papers made provision for different ambulance and crew capabilities.
However, different ambulances are equipped differently, and different crew members have different training, experience, and skills.
Many types of emergencies are distinguished by  Emergency Medical Dispatchers (EMDs), and types of emergencies differ from each other not just in terms of ``priority'' or response time urgency, but also in terms of the ambulance and crew capabilities needed.
In contrast with previous papers, we make provision for different emergency types as distinguished by EMDs as well as for different ambulance and crew capabilities.

\cite{lees:11} made provision for queueing of emergency requests if no ambulances are available, whereas \cite{mayo:13} and \cite{band:14} assumed that emergency requests that arrive when no ambulances are available are transferred to another service.
In contrast with previous papers, we make provision for queueing of emergency requests not only if no ambulances are available, but also as a conscious decision when a low urgency emergency call arrives and few ambulances are available.

Among the papers that make provision for queueing of emergency requests, there are different decision rules for choosing which emergency in queue to dispatch an ambulance to when the ambulance becomes available.
When an ambulance becomes available and there are emergency requests in queue, \cite{lees:11} dispatches the ambulance to the emergency in queue with location that is closest to the ambulance.
\cite{lees:12} proposed a rule that takes into account both the travel times between the ambulance and emergencies in queue, as well as a centrality measure of each emergency in queue.
The rule gives preference to an emergency in queue that is close to other emergencies in queue, so that if the ambulance does not have to take the patient of the first emergency to a hospital, then the ambulance will be close to other emergencies in queue.
In contrast with previous papers, when an ambulance becomes available and there are emergency requests in the queue, we make provision for either dispatching the ambulance to an emergency in the queue or not dispatching the ambulance to an emergency in the queue, and if the ambulance is dispatched to an emergency in queue we make provision for considering the location of the emergency relative to the ambulance, the age of the emergency, the type of the emergency, as well as the ambulance and crew capabilities when deciding which emergency in queue to dispatch the ambulance to.

When an ambulance becomes available and there are no emergency requests in queue, then \cite{lees:11} keeps the ambulance in place (for example, at the hospital where the ambulance delivered a patient).
\cite{mayo:13} and \cite{band:14} assigned a home station to each ambulance, and when an ambulance becomes available and there are no emergency requests in queue, the ambulance returns to its home station.
We use a two-stage approach that chooses where to send an ambulance that has just become available to minimize the expected sum of current cost and future cost under a chosen rollout policy.
We add the following important detail to the problem.
After an ambulance has delivered a patient to a hospital, the ambulance may be cleaned.
If the cleaning can be done at the hospital, then the cleaning is completed before the ambulance is dispatched, either to an emergency in queue, or to an ambulance station.
If the cleaning cannot be done at the hospital, then the ambulance is first dispatched to an ambulance cleaning station.
If the EMS has multiple ambulance cleaning stations, then one of these cleaning stations is chosen.
Then, after the ambulance has been cleaned at the chosen cleaning station, the ambulance either stays at the cleaning station, or the ambulance is dispatched, either to an emergency in queue, or to another ambulance station.

\paragraph{Contributions}
\begin{enumerate}
\item
As pointed out in the literature review, many methods proposed in the literature consider only the ambulance selection decision, and ignore or oversimplify the ambulance reassignment decision, for example by assuming that there are no emergencies in queue and that ambulances always return to pre-assigned home stations after completion of service for an emergency.
In contrast, our methods make provision for both ambulance selection decisions and ambulance reassignment decisions, and make provision for emergencies in queue, and do not restrict ambulance reassignment decisions to returning ambulances to home stations.
\item
Our methods make provision for different types of emergencies, such as the emergency types in the Medical Priority Dispatch System (MPDS).
This is more general than assigning all emergencies to a (small) number of priority levels.
For example, the MPDS classification incorporates both the importance of response time as well as the importance of ambulance and crew capabilities, whereas priority levels refer to response time importance only.
\item
Our methods make provision for different ambulance and crew capabilities.
For example, many EMSs have both basic life support (BLS) and advanced life support (ALS) ambulances, and some have additional types of ambulances too, such as stroke units or helicopters.
Also, different crew members have different qualifications, such as emergency medical technician (EMT), paramedic, or physician.
Different crew members also have different amounts of experience, and different skills, such as crowd control skills or different language skills.
\item
The two-stage versions of our methods consider the trade-off between the current and future consequences of ambulance dispatch decisions.
The current consequences of sending a specific ambulance to an emergency include the penalized response time (penalized according to the type of emergency) and the quality of the match between ambulance and crew capabilities and emergency type.
The future consequences of sending a specific ambulance to an emergency include the impact of the decision on future ambulance availability, as measured by the penalized response times and the quality of the matches for emergencies in several future scenarios.
\end{enumerate}

Our methods satisfy the following requirements:
\begin{itemize}
\item[1)] the decisions are made based on the information available at the time of the decisions as well as a probability model of the process of future emergencies (the exact locations, types, and time instants of future emergencies are not known);
\item[2)] the decisions can be computed quickly;
\item[3)] the heuristics result in good performance, that is, quick response times and responses that are tailored to the emergency types.
\end{itemize}

Instead of using heuristics, one could consider a multistage stochastic optimization problem, or its Sample Average Approximation (SAA).
An example of such a problem is given in \cite{guiklevhn2022}.
The resulting SAA problem can have billions of integer variables, and such a problem cannot be solved sufficiently fast in practice.
The fact that decisions are computed quickly with our heuristics is an important advantage over the SAA approach.

Our code (available both in C++ and Matlab) for all the methods proposed in this paper, as well as our implementation of the closest available ambulance heuristic and of the 7 heuristics from \citealt{ande:07,lees:12,mayo:13,band:14,lees:14,jagt:15,lees:17}, are available on GitHub at \url{https://github.com/vguigues/Heuristics_Dynamic_Ambulance_Management}.

We also developed the webpage \url{http://samu.vincentgyg.com/} to provide management and visualization tools for EMSs and researchers.
On this website, the user can find emergency data of the Rio de Janeiro EMS.  The arrivals of emergencies and movement of ambulances under the various heuristics can also be visualized 
on this website
on a map of the service region, and various performance metrics such as distributions of response times for different emergency types can be visualized.

\paragraph{Organization of the paper}
The paper is organized as follows. 
In Section \ref{sec:detcase}, we explain our notation and define the performance metric for the deterministic case.
In Section \ref{smodel}, we describe the stochastic optimization model.
In the subsequent sections, we explain various heuristics for both ambulance selection decisions and ambulance reassignment decisions.
More precisely, for ambulance selection, we describe the Best Myopic (BM) heuristic in Section~\ref{bestmyopic}, a NonMyopic nonanticipative policy in Section~\ref{nonmyopic}, and two greedy heuristics (GHP1 and GHP2) in Section~\ref{sec:prior}.
In Section \ref{sec:examples}, we illustrate the use of our
heuristics on several simple examples, to better understand how
these heuristics work. The Best Station Rule (BSR) for the ambulance reassignment problem is presented in Section~\ref{sec:newbase}.
Finally, numerical results using data for  Rio de Janeiro EMS are presented in Section~\ref{sec:numsec}, and it is shown that our methods outperform 8 heuristics proposed in the literature.

\if{
\section{Overview of Emergency Medical Service Operations}
\label{sec:EMS operations}

In this section, we give an overview of Emergency Medical Service (EMS) operations.
First, every EMS is unique, and thus some operational details of an EMS may differ from this description.
An EMS operates various ambulances with various crew members.
As mentioned above, different ambulance/crew combinations have different capabilities, and it is important to take these capabilities into account when choosing which ambulance to dispatch to a particular emergency.
An EMS also operates ambulance stations where crew members and ambulances can wait to be dispatched.
Stations can vary from parking spots with minimal additional features, to well-equipped facilities where crew members can eat, relax, and sleep, and where ambulances can undergo thorough cleaning and maintenance.
An EMS also has a call center that receives emergency calls.
The organization that operates the call center and the organization that operates the ambulances may be the same or different, and in the USA either can be part of a local government or can be outsourced to a private contractor.
The Emergency Medical Dispatchers (EMDs) at a call center answer the calls and ask the callers a sequence of scripted questions.
More specifically, the questions form a tree that branches according to the answers received.
Two widely used systems of questions are the Medical Priority Dispatch System (MPDS), and the Association of Public-Safety Communications Officials (APCO) system.
These systems classify emergencies based on the caller's answers into 30--40 chief complaint types, each of which is further subdivided into 4--17 subtypes that determine the combination of ambulance/crew capability needed and response time urgency.
After an EMD has classified an emergency, the EMD dispatches ambulances if deemed appropriate (or instructs someone else to dispatch ambulances), and gives pre-arrival instructions to the caller.
Any ambulance can be dispatched to an emergency, even if the ambulance is on its way to another emergency or on its way to an ambulance station, but typically an ambulance that is busy with an emergency is not dispatched to another emergency.
Not all calls require an ambulance to be dispatched (calls that do not require an ambulance to be dispatched are not included in our model).
Also, as mentioned above, a call may be placed in queue, that is, an ambulance is not dispatched immediately (within a minute or two) after the call has been received, but rather, the dispatcher waits until a sufficient number of ambulances are available before an ambulance is dispatched to the emergency location.

When an ambulance arrives at an emergency location, the crew members perform the tasks at hand according to their training.
If one or more patients need transport to a hospital, then the ambulance(s) and crew transport these patients to the hospital deemed best under the circumstances.
The best hospital depends on both the location of the emergency and the type of emergency.
Sometimes permission is obtained from the hospital's emergency department before the patients are taken to the hospital, and a hospital sometimes denies permission (``diverts'' the ambulance), but permission is not always asked.
After an ambulance arrives with a patient at a hospital, the patient is transferred to the hospital's emergency department.
The transfer may be quick, or may take a substantial amount of time, depending on the emergency department personnel.
After the patient has been transferred, the ambulance crew have to clean the ambulance and complete their report.
If the ambulance requires relatively light cleaning, then the crew can do the cleaning on the spot, and thereafter the ambulance is ready to be dispatched to a call in queue, or to be sent to an ambulance station.
If the ambulance requires major cleaning, then the ambulance is taken to an ambulance station where such major cleaning can be done, and thereafter the ambulance is ready to be dispatched to a call in queue, or it waits at the ambulance station to be dispatched.

}\fi

\section{Notation and objective function for the deterministic case}
\label{sec:detcase}

\subsection{Notation and terminology}

The following notation and terminology will be used in the paper.
Parameter and variable names in this paper match those used in our code, to ease the understanding of the code (available on GitHub at \url{https://github.com/vguigues/Heuristics_Dynamic_Ambulance_Management}).
We also refer to \cite{websiteambrouting24} that contains a detailed description of the simulation used to evaluate policies for ambulance operations.

\begin{figure}
\centering
\begin{tabular}{c}
\includegraphics[scale=1.07]{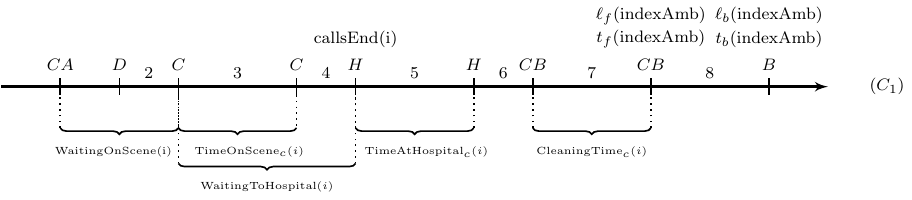}\\
\includegraphics[scale=1.12]{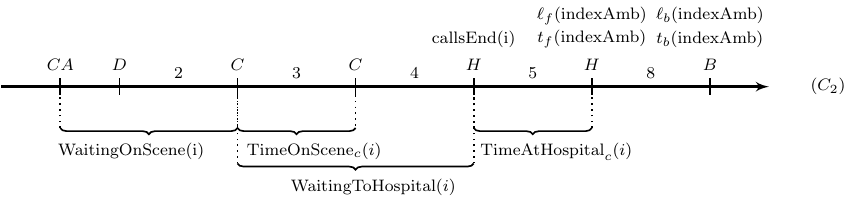}\\
\includegraphics[scale=1.1]{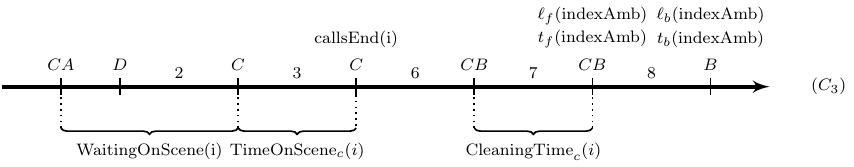}\\
\includegraphics[scale=1.18]{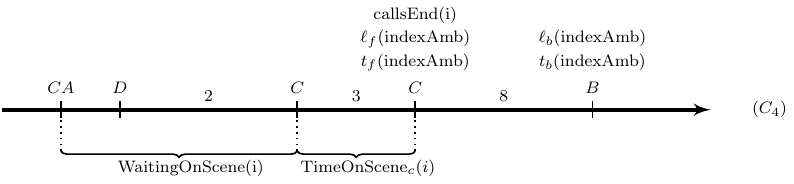}
\end{tabular}
\caption{Sequences of trips for the classes of service $C_{1}$, $C_{2}$, $C_{3}$, $C_{4}$.
The letters above the time axis represent locations for the corresponding time instants.
Specifically, CA: call arrival, D: departure of the ambulance to the emergency scene, C: emergency scene, H: hospital, CB: cleaning station, B: staging station.
Also, the trip type numbers 2, 3, 4, 5, 6, 7, and 8 are shown above the time axis for all classes of service. These trip type numbers correspond to those used in \cite{websiteambrouting24}. indexAmb is the index of the allocated ambulance.}
\label{figgroupcalls}
\end{figure}

The following parameters are associated with emergency of type $c$ indexed with~$i$:
\begin{itemize}
\item
the time instant $t_{c}(i)$ that the emergency call is received;
\item
the location $\ell_{c}(i)$ of the emergency, given by a pair (latitude, longitude);
\item
the type $\mbox{Type}_{c}(i)$ of emergency~$i$;
\item
the time $\mbox{CleaningTime}_{c}(i)$ needed to clean the ambulance at a cleaning station after serving emergency~$i$, if applicable, see Figure \ref{figgroupcalls};
\item
the time $\mbox{TimeOnScene}_{c}(i)$ spent by the ambulance on the scene of emergency~$i$, see Figure \ref{figgroupcalls};
\item
the time $\mbox{TimeAtHospital}_{c}(i)$ that an ambulance spends at the hospital after transporting the patient of emergency~$i$ to a hospital, if applicable, see Figure \ref{figgroupcalls}.
\end{itemize}
The following variables are associated with emergency of type $c$ indexed with~$i$:
\begin{itemize}
\item
the time $\mbox{waitingOnScene}(i)$ between the instant emergency call~$i$ is received and the instant that an ambulance arrives on the scene of the emergency, see Figure \ref{figgroupcalls};
\item
the time $\mbox{waitingToHospital}(i)$ between the instant that an ambulance arrives on the scene of emergency~$i$ and the instant that the ambulance arrives at a hospital with the
patient of emergency~$i$, if applicable, see Figure \ref{figgroupcalls};
\item the time callsEnd(i) when the ambulance finishes service, see Figure \ref{figgroupcalls}.
\end{itemize}

\if{
Each parking station $i$ has a location Station($i$) while cleaningStation($i$).location is the location of cleaning station $i$.
The location (given by a pair (latitude, longitude)) of the closest parking station to cleaning station $i$ is cleaningStation(i).station.

We will also denote by
\begin{enumerate}
\item
hospitals($i$).location: the location (given by a pair (latitude/longitude)) of hospital $i$,
\item
hospitals($i$).station: the station which is the closest to hospital $i$.
\item
hospitals($i$).cleaningStation: the index of the cleaning station which is the closest to the hospital $i$.
\end{enumerate}

We denote by NbAmbulances the number of ambulances, by NbStations the number of ambulance stations, and by NbCalls the number of calls over the planning horizon.
}\fi

For any two locations $A$ and $B$ and any time instant~$t_{0}$, let $\mbox{travelTime}(A, B, t_{0})$ denote the travel time of an ambulance that starts at $t_{0}$ from $A$ and travels to $B$, and let $\mbox{positionBetweenOriginDestination}(A, B, t_{0}, t)$ denote the position at instant~$t$ of an ambulance which starts a trip from $A$ to $B$ at time instant $t_{0}$.

\subsection{Performance metric}

The performance metric for an emergency of type~$c$ attended by an ambulance of type~$a$ with resulting response time~$t$ is specified by function
\begin{equation}
\label{costalloc0}
\mbox{\tt cost\_allocation\_ambulance}(a,c,t)
\end{equation}
As an example, consider the situation where this function has two components:
\begin{equation}
\label{costalloc}
\mbox{\tt cost\_allocation\_ambulance}(a,c,t) \ \ = \ \ \mbox{\tt penalization}(t,c) + M_{ac}
\end{equation}
where $\mbox{\tt penalization}(t,c)$ denotes the penalized response time for an emergency of type~$c$ waiting for~$t$ time units and $M_{ac}$ denotes the ``cost'' (in response time units) of assigning an ambulance of type~$a$ to an emergency of type~$c$. See Table~\ref{tbl:compatibility_matrix} for an example of such costs, used in the numerical examples of Section~\ref{sec:numsec}.
Several penalized response times have been proposed in the literature.
We provide two examples below.
\begin{example}
Linear penalized response times for every emergency type.
For this model,
\begin{equation}
\label{penfunc}
\mbox{\tt penalization}(t,c) \ \ = \ \ \theta_{c} t
\end{equation}
is a linear function of response time~$t$ that depends on the emergency type~$c$.
\end{example}
\begin{example}
Target response time.
This model has a threshold time $\alpha_{c} > 0$, such that it is desirable to have response times below $\alpha_{c}$ (for example, $\alpha_{c} = 8$ minutes for high priority emergencies and $\alpha_{c} = 15$ minutes otherwise).
Then
\[
\mbox{\tt penalization}(t,c) \ \ = \ \ \left\{
\begin{array}{lcl}
0 & \mbox{if } t \leq \alpha_{c},\\
\theta_{c} (t - \alpha_{c}) + \beta_{c} & \mbox{if } t > \alpha_{c},\\
\end{array}
\right.
\]
where $\beta_{c} \geq 0$ is a jump in the cost when the response time goes
above the threshold.
\end{example}

\section{Stochastic optimization model}
\label{smodel}

We consider a stochastic model of ambulance operations that includes random locations, types, and arrival times of emergencies, as well as random travel times of ambulances.
The stochastic model, and its simulation implementation, is described in~\cite{websiteambrouting24}.
In this section, we describe the rollout framework that we used for various ambulance dispatch policies.

\paragraph{The exogenous random variables}
The exogenous random variables include the types, locations, and arrival times of emergency calls, travel times, on-scene times, dwell times at hospitals, and cleaning times.
Notation $\xi$ is used for a sample path of exogenous random variables, and notation $\xi_{t_{0}}$ is used for a sample path of exogenous random variables after time~$t_{0}$.
The methods of this paper make provision for random travel times, on-scene times, dwell times at hospitals, and cleaning times, but the exposition describes only random types, locations, and arrival times of emergency calls.

\paragraph{The ambulance dispatch decisions}
Recall that two types of events trigger decisions: (i)~when an emergency call arrives and (ii) when an ambulance finishes service.
Let $t_{0}$ denote the time instant when one of these decisions has to be made, that is, $t_{0}$ is either the time instant when an emergency call arrives or the time instant when an ambulance finishes service.
Given the event data obtained at time~$t_{0}$, let $\mathcal{S}$ denote the set of feasible decisions.
Typically, $\mathcal{S}$ is of small cardinality, and can easily be enumerated, as explained next.

\paragraph{An emergency call arrives}
When a call arrives, either an available ambulance is dispatched to that emergency immediately, with possible transport of the patient(s) to a hospital thereafter, or the emergency is placed in the queue of emergencies.
Thus, the number of feasible decisions (the cardinality of $\mathcal{S}$) is one plus the number of combinations of available ambulances and hospitals appropriate for the emergency.

\paragraph{An ambulance finishes service}
When an ambulance finishes service, either the ambulance is dispatched to an emergency in queue, or the ambulance is dispatched to an ambulance station.
Thus the number of feasible decisions (the cardinality of $\mathcal{S}$) is the number of ambulance stations plus the number of combinations of emergencies in queue appropriate for the ambulance and hospitals appropriate for the emergency.

\if{
\paragraph{Building the two-stage problems}
We assume we have at hand a model for the emergency calls that allows us to generate scenarios of calls (such will be the case for the numerical experiments we consider in Section~\ref{sec:numsec} where we have calibrated a model for the process of calls based on historical data of emergency calls to an emergency medical service).
Next, let $x_{0}$ be the first stage decision taken at $t_{0}$ (recall that decision $x_{0} \in \mathcal{S}$ has just been described above and tells us  what to do with the call that has just arrived or with the ambulance that has just finished service).
Let also $\xi_{t_{0}}^{i}(x_{0}),i=1,\ldots,N$, be a set of scenarios of calls from $t_{0}$ (the instant the current decision needs to be taken) until the end of the horizon $T$, say a few hours.
Notice that scenario $\xi_{t_{0}}^{i}(x_{0})$ indeed depends on first stage decision $x_{0}$ because it not only has calls arriving at instants $t > t_{0}$ but also the queue of calls at $t_{0}$ after taking
decision $x_{0}$.
More precisely, on top of scenarios of future calls (arriving after $t_{0}$) $\xi_{t_{0}}^{i}(x_{0})$ contains the following additional calls:
\begin{itemize}
\item
if a call arrives at $t_{0}$ and decision $x_{0}$ is to send an available ambulance to that call then $\xi_{t_{0}}^{i}(x_{0})$ contains the calls that were still in queue just before the arrival of the call at $t_{0}$;
\item
if a call arrives at $t_{0}$ and decision $x_{0}$ is to put that call in queue, then $\xi_{t_{0}}^{i}(x_{0})$
contains the call that arrived at $t_{0}$ together with the calls that were still in queue just before the arrival of the call at $t_{0}$.
Additionally, since the call is put in queue to use an ambulance that is not available (otherwise we would have sent immediately an available ambulance), we add to the call that arrived at $t_{0}$ the information/constraint that it needs to be attended by an ambulance that is busy at $t_{0}$.
\item
When an ambulance finishes service and decision $x_{0}$ is to send that ambulance to a call in queue,
say call $c$, then $\xi_{t_{0}}^{i}(x_{0})$ contains the queue of calls just before the ambulance became available at $t_{0}$ minus call $c$.
\item
When an ambulance finishes service and decision $x_{0}$ is to send that ambulance to a station, then $\xi_{t_{0}}^{i}(x_{0})$ contains the queue of calls just before the ambulance became available at $t_{0}$ with the additional constraint that these calls cannot be attended by the ambulance that just finished service at $t_{0}$ (otherwise we would have sent that ambulance to one of the calls in queue at $t_{0}$).
\end{itemize}

The $i$th scenario $\xi_{t_{0}}^{i}(x_{0})$ therefore provides a set of calls given by their locations, types, instants, and possibly, as we have seen, constraints on a subset of ambulances that can attend these calls.
}\fi

\paragraph{The state of the process.}
The state of the process at a point $t_{0}$ in time includes information about the location of ambulances, the remaining tasks in the current service of each ambulance, the emergencies in queue, and the current time (to accommodate time-dependent parameters such as time-varying demand and time-varying travel times).
The state contains sufficient information so that future states can be determined as a function of the current state, ambulance dispatch decisions, and exogenous random variables.
At any time $t$, the state includes the following quantities (see Figure \ref{figgroupcalls} for  a representation of the components of the state vector):
\begin{enumerate}
\item[(1)]
For each ambulance~$j$, a location $\ell_{f}(j)$ and a time $t_{f}(j)$, with meanings depending on whether (1.1) ambulance~$j$ is in service at time~$t$ or (1.2) ambulance~$j$ is available at time~$t$.
\begin{enumerate}
\item[(1.1)]
If ambulance~$j$ is in service at time~$t$, then $t_{f}(j)$ denotes the next time that the ambulance will be available for dispatch and $\ell_{f}(j)$ denotes its location at that time.
Therefore, if ambulance~$j$ is in service at time~$t$, then $t_{f}(j) > t$.
The following four situations can happen:
\begin{enumerate}
\item[(a)]
The ambulance does not need to transport a patient to a hospital and does not need to go to a cleaning station.
Then $\ell_{f}(j)$ is the location of the emergency where the ambulance is serving at time~$t$, and $t_{f}(j)$ is the time when the ambulance will complete service at the scene of the emergency.
\item[(b)]
The ambulance does not need to transport a patient to a hospital but needs to go to a cleaning station.
Then $\ell_{f}(j)$ is the location of the cleaning station where the ambulance will go after serving the current emergency, and $t_{f}(j)$ is the time when cleaning will be completed.
\item[(c)]
The ambulance has to transport a patient to a hospital and does not need to go to a cleaning station.
Then $\ell_{f}(j)$ is the location of the hospital, and $t_{f}(j)$ is the time when the ambulance can leave the hospital.
\item[(d)]
The ambulance has to transport a patient to a hospital and needs to go to a cleaning station.
Then $\ell_{f}(j)$ is the location of the cleaning station where the ambulance will go after leaving the hospital, and $t_{f}(j)$ is the time when cleaning will be completed.
\end{enumerate}
\item[(1.2)]
If the ambulance is available at time~$t$, then $t_{f}(j)$ is the previous (past) time instant when the ambulance became available, i.e., the time when it completed its last service (thus, in this case, $t_{f}(j) \le t$), and $\ell_{f}(j)$ is the corresponding previous location of the ambulance when it became available (either the location of an emergency, a hospital, or a cleaning station).
\end{enumerate}
\item[(2)]
For each ambulance~$j$, a location $\ell_{b}(j)$ and a time $t_{b}(j)$ with the following meanings.
\begin{enumerate}
\item[(2.1)]
If ambulance~$j$ is on its way to an ambulance station, then $\ell_{b}(j)$ is the location of that station, and $t_{b}(j)$ is the time when the ambulance will arrive at that station (thus, in this case, $t_{b}(j) > t$).
\item[(2.2)]
If ambulance~$j$ is not on its way to an ambulance station, then $\ell_{b}(j)$ is the location of the previous station that the ambulance was on its way to, and $t_{b}(j)$ is the previous time that it was scheduled to arrive at this station (thus, in this case, $t_{b}(j) \leq t$).
\end{enumerate}
\end{enumerate}
The state also includes information about emergencies in queue.

\if{
We will also store the discretized trajectories of ambulances, i.e, for every ride of every ambulance the pair (origin, destination), the time the origin was left, the time the destination was reached, and the ride type.
With this information at hand, we will be able to know at every time instant the position of every ambulance, using the function positionBetweenOriginDestination described before, which, given inputs (A, B, $t_{0}$, $t$) returns the position at $t$ of an ambulance that leaves $A$ at $t_{0}$ to go to $B$ (we therefore assume that the travel time of a ridfe is known as long as we know for this ride the origin, the destination as well as  the departure time).

More precisely, we will use the following arrays:
\begin{enumerate}
\item[1)]
ambulancesTrips: ride $i$ (over the planning period) of ambulance $j$ starts at location ambulancesTrip$(j)(i)$ and ends at location ambulanceTrip$(j)(i+1)$ (both origin and destination are given by a pair (latitude/longitude)).
\item[2)]
ambulancesTimes: ride $i$ (over the planning period) of ambulance $j$ starts at time ambulancesTimes$(j)(i)$ and ends at time ambulanceTimes$(j)(i+1)$.
\item[3)]
tripType$(j)(i)$ is the type of $i$th ride for ambulance~$j$.
This ride type takes 8 values:
\begin{itemize}
\item 1: the ambulance stays at the station.
\item 2: the ambulance is on its way to the scene of a call (either coming from a fixed station or a mobile station).
\item 3: the ambulance is on the scene of the call.
\item 4: the ambulance is going to a hospital.
\item 5: the ambulance is at hospital, waiting to be freed.
\item 6: the ambulance is going to a cleaning station.
\item 7: the ambulance is being cleaned at a station.
\item 8: the ambulance is going to a station (not for a cleaning task).
\end{itemize}
Observe that the trip type numbers are in chronological order within the cycle of a call, which starts with a ride of type 2 for an ambulance which was on a ride of type 1 or 8.
For instance, if the patient of the call is transported to a hospital and ambulance cleaning is needed then the ride types for the corresponding call are, in chronological order, 2, 3, 4, 5, 6, 7, and 8.
\end{enumerate}

These arrays are initialized given the ongoing rides of the ambulances when the planning period starts.

}\fi

\paragraph{The rollout approach}
Consider any decision time~$t_{0}$ (associated with the arrival of an emergency or the completion of service by an ambulance), and any feasible decision~$x_{0}$ at time~$t_{0}$.
Choose any time horizon, and consider an i.i.d.\ set $\left\{\xi^{i}_{t_{0}}\right\}_{i=1}^{N}$ of scenarios, from time~$t_{0}$ until the end of the specified time horizon.
A scenario $\xi^{i}_{t_{0}}$ specifies the values of all random variables that may be relevant for the performance of the system over the time horizon, including the arrival times of emergency calls and the type of each emergency, the travel times along different paths in the network (irrespective whether an ambulance ends up traveling on the path or not), the on-site emergency service time, the time to hand over the patient to a hospital if relevant for the emergency, and the time to clean the ambulance after each emergency.
Consider any policy~$H$ for making ambulance dispatch decisions, both when an emergency call arrives and when an ambulance finishes service.
Thus, for any $x_{0}$ and $\xi^{i}_{t_{0}}$, policy~$H$ specifies all decisions after time~$t_{0}$ until the end of the specified time horizon.
Given the state of the process and the decision~$x_{0}$ at time~$t_{0}$, let $f(x_{0})$ denote the immediate cost of decision~$x_{0}$ (for example, the cost given by cost function~\eqref{costalloc} associated with an emergency if decision~$x_{0}$ dispatches an ambulance to the emergency, and the cost of moving the ambulance to the ambulance station if decision~$x_{0}$ dispatches an ambulance to an ambulance station), and let $Q_{H}(x_{0},\xi^{i}_{t_{0}})$ denote the total cost resulting from decision~$x_{0}$ and the use of rollout policy~$H$ after time~$t_{0}$ for scenario $\xi^{i}_{t_{0}}$ over the specified time horizon.
Cost $Q_{H}(x_{0},\xi^{i}_{t_{0}})$ may include an end-of-horizon cost that penalizes undesirable states at the end of the time horizon.
It is desirable to use a rollout policy~$H$ such that $Q_{H}(x_{0},\xi^{i}_{t_{0}})$ can be computed efficiently for any $x_{0}$ and $\xi^{i}_{t_{0}}$.
(This notation does not show the state at time~$t_{0}$ or the time horizon, because these quantities are the same for all dispatch decisions~$x_{0}$ under consideration at time~$t_{0}$.)

Under a rollout approach, given rollout policy~$H$, at each decision time~$t_{0}$ a dispatch decision~$x_{0}$ is chosen that solves
\begin{equation}
\label{secstagepb}
\min_{x_{0} \in \mathcal{S}} \; \left\{f(x_{0}) + \frac{1}{N} \sum_{i=1}^N Q_{H}(x_{0},\xi^{i}_{t_{0}})\right\}.
\end{equation}
Since $\mathcal{S}$ has small cardinality, one can compute the objective function in~\eqref{secstagepb} for every feasible decision $x_{0} \in \mathcal{S}$ and choose the best one.
Note that the resulting policy is not the same as rollout policy~$H$ --- the decision $x_{0}$ that solves~\eqref{secstagepb} may be different from the decision prescribed by rollout policy~$H$ in the same state.
Also note that the objective function in~\eqref{secstagepb} uses an i.i.d.\ sample $\left\{\xi^{i}_{t_{0}}\right\}_{i=1}^{N}$ of scenarios, so that by the law of large numbers $\frac{1}{N} \sum_{i=1}^N Q_{H}(x_{0},\xi_{t_{0}}^{i})$ approximates $\mathbb{E}\left[Q_{H}(x_{0},\xi_{t_{0}})\right]$, and that none of the sampled scenarios $\xi^{i}_{t_{0}}$ have to match the sample path $\xi_{t_{0}}$ that ends up occurring (for example, the sample path $\xi_{t_{0}}$ that ends up occurring in a simulation used to compare various policies).
The resulting policy is feasible and nonanticipative, even if rollout policy~$H$ is infeasible or anticipative.

In Sections~\ref{bestmyopic}, \ref{nonmyopic}, \ref{sec:prior}, and~\ref{sec:newbase}, we propose heuristics that can be used as rollout policies~$H$ to quickly compute the rollout cost $Q_{H}(x_{0},\xi^{i}_{t_{0}})$ for any decision time~$t_{0}$, state at time~$t_{0}$, decision~$x_{0}$ at time~$t_{0}$, and scenario $\xi^{i}_{t_{0}}$.

\section{Best Myopic heuristic}
\label{bestmyopic}

The \emph{Best Myopic} (BM) heuristic uses a simple improvement of the \emph{closest-available-ambulance} heuristic.
This heuristic is ``best'' in the sense that the immediate allocation costs given by~\eqref{costalloc0} are minimized.
The heuristic is myopic in the sense that it uses no information (either deterministic or probabilistic) about future events when decisions are made.
The heuristic requires as input estimates of the service times, including the associated travel times, for an emergency.

The BM heuristic allocates ambulances to emergencies sequentially.
When the $i$th emergency call arrives, ambulances have already been allocated to all previous emergencies $1, 2, \ldots, i-1$.
Therefore, the trajectory of each ambulance involved in serving emergencies $1, 2, \ldots, i-1$ is known (or estimated, as mentioned above) until completion of the service for these emergencies.
In particular, if ambulance~$j$ is at an ambulance station at the time $t_{c}(i)$ when the $i$th emergency call arrives, then we know the station $\ell_{b}(j)$ at which the ambulance is, and if ambulance~$j$ is on its way to an ambulance station at time $t_{c}(i)$, then we know where the ambulance currently is and the time $t_{b}(j)$ when it will arrive at the station.
Also, if ambulance~$j$ is in service at time $t_{c}(i)$, then we know the location $\ell_{f}(j)$ where and the time $t_{f}(j)$ when it will become available.
This allows us to calculate the response time of emergency~$i$ for every possible ambulance~$j$ that can be assigned to that emergency, as well as the corresponding allocation cost
given by~\eqref{costalloc0}.
The heuristic is then based on the following four main ingredients to choose the ambulance for emergency~$i$, and if necessary the hospital and cleaning station:
\begin{itemize}
\item[(A)]
Given the current status of all ambulances, the allocation cost~\eqref{costalloc0} for every possible ambulance~$j$ that can be assigned to emergency~$i$ is computed.
Next, we explain the computation of the response time for a given ambulance $j$.
We consider three cases: (A1) $t_{b}(j) \leq t_{c}(i)$, (A2) $t_{f}(j) \leq t_{c}(i) < t_{b}(j)$, and (A3) $t_{c}(i) < t_{f}(j)$.
\begin{enumerate}
\item[(A1)]
$t_{b}(j) \leq t_{c}(i)$, i.e., ambulance~$j$ is at a station when call~$i$ arrives:
In this case, the response time is the time for the ambulance to go from station $\ell_{b}(j)$ (since it is immediately available) to the emergency location $\ell_{c}(i)$ departing at time $t_{c}(i)$, so it is given by $\mbox{travelTime}(\ell_{b}(j),\ell_{c}(i),t_{c}(i))$.
\item[(A2)]
$t_{f}(j) \leq t_{c}(i) < t_{b}(j)$, i.e., the ambulance is available but is traveling to a station:
We first compute its position $P$ on its current trip, given by 
\[
P \ \ = \ \ \mbox{positionBetweenOriginDestination}(\ell_{f}(j), \ell_{b}(j), t_{f}(j), t_{c}(i)).
\]
The response time is then the travel time for ambulance~$j$ from $P$ to $\ell_{c}(i)$ departing at time $t_{c}(i)$ (since it is immediately available), so it is given by $$\mbox{travelTime}(P, \ell_{c}(i), t_{c}(i)).$$
\item[(A3)]
$t_{c}(i) < t_{f}(j)$, i.e., ambulance~$j$ is in service when call~$i$ arrives:
First we compute the time until the ambulance finishes its service, which is given by $t_{f}(j) - t_{c}(i)$.
Then we compute the time for the ambulance to travel from $\ell_{f}(j)$ (where it will be when it becomes available) to the emergency location $\ell_{c}(i)$ departing at time $t_{f}(j)$.
This time is given by $\mbox{travelTime}(\ell_{f}(j),\ell_{c}(i),t_{f}(j))$.
Therefore, if ambulance $j$ is used to serve call~$i$, the the response time is given by $$t_{f}(j) - t_{c}(i) + \mbox{travelTime}(\ell_{f}(j),\ell_{c}(i),t_{f}(j)).$$
\end{enumerate}
\item[(B)]
The set of ambulances with the least allocation cost is determined.
(All ambulances that can be dispatched from the same ambulance station will achieve the same response time for a given emergency, so there may be multiple ambulances with the same allocation cost.)
\item[(C)]
If multiple ambulances achieve the least allocation cost, then choose among them one of the least advanced ambulances.
For instance, if a BLS ambulance and an ALS ambulance both achieve the least allocation cost, then we choose a BLS ambulance (to save the ALS ambulance for emergencies that require more advanced capabilities).
\end{itemize}

\begin{example}
\label{exrj}
Suppose that there are three types of emergencies: high, intermediate, and low priority emergencies.
Also, suppose that there are three types of ambulances: advanced ambulances that are good for all types of emergencies, intermediate ambulances that are good for intermediate and low
priority emergencies but that are inferior to advanced ambulances for high priority emergencies, and basic ambulances that are good for low priority emergencies but that are inferior for intermediate priority emergencies and especially for high priority emergencies.
In this situation,
\begin{itemize}
\item[(i)]
if the smallest allocation cost is achieved by at least one basic ambulance, then a basic ambulance is chosen among these;
\item[(ii)]
if the smallest allocation cost is achieved by no basic ambulance but by at least one intermediate ambulance, then an intermediate ambulance is chosen among these;
\item[(iii)]
if the smallest allocation cost is only achieved by advanced ambulances, then one of these ambulances is chosen.
\end{itemize}
\end{example}

Note that, under the BM heuristic, emergencies may wait for its allocated ambulance to complete a previous service, but there is no queue of emergencies that have not yet been allocated an ambulance.
Observe that the BM heuristic is not the same as the \emph{closest available ambulance heuristic}, because no ambulance may be available when the call arrives, or the ambulance 
with the least allocation cost may not be available when the call arrives, and may not even be the next available ambulance when the call arrives.
Moreover, the closest available ambulance heuristic does not take emergency types and ambulance types into account.

\section{A NonMyopic nonanticipative heuristic}
\label{nonmyopic}

The \emph{NonMyopic (NM) nonanticipative}  heuristic takes into account calls that may arrive in a near future, but it is still a nonanticipative policy, meaning that an ambulance cannot be dispatched to an emergency call that has not arrived yet.

The steps of the NM heuristic are as follows. 
We consider a call index $i$ initialized to 1.\\
\textbf{Step 1.}
If an ambulance has already been allocated to emergency~$i$, then
either $i$ is the last call and the heuristic ends
or we increment the call index: $i \leftarrow i+1$ and go back to Step~1.
If no ambulance has been allocated to emergency $i$, we go to Step 2.\\
\textbf{Step 2.}
Given the current status of all ambulances, the allocation cost~\eqref{costalloc0} for every possible ambulance~$j$ that can be assigned to emergency~$i$ is computed, as for the BM heuristic. 
In this computation,  the allocation cost for busy ambulances is done knowing the service times (on the scene of the call and waiting at hospital).\\
\textbf{Step 3.}
The set~$\mathcal{S}(i)$ of ambulances with the smallest allocation cost 
for call $i$
is determined.
If any ambulance $j \in \mathcal{S}(i)$ is available, then such an ambulance~$j$ is dispatched to emergency~$i$, and the heuristic goes to Step~1.
Otherwise, if all ambulances in $\mathcal{S}(i)$ are busy, then
take  $j \in \mathcal{S}(i)$ a least advanced ambulance from the
set $\mathcal{S}(i)$
and go to Step~4. \\
\textbf{Step 4.}
To choose an ambulance  to allocate to emergency~$i$, the heuristic takes into account alternative allocations of ambulance~$j$ to future emergencies. We consider all calls $i'$ which arrive not later
than  time $t_{f}(j)$ when ambulance~$j$ 
will next finish service, i.e., calls $i'$ such that
$t_c(i') \leq t_f(j)$ (note that call $i$ is one of them) and such that no allocation has been made for these calls. For these calls, given previous allocation of ambulances,
we compute the allocation costs of all ambulances
to all these calls $i'$ which allows us to compute
for each call $i'$ with $t_c(i') \leq t_f(j)$ the minimal
allocation cost denoted by minAlloc($i'$).
We denote by $\mathcal{T}(j)$ the set of calls
$i'$ satisfying $t_c(i') \leq t_f(j)$
and such that the minimal allocation cost
minAlloc($i'$) is obtained with ambulance $j$, meaning that
$j \in \mathcal{S}(i')$.
Let now $i^* \in \arg\max\{\mbox{minAlloc}(i') \, : \, i' \in \mathcal{T}(j)\}$. We allocate ambulance j to emergency 
$i^*$. If $i^*=i$ we go to Step 1. Otherwise, we go to Step 2.


\section{Greedy heuristics}
\label{sec:prior}

The heuristics described in this section are variants of the BM heuristic that make provision for a queue of emergencies that have not yet been allocated to an ambulance.
More specifically, allocation costs are computed in the same way as for the BM heuristic.
However, for the heuristics considered in this section, an ambulance has to be available when it is allocated to an emergency, and therefore these heuristics have to make provision for a queue of emergencies that have not yet been allocated to an ambulance.
These heuristics use the state of the BM heuristic augmented with the queue of emergencies.
We describe two heuristics: Greedy Heuristic with Priorities 1 (GHP1) and Greedy Heuristic with Priorities 2 (GHP2).
We will use the following
notation.
Let $\mbox{queue}(k)$ denote the index of the $k$th emergency in the queue (ordered as described below), and let $\mbox{queueSize}$ denote the size of the current queue.
Finally, we denote by currentTime the current time when a decision has to be made.

\subsection{GHP1}
\label{sec:firstgreedy}

The steps of GHP1 are described below. In these computations, one of our three reassignement heuristics for the choice of a station (described in Section \ref{sec:newbase}) is assumed to be used.\\
\textbf{Step 0.} We compute the next event which is
"a call arrives" or
"an ambulance finishes service".
Go to Step 1.\\
\textbf{Step 1.} If the event is 
"a call arrives" then add the index of the call  to the queue of calls and increment 
by one the size $\mbox{queueSize}$ of the queue of calls: $\mbox{queueSize} \leftarrow \mbox{queueSize} + 1$. Go to Step 2. \\
\textbf{Step 2.}
Sort the emergencies in the queue in decreasing order of the penalized time ${\tt penalization}(\mbox{currentTime} - t_{c}(i), \mbox{Type}_{c}(i))$ that has elapsed since the call arrived, where {\tt{penalization}} is the function in \eqref{costalloc}. \\
\textbf{Step 3.}
In this step, available ambulances (if any) are allocated to some of the emergencies in the queue. \\
For each $k = 1, \ldots, \mbox{queueSize}$, do the following: given the current status of all ambulances, the allocation cost~\eqref{costalloc} for every possible ambulance~$j$ that can be assigned to emergency~$\mbox{queue}(k)$ is computed, as for the BM heuristic.
Let $\mathcal{S}(\mbox{queue}(k))$ denote the set of ambulances with the least allocation cost for emergency~$\mbox{queue}(k)$.
If none of the ambulances in $\mathcal{S}(\mbox{queue}(k))$ is available, then increment $k \leftarrow k + 1$, and consider the next emergency in queue.
Otherwise, among the ambulances in $\mathcal{S}(\mbox{queue}(k))$ that are available, allocate one of the least advanced ambulances to emergency~$\mbox{queue}(k)$, remove emergency~$\mbox{queue}(k)$ from the queue, decrement $\mbox{queueSize} \leftarrow \mbox{queueSize} - 1$, and consider the next emergency in queue.\\
Go to Step 4.\\
{\textbf{Step 4.}} We update the next event: either "a call arrives"
or "an ambulance finishes service". Go to Step 1.

\if{
We compute the cost, given by~\eqref{costalloc},  of allocating every ambulance to the call with index queue$(k)$ considering the current scheduled rides for all ambulances.
Among the ambulances that are the best for this call (by best, we mean the ambulances providing the smallest value of the allocation cost given by \eqref{costalloc}), we check if some of them are available at instant currentTime.
If this is the case, then we select among these "best" available ambulances one of the least advanced
(for instance, if there are 2 "best" available ambulances, one advanced (ALS)
and one basic (BLS), we choose the basic one).
We denote again by indexAmb$(k)$ the index of the 
corresponding
best ambulance for call with index queue($k$).
If this is not the case, i.e., if 
no "best"
ambulance is available, then we put the corresponding call index in the queue queueAux of the calls for which an ambulance will be sent at a later time.
This means that queue(k) is added to queueAux.
Otherwise, we allocate to call with index queue(k) the ambulance with index indexAmb($k$) and we update the state vector and corresponding rides for that ambulance as was done in the BM heuristic.
We increase callsAttended by 1.
\par {\textbf{Step 3.b.}}
If $k = \mbox{length(queue)}$, then we go to {\textbf{Step 4}} otherwise we increase $k$ by 1 and go to {\textbf{Step 3.a}}.
\par {\textbf{Step 4: update the next value of currentTime, indexCall, and eventCall.}}
We update the new queue of calls to queueAux, and we now need to compute the type of the next event (either a new call arrives or an ambulance finishes service) updating the value of the boolean variable eventCall which, as we recall, will be 1 if the next event is a call and 0 if the next event is {\em{the ambulance finishes service.}}
We also compute the next value of currentTime which will be the time corresponding to this next event and the value of the index of the next call to treat if the next event is the arrival of a call.
For that, we collect in an array futureArrivalTimes the list of future times $t_{f}(k)$ when an ambulance $k$ will become available again.
We have two cases: a) either this list is empty or b) it is not empty.
\par In case a), if indexCall is strictly less than the total number of calls to treat we still have future calls, and we set eventCall=1, we increase indexCall by 1, and set currentTime to
$t_{c}(\mbox{indexCall})$.
\par In case b), we compute the minimum value minTime in futureArrivalTimes.
If indexCall is equal to the total number of calls then ambulances have already been allocated to all calls, and we set eventCall=0, currentTime=minTime.
Otherwise, we have two possibilities: b1) $t_{c}(\mbox{indexCall}+1) \leq \mbox{minTime}$ and b2) $t_{c}(\mbox{indexCall}+1) > \mbox{minTime}$.
In case b1) the next event is {\em{a call arrives}} and we set eventCall=1, we increase indexCall by 1, and we set currentTime to the instant $t_{c}(\mbox{indexCall})$ of the next call (recall that in this latter notation, indexCall has already been increased by 1).
In case b2), the next event is {\em{an ambulance finishes service}} and therefore we set eventCall = 0 and the time instant of the next event will be currentTime = minTime.
\par This closes the treatment of the calls in the queue at time instant currentTime.
We then repeat the above computations as long as there are calls for which an ambulance has not been allocated to that call.
}\fi

\subsection{GHP2}
\label{sec:secondgreedy}

In this section, we describe a variant GHP2 of GHP1.
Whereas GHP1 orders the emergencies~$i$ in the queue in decreasing order of the penalized time ${\tt penalization}(\mbox{currentTime} - t_{c}(i), \mbox{Type}_{c}(i))$ that has elapsed since the call arrived, GHP2 orders the emergencies in the queue in decreasing order of the minimum allocation cost~\eqref{costalloc0} over all ambulances, given the current status of all ambulances.
\if{
Each time a decision has to be taken, in GHP2, we compute for every call $i$ in queue for which no ambulance has been sent yet and for each ambulance $j$, the immediate cost of allocating call $i$ to ambulance $j$ (immediately after ambulance $j$ finishes rides already forecast for that ambulance at the current time currentTime) where the immediate allocation cost is given by
\eqref{costalloc}.
We then consider a call with the worst (largest) value of the  allocation cost.
We then have two possibilities:
(a) a {\em{best}} ambulance (we call a best ambulance an ambulance achieving the smallest value of the allocation cost, i.e., minimizing the immediate allocation cost given previous allocations) for that call is available or
(b) none of the best ambulances for that call are available.
In case (b), we put the call in a new queue, denoted by queueAux as in GHP1, of calls to which an ambulance will be allocated after the current time currentTime.
In case (a), we allocate a best ambulance to this call updating state variables and ambulance rides correspondingly.
Moreover, if several ambulances achieve the same allocation cost, we choose among them an ambulance as basic as possible.
We now provide the detailed computations for that heuristic.

Same as for GHP1, we denote by queue the variable storing the indexes of the calls for which no ambulance has been allocated so far and by queueSize the size of this queue.
We also use variables eventCall, currentTime, indexCall, and callsAttended which have the same meaning as in GHP1.
}\fi

The steps of GHP2 are described below. In these computations, one of our three reassignement heuristics for the choice of a station (described in Section \ref{sec:newbase}) is assumed to be used.\\
\textbf{Step 0.} We compute the next event which is
"a call arrives" or
"an ambulance finishes service".\\
\textbf{Step 1.} If the event is "an ambulance finishes service" and the queue is empty, call the reassignment heuristic and
go to Step 4; otherwise go to Step 2.
If the event is 
"a call arrives" then add the index of the call  to the queue of calls and increment 
by one the size $\mbox{queueSize}$ of the queue of calls: $\mbox{queueSize} \leftarrow \mbox{queueSize} + 1$. Go to Step 2. \\
\textbf{Step 2.}
Given the current status of all ambulances, the allocation cost~\eqref{costalloc0} for every emergency~$i$ in queue and every possible ambulance~$j$ that can be assigned to emergency~$i$ is computed, as for the BM heuristic.
For each emergency~$i$ in queue, let $\mbox{minAlloc}(i)$ denote the least allocation cost over all ambulances for emergency~$i$, and let $\mathcal{S}(i)$ denote the set of ambulances with the least allocation cost for emergency~$i$.
Sort the emergencies~$i$ in the queue in decreasing order of $\mbox{minAlloc}(i)$. \\
\textbf{Step 3.}
For each $k = 1, \ldots, \mbox{queueSize}$, do the following:
If none of the ambulances in $\mathcal{S}(\mbox{queue}(k))$ are available, then increment $k \leftarrow k + 1$, and consider the next emergency in queue.
Otherwise, among the ambulances in $\mathcal{S}(\mbox{queue}(k))$ that are available, allocate one of the least advanced ambulances to emergency~$\mbox{queue}(k)$, update the status of the allocated ambulance, remove emergency~$\mbox{queue}(k)$ from the queue, decrement $\mbox{queueSize} \leftarrow \mbox{queueSize} - 1$, and go to Step~2.\\
\textbf{Step 4.} We compute the next event:
"a call arrives" or "an ambulance finishes service" and go to Step 1.

\ignore{
then we put this call with index indexCall in queue queueAux, and we go to {\textbf{Step 3.4}}. Otherwise, let indexAmb be the index of the first ambulance available in 
$$
\mbox{indexAmbs(remainingIndexes(currentCall))}
$$
(this corresponds to an ambulance as basic as possible,
among the best ambulances for that call).
We allocate the ambulance with index indexAmb to call with index indexCall as in the BM heuristic, updating the state vector and the rides of ambulance with index indexAmb.
\par {\textbf{Step 3.3.}} If ambulance with index indexAmb was the best (i.e., provided the smallest allocation cost) for some of the calls in the queue not already treated (that is to say calls with index queue(remainingIndexes$(k)$) for some $k$), we need to check, now that this ambulance has just been scheduled for another ride, if it is still among the best for these calls updating entries of arrays travelTimes, indexAmbs, and minAlloc as in Algorithm~\ref{algste33} below.
\par {\textbf{Step 3.4.}}
If (nbTreateadCalls+1) $<$ totalInQueue, then we suppress line currentCall of minAlloc, and remainingIndexes.
Increase nbTreatedCalls by 1.
If nbTreatedCalls is equal to length(queue) we go to {\textbf{Step 4}} otherwise we go to {\textbf{Step 3.1}}.
\par {\textbf{Step 4.}} Same as GHP1.
}

\if{
\begin{algorithm}
\begin{algorithmic}[1]
{\tt{
\FOR{$k \in \{1,\ldots,\mbox{totalInQueue}-\mbox{nbTreateadCalls}\}$}

\IF{k$\neq$currentCall and indexAmb $\in$ indexAmbs(remainingIndexes(k))}
\STATE //Compute the time needed for ambulance with index indexAmb to
\STATE //become available again:
\STATE timeToFree=$t_{f}$(indexAmb)-currentTime.
\STATE //New time for ambulance with index indexAmb to arrive on the scene
\STATE //of call with index queue(remainingIndexes(k)):
\STATE 
$$
\begin{array}{lcl}
\mbox{timeFromFreeToCall}&=&
\mbox{travelTime}(\ell_{f}(\mbox{indexAmb}),\\
&&\hspace*{2.3cm}\ell_{c}\mbox{(queue(remainingIndexes(k)))},\\
&&\hspace*{2.3cm}t_{f}({\mbox{indexAmb}}))
\end{array}
$$
\STATE //Update travelTimes:
\STATE 
\begin{tabular}{l}
travelTimes(remainingIndexes(k),indexAmb)=\\
timeToFree+timeFromFreeToCall
\end{tabular}
\STATE //Update immediate cost allocations:
\STATE costAllocations(remainingIndexes(k),indexAmb)=\\ 
\hspace*{2cm}cost\_allocation(Type$_{a}$(indexAmb),  \\
\hspace*{2cm}Type$_{c}$(queue(remainingIndexes(k))), \hspace*{2cm}travelTimes(remainingIndexes(k),indexAmb)\\ 
\hspace*{2cm}+currentTime-$t_{c}$(queue(remainingIndexes(k))))
\STATE //Update the best ambulance for that call
\STATE indexAmbs(remainingIndexes(k))=index
\STATE minAlloc(k) = $\infty$
\STATE indexAmbs(remainingIndexes(k)) = $\emptyset$
\FOR{$j \in \{1,...,\text{nbAmbulances}\}$}
    \IF{costAllocations(remainingIndexes(k),j) $<$ minAlloc(k)}
        \STATE \hspace*{-0.2cm}\mbox{indexAmbs(remainingIndexes(k))} = $\{j\}$
        \STATE \hspace*{-0.2cm}\mbox{minAlloc(k) = costAllocations(remainingIndexes(k),j)}
    \ELSIF{costAllocations(remainingIndexes(k),j) = minAlloc(k)}
        \STATE \hspace*{-0.2cm}indexAmbs(remainingIndexes(k))=indexAmbs(remainingIndexes(k))$\cup$\{j\}
    \ENDIF
\ENDFOR
\ENDIF
\ENDFOR
}}
\end{algorithmic}
\caption{Travel times, indexAmbs, and minAlloc updates.}
\label{algste33}
\end{algorithm}

}\fi

\section{Comparison of the heuristics on some examples}
\label{sec:examples}

To better understand the proposed heuristics and to show that
none is always better than the others, we apply these heuristics
on a few simple examples.

In these examples, calls
can have three types
(calls of high priority, of intermediate
priority, and of low priority).
Correspondingly there are
three types of ambulances:
advanced that can attend
all types of calls, 
intermediate that can
attend calls of intermediate
and low priorities, and basic
ambulances that can only attend
calls of low priorities.
The allocation cost of an ambulance
of type a to a call of type c
is given by \eqref{costalloc}
with function
{\tt{penalization}} of form
\eqref{penfunc} and $M_{a c}$ being
0 if type a ambulance can take type c call
and 0 otherwise. 
For function 
{\tt{penalization}}, we take
$\theta_c=4$ for high priority calls,
$\theta_c=2$ for intermediate priority calls,
and
$\theta_c=1$ for low priority calls.
The steps in these examples
correspond to the steps given in the general description of the heuristics.

The following example
provides an instance
where the
penalized response
times are larger with
the BM and
Closest Available (CA) ambulance
heuristics 
than with the NonMyopic
and  greedy heuristics.

\begin{example}\label{ex1toy} 
Consider an instance
in the plane parameterized
by some integer $t \geq 1$
$t$\footnote{We can easily
modify this example replacing
the plane by the surface of earth}:
there is
one hospital at $(5,5)$,
a call of low priority at
time 1+$20t$
at $(5,4.5)$,
a call of low priority  at
time 6.1+$20t$
at $(4.8,5.2)$,
a call of intermediate priority  at
time 6.3+$20t$
at $(5,5.2)$, and
a call of high priority  at
time 6.5+$20t$
at $(5.2,5.2)$. The time spent
on the scene of the calls
is 0.1, all calls will be sent
to the hospital and the ambulances
wait at hospital 0.4 for the
calls arriving at 1+20$t$
and 0.1 for the remaining calls.
No ambulance cleaning is necessary
and there are four ambulance stations
at $(0,0)$, $(10,0)$,
$(10,10)$, and $(0,10)$.
After service, ambulances
return to the closest station.
There are four ambulances,
all advanced and at $t=0$
there is one ambulance at each
station, say ambulance 1 at
$(0,0)$, ambulance 2 at $(10,0)$,
ambulance 3 at $(0,10)$, and 
ambulance 4 at $(10,10)$. Time is penalized
linearly with a factor
1 for low priority calls,
2 for intermediate priority
calls, and 4 for high priority calls.
All ambulances have speed 1.
Below, all allocation costs are computed at any time given the state of the system.
The computations are first done with
$t=1$ and then a comment is made
when we have sequences of four calls
arriving at $1+20t$,  $6.1+20t$,
$6.3+20t$, and
$6.5+20t$
for $t=1,2,3,\ldots$
\par {\textbf{Best Myopic.}} For BM, when the first call arrives, we compute the allocation costs which
are 6.73  for the first ambulance, 6.73 for the second ambulance, 7.43 for the third ambulance, and 7.43 
for the fourth ambulance. The smallest allocation costs
are obtained with ambulances 1 and 2. We therefore need to choose among these ambulances an ambulance which is one of the least advanced. Since both are
advanced, we can either choose ambulance 1 or 2 for BM. We choose ambulance 1. 
When the second call arrives, the allocation costs for this call
are respectively 2.91, 7.35, 6.79, and 7.08  for ambulances 1, 2, 3, and 4. 
The smallest allocation cost is for ambulance 1 which is chosen for the second call.
When the third call arrives, the allocation costs for this call
are respectively 6.78, 14.43, 13.86, and 13,86  for ambulances 1, 2, 3, and 4.
The smallest allocation cost is again obtained with ambulance 1 which is chosen for the third call.
When the fourth call arrives, the allocation costs for this call
are respectively 15.50, 28.31, 28.31, and 27.15 for ambulances 1, 2, 3, and 4. 
The smallest allocation cost is again obtained with ambulance 1 which is also chosen for the fourth call.
We obtain a total allocation cost of 31.92.

\par {\textbf{NonMyopic.}} Step 1. As with BM, we initialize the call counter to 1.
Step 2. We compute the allocation costs for the first call
which are respectively 6.73, 6.73, 7.43, and 7.43 for ambulances 1, 2, 3, and 4.
Step 3. We denote by $\mathcal{S}(1)=\{1,2\}$ the set of ambulances providing the smallest allocation cost
for call 1. Since one of these ambulances is available we send a least advanced ambulance
among ambulances 1 and 2. We choose to send ambulance 1 to call 1 and increase the call index by 1.
We go back to Step 1.
Step 2. When the second call arrives (the call index is 2), we compute the allocation costs which are
respectively 2.91, 7.35, 6.79, and 7.08
with ambulances 1, 2, 3 and 4 for that call. 
Step 3. The smallest allocation cost  is $2.91$ obtained with 
ambulance 1, i.e., $\mathcal{S}(1)=\{1\}$ but this ambulance is busy until time $t_f(1)=28.73$. 
Step 4. We need to check if it would not be better (according to a criterion we specify next) to send this ambulance 1 to another nonallocated call 
that arrives not later than $t_f(1)=28.73$. There are two candidate such calls that arrive not later than $t_f(1)$: calls 3 and 4. The smallest allocation cost for call 3 (among all ambulances)
is $5.25$ obtained with ambulance 1 and the smallest allocation cost for call 4 is $10.04$ also obtained with ambulance 1. We therefore have 
$\mathcal{T}(1)=\{2,3,4\}$ (set of calls that, same as call with index 2, prefer ambulance 1).
Step 5. The maximum of the smallest allocation costs
among the calls in $\mathcal{T}(1)$   
is $\max(2.91,5.25,10.04)=10.04$ obtained for call 4 and we allocate
ambulance 1 to call 4. We go back to Step 1 of NM with call index
$i=2$. Call 2 needs to be allocated and we go to Step 2. 
Step 2. We compute the allocation costs for call with index 2 which are
now respectively 3.67, 7.35, 6,79, and 7.08 for ambulances 1, 2, 3, and 4.
Step 3. The minimal allocation cost is obtained with ambulance 1 and therefore 
$\mathcal{S}(2)=\{1\}$.
This ambulance is busy when call
2 arrives and therefore we check if it would not be better to send that
ambulance to nonallocated calls that arrive not later than 
$t_f(1)=29.49$ (time when ambulance 1 will finishe its service after serving call 4). There is only one such call which is call 3 arriving at 26.3.
Step 4. The minimal allocation costs for this call 3 is 6.79 and the minimal allocation cost
is obtained with ambulance 1 and we therefore have 
$\mathcal{T}(1)=\{2,3\}$.
Step 5. The maximal allocation cost among all
calls in $\mathcal{T}(1)$ is
$\max(3.67,6.79)=6.79$ obtained with call 3 and
we therefore allocate ambulance 1 to call 3. 
We go back to Step 1.
Step 1. Call with index 2 is not allocated.
We go  to Step 2.
Step 2. We compute the allocation costs
for call 2 which are 4.27, 7.35, 6.79, 7.08
with ambulances 1, 2, 3, and 4.
Step 3. The set of ambulances providing the smallest
allocation cost with call with index 2 is
$\mathcal{S}(2)=\{1\}$ (1 is still the best ambulance
for call 2).
Step 4. The time ambulance 1 will finish service
with calls 1, 3, and 4 is now \(t_f(1) = 30.09\). 
There is only one nonallocated call that arrives not later than \(t_f(1) = 30.09\), this is call with index
2. The smallest allocation cost with this call
is 4.27 obtained with ambulance 1.
It follows that $\mathcal{T}(1)$ is now
$\mathcal{T}(1)=\{2\}$. 
Step 5. There is only one element in 
$\mathcal{T}(1)$ and therefore 
we allocate ambulance 1 to call 2 in this set.
We obtain a total allocation cost of 27.82.

\par {\textbf{GHP1.}} Step 1: when the first call arrives at 
\(t = 21\), we  add this call 1 to the queue. The size of the 
queue of calls is 1. Step 2: The queue of calls only contains call 1 with penalized waiting time of zero. 
Step 3: the allocation costs for call 1
are respectively 6.73, 6,73, 7.43, and 7.43
with ambulances 1, 2, 3, and 4 and the set $\mathcal{S}(1)$ of best ambulances for call 1 is
\(\mathcal{S}(1) = \{1,2\}\). Since both ambulances in \(\mathcal{S}(1)\) are advanced, we arbitrarily dispatch one of them, say ambulance 1, to call 1. 
Step 4: We have \(t_f(1) = 28.73\) and compute the type of the next event: a call arrives
or an ambulance finishes service. Call 2 arrives at 26.1, which is less than the minimum of times \(t_f(a)\) over all ambulances \(a\) that are still in service: only 
ambulance 1 is in service and  \(t_f(1) = 28.73\), so the next event will be
"a call arrives".  We go to Step 1.
Steps 1 and 2. The queue only contains call 2 with penalized time 0. 
Step 3. The 
allocation costs for call 2 are 2.91, 7.35, 6.79, 7.08 for ambulances 1, 2, 3, and 4 respectively, implying \(\mathcal{S}(2) = \{1\}\). Ambulance 1 is busy, so call 2 stays in the queue. Step 4. Call 3 arrives at $26.3 < 28.73$, so again the next event is "a call arrives". We go to Step 1.
Step 1. We add call 3 to the queue which is now \(queue = (2,3)\) with penalized times 0.2 and 0 for calls 2 and 3, respectively. Step 2. We order the queue as 
\(queue = (2,3)\), first considering call 2,
 with largest
penalized waiting time. 
Step 3. The allocation costs for call 2 are 2.91, 7.55, 6.99 and 7.28, so \(\mathcal{S}(2) = \{1\}\). The best ambulance for call 2 is 1 which is busy when the
call arrives. For call 3, allocation costs are 5.25 for ambulance 1, 14.42 for ambulance 2, and 13.86 for ambulances 3 and 4, implying \(\mathcal{S}(3) = \{1\}\). 
The best ambulance for call 3 is also busy and
both calls 2 and 3 stay in the queue. 
Step 4. We compute the next event which is "a call arrives" since the next call 4 arrives
at 26.5 and the next ambulance in service will be available again at 
$t_f(1)=28.73$. We go to Step 1.
Steps 1 and 2. We add call 4 to the queue of calls so that the queue of calls
is now \(queue = (3, 2, 4)\) with penalized waiting times 0.4, 0.4, and 0. 
Step 3. The allocation costs for call 3 are 5.25, 14.82, 14.26, and 14.26
for respectively ambulances 1, 2, 3, and 4. The allocation costs for call 2 are 2.91, 7.75, 7.19 and 7.48 for ambulances 1, 2, 3, and 4. The allocation costs for call 4 are 10.04 for ambulance 1, 28.31 for ambulances 2 and 3 and 27.15 for ambulance 4.
Note that ambulance 1 is the best for all calls, but it is still busy, so no dispatch is made.
Step 4. The next event is "ambulance 1 finishes service"
at time 28.73. We go to step 1. 
Steps 1 and 2. The queue of calls is \(queue = (4,3,2)\), sorted in decreasing order of penalized waiting time 
(the penalized waiting times are 8.91, 4.85, and 2.63 for calls 
4, 3, and 2, respectively). Step 3. We consider call 4 which has
the highest penalized waiting time. 
Allocation costs for this call are 10.04 for ambulance 1, 37.21 for ambulances 2 and 3, and 36.06 for ambulance 4, implying $\mathcal{S}(4)=\{1\}$. Ambulance 1 is available and is dispatched to call 4, with 
\(t_f(1) = 29.49\). Step 4. The next event is 
"ambulance 1 finishes service" at \(t_f(1) = 29.49\). 
We go to Step 1.
Step 1. There are sill calls 2 and 3 in the queue
with penalized waiting times (at time \(t_f(1) = 29.49\))
  3.39, and 6.38, respectively.
Step 2. The queue is ordered as \(queue = (3,2)\).
Step 3. We consider call 3 with the largest penalized
waiting time. Allocation costs for call 3 are 6.78 with ambulance 1, 19.28 with ambulance 2, and 18.71 with ambulances 3 and 4.
We have  \(\mathcal{S}(3)=\{1\}\) and we
send ambulance 1 to call 3. 
Step 4. Ambulance 1 will finish its service with call 3 at
\(t_f(1) = 30.09\) and the next event is 
"ambulance 1 finishes service" at \(t_f(1) = 30.09\).
Steps 1 and 2. The queue of calls is now \(queue = (2)\).
Step 3. Allocation costs for call 2 are 4.27, 11.34, 10.78, and 11.06 for ambulances 1, 2, 3, and 4. Ambulance 1 is the best for call 2, is available
at \(t_f(1) = 30.09\), and is sent at this time to call 2.
The queue of calls is now empty and there are no more calls.
Ambulance 1 finishes service at \(t_f(1) = 30.86\). Step 4. The next and final event is "ambulance 1 finishes service" at 30.86.  The total allocation cost is 27.82

\par {\textbf{GHP2.}} \textbf{Step 1.} Call 1 arrives at time 21, and there are no busy ambulances to be released at time 21, so the next event is
"a call arrives at time 21" and we add call 1 to the queue of calls. 
Step 2. Allocation costs for call 1 are  6.73, 6,73, 7.43, and 7.43
for respectively ambulances 1, 2, 3, and 4.
and \(\mathcal{S}(1) = \{1,2\}\).
Step 3. We send to call 1 a least advanced ambulance in \(\mathcal{S}(1) = \{1,2\}\).
Since both are advanced ambulances, we can send 1 or 2 and choose ambulance 1.
Step 4. Ambulance 1 finishes service at 
 \(t_f(1) = 28.73\) and the next event is "Call 2 arrives".
 We go to Step 1.
 Step 1. We add call 2 to the queue of calls.
 Step 2. Allocation costs  for call 2 are 2.91, 7.35, 6.78, and 7.08
 for respectively ambulances 1, 2, 3, and 4 and we have 
 \(\mathcal{S}(2) = \{1\}\). The best ambulance for call 2 is not available
 and we go to Step 4.
 Step 4. The next event is "Call 3 arrives" at time 26.3. We go to Step 1.
 Step 1. We add call 3 to the queue which now contains calls 2 and 3.
 Step 2. Allocation costs for call 2 are 2.91, 7.55, 6.99, 7.28
 with respectively ambulances 1, 2, 3, and 4 while
 for call 3, allocation costs are 5.25, 14.43, 13.86, 13.86
  with respectively ambulances 1, 2, 3, and 4.
The minimal allocation costs are therefore 2.91 for call 2 and 5.25 for call 3.
The maximum of these allocation costs is 5.25 and we handle call 3 first.
The reordered queue of calls is queue=(3,2).
Step 3. \(\mathcal{S}(3) = \{1\}\) and ambulance 1 is still busy  so call 3 stays in queue. 
We also have \(\mathcal{S}(2) = \{1\}\) so call 2 also stays in queue.
Step 4. Call 4 arrives at 26.5 before ambulance 1 is released at 28.73 so 
the next event is "Call 4 arrives". We go to Step 1.
Step 1. We add call 4 to the queue of calls.
Step 2. Allocation costs for call 2 are 2.91, 7.75, 7.18, and 7.48  with respectively ambulances 1, 2, 3, and 4. For call 3, allocation costs are 5.25, 14.83, 14.26, 14.26  with respectively ambulances 1, 2, 3, and 4. For call 4, allocation costs are 10.04, 28.31, 28.31, and 27.15  with respectively ambulances 1, 2, 3, and 4. 
The minimal allocation costs are 2.91, 5.25, and 10.04 for calls 2, 3, and 4.
Step 3. The maximum of these costs is 10.04, obtained for call 4 so call 4 is handled first, then call
3, and then call 2. The best ambulance for these calls is ambulance 1 which is busy so no dispatch
is made at that instant. 
Step 4. The next event is "Ambulance 1 finishes service" at time \(t_f(1) = 28.73\).
Steps 1 and 2. Allocation costs for calls 2, 3, and 4 in queue are given in the matrix
\[
\begin{bmatrix}
2.91 & 9.98 & 9.41 & 9.70 \\
5.25 & 19.28 & 18.71 & 18.71 \\
10.04 & 37.21 & 37.21 & 36.06
\end{bmatrix}
\]
with column j corresponding to ambulance j and line 1 to call 2, line 2 to call 3, and
line 3 to call 4. 
Step 3. The minimal allocation costs are
2.91, 5.25, and 10.04 for respectively
calls 2, 3, and 4.
The maximum of these minimal allocation costs
is 10.04 and
again, call 4 is handled first. Ambulance 1 is dispatched to call 4, with \(t_f(1) = 29.49\), and the allocation costs of calls 2 and 3 are updated correspondingly: call 2 allocation costs are now 3.67, 9.98, 9.41, 9.70 with respectively ambulances 1, 2, 3, and 4 and call 3 allocation costs are now 6.79, 19.28, 18.71, and 18.71 with respectively ambulances 1, 2, 3, and 4. Call 3 and 2 are handled in this order, but since ambulance 1 is the best for both and is not available at $t_f(1)$, no dispatch is made. The next event is
"Ambulance 1 finishes service" at time $t_f(1)=29.49$.
We go to Step 1. Steps 1 and 2. Allocation costs for calls 2 and 3 are:
\[
\begin{bmatrix}
3.67 & 10.74 & 10.18 & 10.46 \\
6.79 & 20.81 & 20.25 & 20.25
\end{bmatrix}
\]
The minimal allocation costs are 3.67 for call 2
and 6.78 for call 3 and the maximum of 3.67 and 6.78
(corresponding to call 3) is 6.78 so we handle call 3
first.
Step 3. We dispatch ambulance 1 to call 3, with \(t_f(1) = 30.09\). The best ambulance for call 3 is ambulance
1 which is busy so no dispatch is made for call 3.
Step 4. We compute the next event which is
"Ambulance 1 finishes service" at time \(t_f(1) = 30.09\).
Steps 1 and 2. Updated allocation costs for call 2 are 
4.27, 11.34. 10.78 and 11.07 with respectively ambulances
1, 2, 3, and 4.  \textbf{Step 3.} Ambulance 1 is dispatched to call 2.
There are no more calls in queue and we compute the
total allocation cost which is 27.82.

{\textbf{Closest available.}} At the first call, arriving at time 21, all ambulances are available. The following allocation costs are computed: 6.73 for ambulances 1 and 2, 7.43 for ambulances 3 and 4. Ambulance 1 is dispatched, with \(t_f(1) = 28.73\). The second call arrives at time 26.1, and ambulances 2, 3, and 4 are available. The allocation costs are 7.35 for ambulance 2, 6.79 for ambulance 3 and 7.08 for ambulance 4. Ambulance 3 is dispatched, with \(t_f(3) = 33.37\). The third call arrives at 26.3 and ambulances 2 and 4 are available. The allocation costs are 14.43 for ambulance 2 and 13.86 for ambulance 4. Ambulance 4 is dispatched, with \(t_f(4) = 33.63\). The last call arrives at 26.5, and ambulance 2 - the only available - is dispatched, with \(t_f(2) = 34.05\). The total allocation cost is 55.68.

In this example, the difference
in total penalized response time
between the best three heuristics
(which in this example provide
the same allocations) and the Best Myopic heuristic, as well as the difference in total penalized
response time between the Best Myopic heuristic and the 
CA heuristic 
 are increasing functions of
$t$ and can therefore be made
arbitrarily large taking $t$
sufficiently large.
\end{example}

\begin{example}
We can modify
the previous Example
\ref{ex1toy} in such a way
that the NonMyopic
heuristic is better than
all others. Consider the
previous example 
with $t=1$
replacing the time
the ambulance waits at hospital
by 0.4, 2, 4, 2,
the time the ambulance stays
on the scene of the call by
1, 3, 5, and 2, and the priorities
of the calls by low, low, high,
and intermediate for calls
arriving at 1, 6.1, 6.3, and 
6.5 respectively.
For this new instance,
the total penalized
waiting time
is 41.2 with the NonMyopic
heuristic, 50.9 with the
GHP1 and GHP2 heuristics,
51.8 with the
Best Myopic heuristic,
55.4 with the CA heuristic.
\end{example}

\begin{example}
The previous two examples
are examples where
the NonMyopic heuristic
is the best. However,
it can be worse than some of the  other
heuristics we presented.
Indeed, consider
an instance in the plane
with two calls of low priority:
the first one
arriving at $t=0.5$ at location 
$(5,9)$, the second one
arriving at $t=0.6$ at location 
$(8,10)$.
There is one hospital
at location $(5,10)$.
All ambulances are advanced
and ambulance 1 is at
its station $(3,0)$
at $t=0$, ambulance 2
is at
its station $(15,10)$
at $t=0$, and ambulance
3 is at the hospital at
$t=0$. This ambulance will
be available again at
$t=1$. The time the ambulance
spends on the scene of the calls
is 0.1 and the ambulances
wait for 10 time units
at hospital before being available
again.
On this example,
the Best Myopic heuristic
and GHP1
provide the smallest
total penalized response
time (7.1) obtained sending ambulance
3 to the first call and then
to the second call.
The NonMyopic heuristic
sends ambulance 3 first to
the second call and then
to the first call, yielding
a total penalized response
time of 11.1.
GHP2
gives a total penalized waiting
time of 11.1 obtained sending
the first ambulance to the first
call and ambulance 3 to the second call (ambulance 3 to the first and  second call). The remaining heuristic
CA is the worst on this example
with a total penalized response
time of 16.2 obtained sending
ambulance 1 to the first call
and ambulance 2 to the second call.
\end{example}

\if{\section{A New Preparedness Metric}
    \label{sec:preparedness}
    
    In this section we propose a continuous-time Markov chain model to measure the preparedness of an ambulance station with a given set of ambulances at the station.
    The model is a local approximation in the sense that it considers only a single station, and it considers only the currently available ambulances at the station (thus it is local in both space and time).
    The benefit of this local approximation is that it is tractable while providing a fairly accurate metric of the ability to dispatch ambulances to emergencies in the zone of the station without delay.
    
    The model considers a single station at a time, and therefore the notation does not indicate the station being considered.
    For each emergency type $c \in C$, let $\lambda(c)$ denote the arrival rate of type~$c$ emergencies to the zone of the station, and let $A(c)$ denote the set of ambulance types that can serve a type~$c$ emergency.
    Conversely, for each ambulance type $a \in A$, let $C(a) \defi \{c \in C \, : \, a \in A(c)\}$ denote the set of emergency types that can be served by ambulance type~$a$.
    We assume that associated with each emergency type $c \in C$ there is a total order $\succ_{c}$ on $A(c)$ that governs the type of ambulance dispatched to serve a type~$c$ emergency, as follows.
    If a type~$c$ emergency arrives, and an ambulance of type $a \in A(c)$ is available, then such an ambulance is dispatched to the emergency if and only if no ambulance of type $a' \in A(c)$ such that $a' \succ_{c} a$ is available, that is, the most preferred ambulance type for the emergency type that is available is dispatched.
    Note that $a' \succ_{c} a$ does not mean that type~$a'$ ambulances are more capable than type~$a$ ambulances.
    For example, for a minor emergency type~$c$ it may hold that BLS $\succ_{c}$ ALS.
    If a type~$a$ ambulance is dispatched to a type~$c$ emergency in the zone of the station, then the ambulance becomes available again after an exponentially distributed amount of time with mean $1 / \mu(a,c)$.
    
    For each ambulance type $a \in A$, let $m(a)$ denote the number of type~$a$ ambulances at the station, and let $m \defi (m(a), a \in A)$.
    The state of the Markov chain is specified by the number of ambulances of each type that is busy serving emergencies of each type.
    For each emergency type $c \in C$ and ambulance type $a \in A(c)$, let $x_{a,c}$ denote the number of type~$a$ ambulances currently busy serving type~$c$ emergencies.
    It is assumed that each ambulance not busy serving an emergency is available, and thus the number of available type~$a$ ambulances is given by $m(a) - \sum_{c \in C(a)} x_{a,c}$.
    Let $x \defi (x_{a,c}, a \in A, c \in C(a))$ denote the state of the Markov chain.
    
    Let $q^{m}_{x,x'}$ denote the transition rate from state~$x$ to state~$x'$.
    There are two types of state transitions: transitions caused by emergency arrivals and subsequent ambulance dispatches, and transitions caused by ambulances completing tasks and becoming available.
    Consider any current state~$x$ and any emergency type~$c$ such that there is an ambulance available to serve an emergency of that type, that is, $\sum_{a \in A(c)} \left[m(a) - \sum_{c' \in C(a)} x_{a,c'}\right] > 0$.
    Let
    \[
    \hat{a} \ \ \in \ \ \left\{a \in A(c) \; : \; m(a) > \sum_{c' \in C(a)} x_{a,c'}, \; \nexists \, a' \in A(c) \mbox{ s.t.\ } m(a') > \sum_{c' \in C(a')} x_{a',c'}, \; a' \succ_{c} a\right\}
    \]
    denote the available ambulance type that is preferred for type~$c$ emergencies.
    Note that $\hat{a}$ depends on~$x$ and~$c$, but the notation does not show the dependence.
    Then $q^{m}_{x,x'} = \lambda(c)$ where $x'$ is determined as follows:
    \begin{eqnarray*}
    x'_{\hat{a},c} & = & x_{\hat{a},c} + 1 \\
    x'_{a',c'} & = & x_{a',c'} \qquad \forall \ (a',c') \neq (\hat{a},c).
    \end{eqnarray*}
    In addition, for each emergency type $c \in C$ and ambulance type $a \in A(c)$ such that $x_{a,c} > 0$,
    $q^{m}_{x,x'} = \mu(a,c) x_{a,c}$ where $x'$ is determined as follows:
    \begin{eqnarray*}
    x'_{a,c} & = & x_{a,c} - 1 \\
    x'_{a',c'} & = & x_{a',c'} \qquad \forall \ (a',c') \neq (a,c).
    \end{eqnarray*}
    Let $q^{m}_{x,x} \defi - \sum_{x'} q^{m}_{x,x'}$.
    Note that all states communicate, for example, any state~$x'$ can be reached from state~$x = 0$ in a finite number of transitions and vice versa.
    Thus, the Markov chain has a unique stationary distribution~$\nu^{m} = (\nu^{m}_{x})$ given by
    \begin{eqnarray}
    \label{eqn:balance equations}
    \sum_{x} \nu^{m}_{x} q^{m}_{x,x'}  & = & 0 \qquad \forall \ x' \\
    \sum_{x} \nu^{m}_{x} & = & 1.
    \label{eqn:probability scaling}
    \end{eqnarray}
    
    The system of equations~\eqref{eqn:balance equations}--\eqref{eqn:probability scaling} may be quite large, and thus care should be taken to solve the system.
    Note that the system of equations~\eqref{eqn:balance equations}--\eqref{eqn:probability scaling} has one more equation than unknown, and the equation $\sum_{x} \nu^{m}_{x} q^{m}_{x,x'} = 0$ for one $x'$, such as $x' = 0$, can be dropped.
    Let $Q^{m} \defi (q^{m}_{x,x'})$ denote the transition rate matrix, let $\mathbf{0}$ denote the vector $(0,\ldots,0)$, and let $\mathbf{1}$ denote the vector $(1,\ldots,1)$.
    Then the system~\eqref{eqn:balance equations}--\eqref{eqn:probability scaling} can be written as ${Q^{m}}^{\top} \nu^{m} = \mathbf{0}, \mathbf{1}^{\top} \nu^{m} = 1$.
    Let $D^{m} \defi (\mathbf{1}, (q^{m}_{x,x'}, x' \neq 0))$ denote the reduced matrix.
    Then the system~\eqref{eqn:balance equations}--\eqref{eqn:probability scaling} is equivalent to ${D^{m}}^{\top} \nu^{m} = (1,0,\ldots,0)$.
    The transition rate matrix~$Q^{m}$ is very sparse, and thus~$D^{m}$ is sparse.
    A good solver should exploit the sparsity of~$D^{m}$.
    One possibility to solve ${D^{m}}^{\top} \nu^{m} = (1,0,\ldots,0)$ is to consider $D^{m} {D^{m}}^{\top} \nu^{m} = D^{m} (1,0,\ldots,0) = \mathbf{1}$, and to use the conjugate gradient method to solve $D^{m} {D^{m}}^{\top} \nu^{m} = \mathbf{1}$.
    The matrix $D^{m} {D^{m}}^{\top}$ should be diagonally dominant, and thus the diagonal of $D^{m} {D^{m}}^{\top}$ can be used as a preconditioner.
    Also, the matrix $D^{m} {D^{m}}^{\top}$ should still be quite sparse (but not as sparse as $D^{m}$).
    Alternatively, a variant of the conjugate gradient method can be used for the system ${D^{m}}^{\top} \nu^{m} = (1,0,\ldots,0)$ with asymmetric~$D^{m}$.
    
    Given stationary distribution~$\nu^{m} = (\nu^{m}_{x})$, one can compute steady-state costs such as the following.
    Let $\phi(c)$ denote the penalty if there is no ambulance available to serve a type~$c$ emergency.
    Then the cost rate $\psi_{x}$ while in state~$x$ is given by
    \[
    \psi^{m}_{x} \ \ = \ \ \sum_{c \in C} \mathds{1}\left\{\sum_{a \in A(c)} \left[m(a) - \sum_{c' \in C(a)} x_{a,c'}\right] = 0\right\} \lambda(c) \phi(c)
    \]
    and the steady state cost rate $\bar{\psi}^{m}$ given ambulance supply~$m$ is $\bar{\psi}^{m} = \sum_{x} \nu^{m}_{x} \psi^{m}_{x}$.
    
    For each station~$b$, $\bar{\psi}^{m}_{b}$ can be computed in advance for each $m \defi (m(a), a \in A)$.
    Then $\bar{\psi}^{m}_{b}$ can be used as a ``preparedness'' metric for ambulance supply~$m$ at station~$b$.
    For example, when an ambulance becomes available after completing a task, it can be sent to the station~$b$ where it will improve $\bar{\psi}^{m}_{b}$ the most.
    More specifically, for each station~$b$, let $m_{b}$ denote the current ambulance supply at station~$b$,
    and let $m^{+}_{b}$ denote the ambulance supply at station~$b$ if the newly available ambulance would be added to the ambulance supply at station~$b$.
    Then send the newly available ambulance to the station $b^* \in \arg\max\left\{\bar{\psi}^{m_{b}}_{b} - \bar{\psi}^{m^{+}_{b}}_{b} \, : \, b \in B\right\}$.

    \section{A New Ambulance Selection Heuristic}
    \label{sec:selection}
    
    Suppose that a type~$c$ emergency arrives at location~$\ell$.
    Let $i_{0}$ be the index of the newly arrived emergency.
    For any station~$b$, let $\mathcal{A}_{b}$ denote the set of ambulances currently at station~$b$ or en route to station~$b$, and let $\mathcal{A}'$ denote the set of currently busy ambulances.
    Thus the set of ambulances is $\mathcal{A} = \cup_{b \in B} \mathcal{A}_{b} \cup \mathcal{A}'$.
    For each ambulance $a \in \mathcal{A}$, let $t(a) \in A$ denote the type of the ambulance.
    For each station~$b$, let $m_{b} = (m_{b}(\tilde{a}), \tilde{a} \in A)$ denote the current ambulance supply at station~$b$, with the understanding that any ambulance currently en route to a station is included in the ambulance supply of that station.
    Let $\mathcal{Q}$ denote the set of emergencies currently in queue.
    For each emergency $i \in \mathcal{Q} \cup \{i_{0}\}$ and each ambulance $a \in \mathcal{A}$, let $r(a,i)$ denote the weighted response time if ambulance~$a$ is dispatched to emergency~$i$.
    Note that if ambulance~$a$ is currently busy, then the calculation of $r(a,i)$ is based on ambulance~$a$ first completing its current task and then traveling to emergency~$i$.
    
    For each emergency $i \in \mathcal{Q} \cup \{i_{0}\}$ and each ambulance $a \in \mathcal{A}$, let $x(a,i)$ denote a decision variable that is~$1$ if ambulance~$a$ is dispatched to emergency~$i$, and is~$0$ otherwise.
    For each busy ambulance $a \in \mathcal{A}'$ and station~$b \in B$, let $y(a,b)$ denote a decision variable that is~$1$ if ambulance~$a$ is sent to station~$b$, and is~$0$ otherwise.
    Let $x \defi (x(a,i), a \in \mathcal{A}, i \in \mathcal{Q} \cup \{i_{0}\})$ and $y \defi (y(a,b), a \in \mathcal{A}', b \in B)$.
    For any ambulance type~$\tilde{a} \in A$ and station~$b \in B$, let $m^{+}_{b}(\tilde{a},x,y) \defi m_{b}(\tilde{a}) - \sum_{a \in \mathcal{A}_{b}} \mathds{1}\left\{t(a) = \tilde{a}\right\} \sum_{i \in \mathcal{Q} \cup \{i_{0}\}} x(a,i) + \sum_{a \in \mathcal{A}'} \mathds{1}\left\{t(a) = \tilde{a}\right\} y(a,b)$, and let $m^{+}_{b}(x,y) \defi (m^{+}_{b}(\tilde{a},x,y), \tilde{a} \in A)$.
    
    Let $\Gamma$ be a parameter that weighs the effects of preparedness (uncertain future response times) relative to current response times.
    For each emergency $i \in \mathcal{Q} \cup \{i_{0}\}$, let $\gamma(i)$ be a parameter that weighs the effect of queueing emergency~$i$; $\gamma(i)$ should depend on the type of emergency~$i$.
    Version~1 of the ambulance selection heuristic solves the following optimization problem:
    \begin{eqnarray}
    \min_{x,y} & & \sum_{a \in \mathcal{A}} \sum_{i \in \mathcal{Q} \cup \{i_{0}\}} r(a,i) x(a,i) + \sum_{i \in \mathcal{Q} \cup \{i_{0}\}} \gamma(i) \left[1 - \sum_{a \in \mathcal{A}} x(a,i)\right] + \Gamma \sum_{b \in \mathcal{B}} \bar{\psi}^{m^{+}_{b}(x,y)}_{b}
    \label{eqn:nonlinear selection heuristic objective} \\
    \mbox{s.t.} & & \sum_{i \in \mathcal{Q} \cup \{i_{0}\}} x(a,i) \ \ \le \ \ 1 \qquad \forall \ a \in \cup_{b \in B} \mathcal{A}_{b} \\
    & & \sum_{i \in \mathcal{Q} \cup \{i_{0}\}} x(a,i) + \sum_{b \in B} y(a,b) \ \ = \ \ 1 \qquad \forall \ a \in \mathcal{A}' \\
    & & \sum_{a \in \mathcal{A}} x(a,i) \ \ \le \ \ 1 \qquad \forall \ i \in \mathcal{Q} \cup \{i_{0}\}
    \end{eqnarray}
    
    The objective function~\ref{eqn:nonlinear selection heuristic objective} is not linear in $(x,y)$.
    Version~2 of the heuristic solves a linear optimization problem.
    For each ambulance type~$\tilde{a} \in A$, let $e(\tilde{a})$ denote the unit vector with the component for $\tilde{a}$ equal to~$1$ and the components for all $\tilde{a}' \in A \setminus \{\tilde{a}\}$ equal to~$0$.
    For each ambulance type~$\tilde{a} \in A$ and station~$b \in B$, let $s^{+}(\tilde{a},b) \defi \bar{\psi}^{m_{b} + e(\tilde{a})}_{b} - \bar{\psi}^{m_{b}}_{b}$ and $s^{-}(\tilde{a},b) \defi \bar{\psi}^{m_{b} - e(\tilde{a})}_{b} - \bar{\psi}^{m_{b}}_{b}$.
    Then Version~2 of the ambulance selection heuristic solves the following optimization problem:
    \begin{eqnarray}
    \min_{x,y} & & \sum_{a \in \mathcal{A}} \sum_{i \in \mathcal{Q} \cup \{i_{0}\}} r(a,i) x(a,i) + \sum_{i \in \mathcal{Q} \cup \{i_{0}\}} \gamma(i) \left[1 - \sum_{a \in \mathcal{A}} x(a,i)\right] \nonumber \\
    & & \ {} + \Gamma \sum_{b \in \mathcal{B}} \sum_{a \in \mathcal{A}_{b}} \sum_{i \in \mathcal{Q} \cup \{i_{0}\}} s^{-}(t(a),b) x(a,i) + \Gamma \sum_{a \in \mathcal{A}'} \sum_{b \in \mathcal{B}} s^{+}(t(a),b) y(a,b)
    \label{eqn:linear selection heuristic objective} \\
    \mbox{s.t.} & & \sum_{i \in \mathcal{Q} \cup \{i_{0}\}} x(a,i) \ \ \le \ \ 1 \qquad \forall \ a \in \cup_{b \in B} \mathcal{A}_{b} \\
    & & \sum_{i \in \mathcal{Q} \cup \{i_{0}\}} x(a,i) + \sum_{b \in B} y(a,b) \ \ = \ \ 1 \qquad \forall \ a \in \mathcal{A}' \\
    & & \sum_{a \in \mathcal{A}} x(a,i) \ \ \le \ \ 1 \qquad \forall \ i \in \mathcal{Q} \cup \{i_{0}\}
    \end{eqnarray}
    
    Let $(x^*,y^*)$ denote an optimal solution of Version~1 or Version~2 of the ambulance selection heuristic.
    With both versions of the heuristic, the ambulance selection decision dispatches ambulance~$a$ to newly arrived emergency~$i_{0}$ if and only if $x^*(a,i_{0}) = 1$.

    \section{A New Ambulance Reassignment Heuristic}
    \label{sec:reassignment}
    
    Suppose that ambulance~$a_{0} \in \mathcal{A}$ becomes available after completing a task.
    The rest of the notation is the same as in Section~\ref{sec:selection}.
    Version~1 of the ambulance reassignment heuristic solves the following optimization problem:
    \begin{eqnarray}
    \min_{x,y} & & \sum_{a \in \mathcal{A}} \sum_{i \in \mathcal{Q}} r(a,i) x(a,i) + \sum_{i \in \mathcal{Q}} \gamma(i) \left[1 - \sum_{a \in \mathcal{A}} x(a,i)\right] + \Gamma \sum_{b \in \mathcal{B}} \bar{\psi}^{m^{+}_{b}(x,y)}_{b}
    \label{eqn:nonlinear reassignment heuristic objective} \\
    \mbox{s.t.} & & \sum_{i \in \mathcal{Q}} x(a,i) \ \ \le \ \ 1 \qquad \forall \ a \in \cup_{b \in B} \mathcal{A}_{b} \\
    & & \sum_{i \in \mathcal{Q}} x(a,i) + \sum_{b \in B} y(a,b) \ \ = \ \ 1 \qquad \forall \ a \in \mathcal{A}' \cup \{a_{0}\} \\
    & & \sum_{a \in \mathcal{A}} x(a,i) \ \ \le \ \ 1 \qquad \forall \ i \in \mathcal{Q}
    \end{eqnarray}
    Then Version~2 of the ambulance reassignment heuristic solves the following optimization problem:
    \begin{eqnarray}
    \min_{x,y} & & \sum_{a \in \mathcal{A}} \sum_{i \in \mathcal{Q}} r(a,i) x(a,i) + \sum_{i \in \mathcal{Q}} \gamma(i) \left[1 - \sum_{a \in \mathcal{A}} x(a,i)\right] \nonumber \\
    & & \ {} + \Gamma \sum_{b \in \mathcal{B}} \sum_{a \in \mathcal{A}_{b}} \sum_{i \in \mathcal{Q}} s^{-}(t(a),b) x(a,i) + \Gamma \sum_{a \in \mathcal{A}' \cup \{a_{0}\}} \sum_{b \in \mathcal{B}} s^{+}(t(a),b) y(a,b)
    \label{eqn:linear reassignment heuristic objective} \\
    \mbox{s.t.} & & \sum_{i \in \mathcal{Q}} x(a,i) \ \ \le \ \ 1 \qquad \forall \ a \in \cup_{b \in B} \mathcal{A}_{b} \\
    & & \sum_{i \in \mathcal{Q}} x(a,i) + \sum_{b \in B} y(a,b) \ \ = \ \ 1 \qquad \forall \ a \in \mathcal{A}' \cup \{a_{0}\} \\
    & & \sum_{a \in \mathcal{A}} x(a,i) \ \ \le \ \ 1 \qquad \forall \ i \in \mathcal{Q}
    \end{eqnarray}
    Let $(x^*,y^*)$ denote an optimal solution of Version~1 or Version~2 of the ambulance reassignment heuristic.
    With both versions of the heuristic, the ambulance reassignment decision dispatches ambulance~$a_{0}$ to emergency~$i \in \mathcal{Q}$ if and only if $x^*(a_{0},i) = 1$, and dispatches ambulance~$a_{0}$ to station~$b \in B$ if and only if $y^*(a_{0},b) = 1$.
    
    \section{Optimal allocation of ambulances to calls in queues} \label{sec:optqueue}
    
    In this section, we propose
    an optimization problem for the optimal allocation
    of ambulances to a set
    of calls that has just arrived. 
    It can be used by the heuristics we have
    presented when a set of calls
    arrives simultaneously.
    This allocation model is an adaptation of the model presented in
     with the following modifications:
    (i) instead of being initially available and located at hospitals,
    the ambulances can be anywhere, at a station, going back to a station,
    or in service; (ii) we consider 
    three call types 
    instead of 2;
    (iii) we consider different ambulance types; (iv) we consider
    the possibility to clean an ambulance after service. These modifications require modifying the set of variables
    and constraints.
    By a slight abuse of notation,
    for $i=1,2,3,4$, we denote by
    $C_{i}$ the set of calls requiring 
    a ride of type ($C_{i}$) detailed
    in the introduction.\\
    
    \par {\textbf{Decision variables.}}\\ 
    \par For every call $i$ and 
    ambulance $k$ with 
    $k \in \mathcal{A}(\mbox{Type}_{c}(i))$
    (recall that $\mathcal{A}(\mbox{Type}_{c}(i))$
    is the set of ambulances that can attend
    calls of type Type$_{c}(i)$), let $x(i,k) \in \{0,1\}$ (resp. $z(i,k) \in \{0,1\}$)  be a
    binary variable which is 1
    if $i$ is the first (resp. last) call served by ambulance $k$ and 0 otherwise.
    For calls $i,j$, and 
    ambulance $k \in 
    \mathcal{A}(\mbox{Type}_{c}(i))\cap \mathcal{A}(\mbox{Type}_{c}(j))$ let  $x(i,j,k) \in \{0,1\}$ be a binary variable which is 1
    if call $i$ is served by ambulance $k$ immediately before $j$ and 0 otherwise.
    For every hospital h,
    and cleaning station cb,
    \begin{itemize}
    \item for every $i \in C_{1}$, let 
    $y(i,h,cb)$ be a binary variable which 
    is 1 if the patient of  call  $i$
    is sent to hospital $h$ and the corresponding ambulance 
    is then sent to 
    cleaning station cb;
    \item for every call $i \in C_{2}$, let $y(i,h)$
    be a binary variable which is 1 if the patient of  call  $i$
    is sent to hospital $h$;
    \item for every call $i \in C_{3}$, let  $y(i,cb)$
    be a binary variable which is 1 if the 
    ambulance serving call $i$
    is sent to cleaning station cb after leaving the scene
    of call $i$.
    \end{itemize}
    We will also use continuous variable $t_{i}$ which for
    every call $i$ represents the arrival time of the ambulance
    serving that call on the scene of the call.
    We will also denote by $C$ the set of calls, by H the set of hospitals,
    by CB the set of cleaning stations, and by 
    $A$ the set of ambulances.
    We will use the simplifying assumption
    that function travelTime now
    only depends on origin and
    destination: travelTime(A,B)
    is the time for an ambulance to go
    from A to B.
    Finally, $M_{i}$ is the maximal completion time 
    among all calls of type $i$ (the completion time is the time the ambulance
    arrives on the scene of the call for calls in
    $C_{3} \cup C_{4}$ and the time the ambulance arrives at hospital
    for calls in $C_{1} \cup C_{2}$).\\
    
    \par {\textbf{Constraints on the rides.}}\\
    
    \par An ambulance has at most one first visited patient which can be written:
    \begin{equation}\label{eq1}
    \sum_{i \in C: k \in \mathcal{A}(\mbox{Type}_{c}(i))} x(i,k) \leq 1,\;\forall k \in A.
    \end{equation}
    
    \par Every call $i \in C$ must be attended by one ambulance
    which gives
    \begin{equation}\label{eq2}
    \sum_{k: 
    k \in \mathcal{A}(
    \mbox{Type}_{c}(i))} x(i,k) +
    \sum_{k: k \in \mathcal{A}(\mbox{Type}_{c}(i))}
    \sum_{j \in C, j\neq i,
    k \in \mathcal{A}(\mbox{Type}_{c}(j))}  x(j,i,k) =1,\;
    \forall i \in C.
    \end{equation}
    
    \par If ambulance  $k$ arrives on the scene of call 
    $i$ it will have to leave that scene:
    for all $i \in C$ and 
    $k \in \mathcal{A}(\mbox{Type}_{c}(i))$ we have
    \begin{equation}\label{eq3}
    x(i,k) + \sum_{j \in C, j \neq i: k \in \mathcal{A}(\mbox{Type}_{c}(j))} x(j,i,k) =z(i,k) + 
    \sum_{j \in C, j \neq i: k \in \mathcal{A}(\mbox{Type}_{c}(j))} x(i,j,k).
    \end{equation}
    
    \par Every patient of every call $i \in C_{1}$ must go to a hospital
    and the corresponding ambulance must go to a cleaning station:
    \begin{equation}\label{eq6}
    \sum_{h \in H} \sum_{cb \in CB} y(i,h,cb) =1,\;\forall  i \in C_{1}.
    \end{equation}
    
    \par Every patient of every call $i \in C_{2}$ must go to a hospital:
    \begin{equation}\label{eq5}
    \sum_{h \in H} y(i,h) =1,\;\forall  i \in C_{2}.
    \end{equation}
    \par Every call $i \in C_{3}$ must go to one cleaning station:
    \begin{equation}\label{eq4}
    \sum_{cb \in CB} y(i,cb) =1,\;\forall  i \in C_{3}.
    \end{equation}
    
    \par Hospital $h$ can receive at most $c_h$ patients: if
    the current number of patients in hospital $h$ is $n_h$ then
    \begin{equation}\label{eq7}
    \sum_{i \in C_{2}} y(i,h)
    + \sum_{cb \in CB} 
    \sum_{i \in C_{1}} y(i,h,cb) \leq c_h - n_h,\;\forall  h \in H.
    \end{equation}\\
    
    \par {\textbf{Constraints on the arrival times of the first patients served.}}\\
    \par  For every ambulance $k$ serving a first call $i$
    three situations can happen at $t_{0}:=$currentTime (the instant the set of calls arrived and the ambulance
    allocation problem is solved): the ambulance
    is at a station, i.e., (a) $t_b(k)\leq t_{0}$,
    the ambulance is going to a station, i.e., (b) 
    $t_{f}(k) \leq t_{0} < t_b(k)$, and
    the ambulance is in service (c) $t_{0}<t_{f}(k)$.
    Denoting by $M$ an upper bound on the arrival time
    of the ambulances on the scenes of the calls, we
    get three corresponding sets of 
    constraints
    for the arrival time $t_{i}$ at the first call.
    
    \begin{itemize}
    \item[(a)] For every $i \in C$ and every  ambulance $k$
    such that $t_b(k)\leq t_{0}$ and 
    $k \in \mathcal{A}(\mbox{Type}_{c}(i))$, we have
    \begin{equation}\label{eq8}
    t_{0} + \mbox{travelTime}(\ell_b(k),\ell_{c}(i))
    \leq 
    t_{i} + M(1-x(i,k)).
    \end{equation}
    \item[(b)] For every $i \in C$ and every  ambulance 
    $k \in \mathcal{A}(\mbox{Type}_{c}(i))$
    such that $t_{f}(k)\leq t_{0}<t_b(k)$, we have   
    \begin{equation}\label{eq9}
    t_{0} + \mbox{travelTime}(P_k,\ell_{c}(i)) \leq 
    t_{i} + M(1-x(i,k))
    \end{equation}
    where
    $$
    P_k=
    \mbox{positionBetweenOriginDestination}(\ell_{f}(k), \ell_b(k), t_{f}(k), 
    t_{0}).
    $$
    \item[(c)] For every $i \in C$ and every  ambulance $k$
    such that $t_{0}<t_{f}(k)$ and $k \in \mathcal{A}(\mbox{Type}_{c}(i))$, we have
    \begin{equation}\label{eq10}
    t_{f}(k) + \mbox{travelTime}(\ell_{f}(k),\ell_{c}(i))
    \leq 
    t_{i} + M(1-x(i,k)).
    \end{equation}
    \end{itemize}
    
    \par {\textbf{Constraints on the arrival times for calls that
    are not served first by an ambulance.}}\\
    
    \par For every $i \in C_{1}$, every $h \in H$, every
    $cb \in CB$,
    and every $j \neq i \in C$, we have
    \begin{equation}\label{eq14}
    \begin{array}{l}
    t_{i} + \mbox{TimeOnScene}_{c}(i)\\
     + 
    y(i,h,cb)\Big(
    \mbox{travelTime}(\ell_{c}(i),h)
    +\mbox{TimeAtHospital}_{c}(i) + 
    \mbox{travelTime}(h,cb)\Big)\\
    + y(i,h,cb)\Big(\mbox{CleaningTime}_{c}(i)
    + \mbox{travelTime}(cb,\ell_{c}(j))
    \Big)\\
     \leq 
    t_{j} + M\Big(1- \displaystyle \sum_{k: 
    k \in \mathcal{A}(\mbox{Type}_{c}(i)) \bigcap \mathcal{A}(\mbox{Type}_{c}(j))}x(i,j,k)\Big).
    \end{array}
    \end{equation}
    
    \par By a slight abuse of notation, in function travelTime above, cb should be understood as the location
    of cleaning station cb
    and h as the location of hospital h.\\

    \par For every $i \in C_{2}$, every $h \in H$, and every $j \neq i \in C$, we have
    \begin{equation}\label{eq13}
    \begin{array}{l}
    t_{i} + \mbox{TimeOnScene}_{c}(i)\\
     + 
    y(i,h)\Big(
    \mbox{travelTime}(\ell_{c}(i),h)
    +\mbox{TimeAtHospital}_{c}(i) + 
    \mbox{travelTime}(h,\ell_{c}(j))\Big)\\
     \leq 
    t_{j} + M\Big(1- \displaystyle \sum_{k:k \in \mathcal{A}(\mbox{Type}_{c}(i)) \bigcap \mathcal{A}(\mbox{Type}_{c}(j))}x(i,j,k)\Big).
    \end{array}
    \end{equation}

    \par For every $i \in C_{3}$, every $cb \in CB$, and every $j \neq i \in C$, we have
    \begin{equation}\label{eq12}
    \begin{array}{l}
    t_{i} + \mbox{TimeOnScene}_{c}(i)\\
     + 
    y(i,cb)\Big(
    \mbox{travelTime}(\ell_{c}(i),cb)
    +\mbox{CleaningTime}_{c}(i) + 
    \mbox{travelTime}(cb,\ell_{c}(j))\Big)\\
     \leq 
    t_{j} + M\Big(1- \displaystyle \sum_{k:k \in \mathcal{A}(\mbox{Type}_{c}(i)) \bigcap \mathcal{A}(\mbox{Type}_{c}(j))}x(i,j,k)\Big).
    \end{array}
    \end{equation}
    
    \par For every $i \in C_{4}$ and every $j \neq i \in C$, we have
    \begin{equation}\label{eq11}
    \begin{array}{l}
    t_{i} + \mbox{TimeOnScene}_{c}(i) + 
    \mbox{travelTime}(\ell_{c}(i),\ell_{c}(j))\\
      \leq 
    t_{j} + M\Big(1- \displaystyle \sum_{k:k \in \mathcal{A}(\mbox{Type}_{c}(i)) \bigcap \mathcal{A}(\mbox{Type}_{c}(j))}x(i,j,k)\Big).
    \end{array}
    \end{equation}
    
    \par Assuming a positive service time for each call (for instance
    if TimeOnScene$_{c}(i)>0$) the constraints above imply the following required properties
    of decision variables $x(i,j,k)$,
    $x(i,k)$, $z(i,k)$, of the rides:
    \begin{itemize}
    \item if $x(i,k)=1$, i.e., if call $i$ is the first served by ambulance
    $k$ then from \eqref{eq2} $x(j,i,k)=0$ for all $j \neq i$, i.e,
    no patient j is served before i by k and there is $j$ such that
    $j$ is the last call attended by $k$ meaning that
    $z(j,k)=1$. Indeed, in this situation, either $z(i,k)=1$ or there is 
    $j_{1}$ such that  $x(i,j_{1},k)=1$. Then either 
    $z(j_{1},k)=1$ or there is $j_{2}$ such that
    $x(j_{1},j_{2},k)=1$. Observe that $j_{2} \neq j_{1}$ since
    otherwise from the constraints on the arrival times we
    would have $t_{j_{2}}>t_{j_{1}}=t_{j_{2}}$ which yields
    a contradiction.
    Since there is a finite number of calls, there can only be a finite
    number, say $n$, of calls $j_{1},\ldots,j_n$, all different, such that
    $i$ is served before $j_{1}$, $j_{1}$ is served before $j_{2}$, and so on, until
    the before last patient $j_{n-1}$ is served before $j_n$ by $k$, i.e.,
    $x(j_{1},j_{2},k)=x(j_{2},j_{3},k)=\ldots=x(j_{n-1},j_n,k)=1$
      and by \eqref{eq3} this implies that $z(j_n , k)=1$.
    \item By a similar reasoning, if $i$ is the last call attended
    by ambulance $k$, i.e., $z(i,k)=1$ then
    there is also a first call attended, i.e., there is
    $j$ such that $x(j,k)=1$ and no calls are attended
    by $k$ after $i$, i.e., $x(i,j,k)=0$ by \eqref{eq3}.
    \item Similarly, if $x(i,k)=0$ (resp. $z(i,k)=0$) for all $i$ then 
    the constraints imply $x(i,j,k)=0$ for all $i,j$ and
    $z(i,k)=0$ (resp. $x(i,k)=0$) for all $i$  
    which correctly corresponds to the situation where ambulanc e
    $k$ does not take any call. 
    \end{itemize}
    
    \par {\textbf{Constraints on the maximal
    completion times.}}\\
    
    We will assume we have three types 
    of calls where
    patients of
    calls of type 1 (high priority) and 2 (intermediate
    priority) are transported
    to hospital and patients of calls
    of type 3 (low priority) are not transported
    to hospital (these assumptions
    are satisfied for the real emergency health
    service we consider in our numerical experiments).
    \par For all $i \in C$ with $\mbox{Type}_{c}(i)=1$,
    we must have for all $h \in H$ and $cb \in CB$:
    \begin{equation}\label{eq15}
    \begin{array}{l}
    M_{1} \geq t_{i} + \mbox{TimeOnScene}_{c}(i) + y(i,h) \mbox{travelTime}(\ell_{c}(i),h),\;\forall i \in C_{2},\\
    M_{1} \geq t_{i} + \mbox{TimeOnScene}_{c}(i)
    + y(i,h,cb) \mbox{travelTime}(\ell_{c}(i),h),\;\forall i \in C_{1}.
    \end{array}
    \end{equation}
    
    \par For all $i \in C$ with $\mbox{Type}_{c}(i)=2$,
    we must have for all $h \in H$ and $cb \in CB$:
    \begin{equation}\label{eq16}
    \begin{array}{l}
    M_{2} \geq t_{i} + \mbox{TimeOnScene}_{c}(i) + y(i,h) \mbox{travelTime}(\ell_{c}(i),h),\;\forall i \in C_{2},\\
    M_{2} \geq t_{i} + \mbox{TimeOnScene}_{c}(i)
    + y(i,h,cb) \mbox{travelTime}(\ell_{c}(i),h),\;\forall i \in C_{1}.
    \end{array}
    \end{equation}
    
    \par For all $i \in C$ with $\mbox{Type}_{c}(i)=3$,
    we must have:
    \begin{equation}\label{eq17}
    \begin{array}{l}
    M_{3} \geq 
    t_{i}, \;\forall  i \in C_{3} \cup C_{4}.
    \end{array}
    \end{equation}
    
    The optimal allocation of ambulances to calls in queues
    amounts to minimize
    \begin{equation}\label{objpenalized}
    \sum_{i=1}^3 
     \mbox{penalization}(M_{i}-t_{c}(i),i)
    \end{equation}
    under constraints \eqref{eq1}-\eqref{eq17}
    and $M_{1}, M_{2}, M_{3}, t_{i} \geq 0$ continuous,
    with variables 
    $x(i,k)$, $z(i,k)$, $x(i,j,k)$, $y(i,cb)$, $y(i,h)$, $y(i,h,cb)$
    taking values in $\{0,1\}$, where
    penalization$(t,i)$ is the penalized
    counterpart of $t$ time units for a call of type $i$.
    
    \par Depending on the number of calls in the queue, solving this problem
    quickly could require heuristics as in
     or decomposition techniques such as column generation.
    
    To get a quick solution for queues of small to moderate size,
    we will make the simplifying but reasonable in practice assumption that when a hospital
    is needed the closest hospital to the scene that can attend
    patients of the corresponding type is chosen and 
    when a cleaning station is needed the closest cleaning station is chosen.
    This implies that for a given call $i$ the time $\Delta_{i}$ the ambulance is in
    service and the location $L_{i}$ the ambulance will be when it becomes
    newly available are known in advance and do not depend on hospital
    or cleaning station decision variables. More precisely, 
    \begin{itemize}
    \item if $i \in C_{1}$, $\Delta_{i}$ is the time spent on the scene of the call
    plus the time to go to the closest hospital, the time spent at hospital,
    the time to go to the closest cleaning station, and the time spent
    at the cleaning station while
    $L_{i}$ is the location of the cleaning station which is the closest to the hospital;
    \item if $i \in C_{2}$, $\Delta_{i}$ is the time spent on the scene of the call
    plus the time to go to the closest hospital and the time spent at hospital
    while
    $L_{i}$ is the location of the hospital closest to the scene of the call;
    \item if $i \in C_{3}$,  $\Delta_{i}$ is the time spent on the scene of the call
    plus the time to go to the closest cleaning station plus the time spent
    to clean the ambulance while
    $L_{i}$ is the location of the cleaning station which is the closest to the scene of the call;
    \item if $i \in C_{4}$, $\Delta_{i}$ is the time spent on the scene of the call
    and $L_{i}$ is the scene of the call.
    \end{itemize}
    For $i \in C_{1} \cup C_{2}$ we will also denote by
    $\tau_{i}$ the time spent between the moment the ambulance
    arrives on the scene of the call and the time the ambulance
    arrives at hospital for that call
    while for $i \in C_{3} \cup C_{4}$ we set $\tau_{i}=0$.
    
    With this assumption we need to replace 
    \eqref{eq14}-\eqref{eq11} by
    \begin{equation}\label{eq18}
    t_{i} + \Delta_{i} +   \mbox{travelTime}(L_{i},\ell_{c}(j)) \leq 
    t_{j} + M\Big(1- \sum_{k: k \in \mathcal{A}(\mbox{Type}_{c}(i)) \cap \mathcal{A}(\mbox{Type}_{c}(j))}x(i,j,k)\Big),
    \; \forall i \neq j \in C
    \end{equation}
    and \eqref{eq15}-\eqref{eq17} by
    \begin{equation}\label{eq19}
    \begin{array}{l}
    M_{j} \geq 
    t_{i}+\tau_{i}, \;\forall i \mbox{ with Type}_{c}(i)=j,j=1, 2, 3.
    \end{array}
    \end{equation}

    Therefore the corresponding allocation problem amounts to minimizing 
    \eqref{objpenalized} under the constraints
    \eqref{eq1}-\eqref{eq3}, \eqref{eq8}-\eqref{eq10}, \eqref{eq18}, \eqref{eq19}
    and $M_{1}, M_{2}, M_{3}, t_{i} \geq 0$,
    with variables 
    $x(i,k)$, $z(i,k)$, $x(i,j,k)$ taking values in $\{0,1\}$ 
    (note that variables $y(i,cb)$, $y(i,h)$, and $y(i,h,cb)$
    are not necessary anymore).

}\fi

\section{Choice of a new ambulance station after service}
\label{sec:newbase}

The heuristics presented in the previous sections determine what to do when a new emergency call arrives and what emergency in queue to send an ambulance to when an ambulance becomes available.
When an ambulance finishes service, either on the scene of an emergency or at a hospital, and there are no emergencies in queue, then a decision has to be made where (to what ambulance station) to send the ambulance for staging until it is dispatched again.
We consider the following three rules in numerical experiments:
\begin{itemize}
\item
Home Station Rule (HSR): Each ambulance is assigned a ``home station'', and when an ambulance finishes service and there are no emergencies in queue, then it always returns to its home station.
\item
Closest Station Rule (CSR): When an ambulance finishes service and there are no emergencies in queue, then it goes to the ambulance station that is closest to its current location.
\item
Best Station Rule (BSR): When an ambulance finishes service and there are no emergencies in queue, then it goes to the ambulance station in most urgent need of an additional ambulance, as described below.
\end{itemize}

Next we describe the details of BSR.
Consider an ambulance~$j$ that finishes service at location $\ell_{0}$ at time~$t_{0}$.
For each station~$b$, let $t(j,b)$ denote the time at which ambulance~$j$ would arrive to station~$b$ if it was sent from $\ell_{0}$ to $b$ at $t_{0}$. 
Let $\Delta > 0$ denote the length of a lookahead time interval.
Let $\mathcal{C}$ denote the list of emergency types sorted in decreasing order of priority, that is, emergency type $\mathcal{C}(1)$ has higher priority than emergency type $\mathcal{C}(2)$ and so on.
Each station~$b$ has a service region for which it is planned to cover emergencies.
For station~$b$, emergency type~$c$, and times $t_{1} < t_{2}$, let $q(b,c,t_{1},t_{2})$ denote the planned number of emergencies of type~$c$ served from station~$b$ during time interval $[t_{1},t_{2}]$.
For example, $q(b,c,t_{1},t_{2})$ could be an $\alpha$-quantile (say with $\alpha = 0.9$) of the distribution of the number of emergencies of type~$c$ in the service region of station~$b$ during time interval $[t_{1},t_{2}]$.
If, as in \cite{laspatedpaper} and \cite{laspatedmanual}, emergencies of type~$c$ in the service region of station~$b$ arrive according to a Poisson process with intensity $\lambda_{b,c}$, then $q(b,c,t_{1},t_{2})$ could be the $\alpha$-quantile of the Poisson distribution with mean $\lambda_{b,c} (t_{2} - t_{1})$.
For station~$b$, emergency type~$c$, and time~$t$, let $A(b,c,t)$ denote the number of ambulances (without counting ambulance~$j$) suited to emergencies of type~$c$ that are forecast to be at station~$b$ at time~$t$; this includes ambulances already at station~$b$ at time~$t_{0}$ not forecast to serve another call and ambulances in service forecast to arrive at station~$b$ by time~$t$.
Then, the ambulance deficit at station~$b$ for emergency type~$c$ is given by
\begin{equation}
\label{defamb}
D(b,c) \ \ = \ \ q(b,c,t(j,b),t(j,b)+\Delta) - A(b,c,t(j,b))
\end{equation}

The steps of BSR rule are as follows. \\
\textbf{Step 0.}
Initialize $\ell = 1$. \\
\textbf{Step 1.}
Let $m(\ell) \defi \max\{D(b,\mathcal{C}(\ell)) \, : \, b \in \mathcal{B}\}$ and $b(\ell) \defi \arg\max\{D(b,\mathcal{C}(\ell)) \, : \, b \in \mathcal{B}\}$.
\textbf{Step 2.}
If $m_{\ell} > 0$ then send ambulance~$j$ to station $b(\ell)$, else if $\ell = |\mathcal{C}|$ then send ambulance~$j$ to station $b\big(\arg\max\{m(\ell') \, : \, \ell' \in \{1,\ldots,|\mathcal{C}|\}\}\big)$, otherwise increment $\ell \leftarrow \ell + 1$ and go to Step~1.

\if{
 Additionally, let \(N_{abc}(t)\) be the number of ambulances of type \(a\) currently stationed at station \(b\) at the time \(t\) that can respond to a call of priority \(c\).

Let \(j\) be an ambulance that finished its service at currentTime. Given the highest priority \(c\) that \(j\) can respond to, we compute, for each station \(b\), the arrival time \(t_b(j)\) of \(j\) at station \(b\) and the ambulance deficit \(\Delta_{bc}\) of \(b\) for calls of priority \(c\). Specifically, \(\Delta_{bc} = Q(\bar{\lam}) - N_{abc}(t_{b})\), where \(Q(\bar{\lam})\) is the 0.9-quantile of the Poisson random variable with rate \(\bar{\lam}\) and \(\bar{\lam} = \sum\limits_{t=t_b}^{t_b + t_{ahead}} t\lam_{btc}\) is the arrival rate in the interval \([t_{b}, t_{b} + t_{ahead}]\). Note that we discretize the periods, so the interval \([t_{b}, t_{b} + t_{ahead}]\) is actually the union of \([t_{b},t_{0}]\), \([t_{0},t_{1}]\) ... \([t_{n-1},t_{n}]\), \([t_{n}, t_{b} + t_{ahead}]\), where \(t_{0}\) is the beginning of the first period greater than \(t_{b}\) and \(t_n\) is the beginning of the last period smaller than \(t_{b} + t_{ahead}\).

Next, we compute the station \(b^{*}\) that achieves maximum \(\Delta_{bc}\) and is closest to the current location of \(j\). If \(\Delta_{b^{*}c} \leq 0\), there is no deficit of calls of priority \(c\), so we compute \(\Delta_{b(c+1)}\) for the next (less advanced) priority type \(c+1\). If \(\Delta_{b^{*}c} > 0\), then we reassign \(j\) to \(b^{*}\).

}\fi

\section{Numerical results}
\label{sec:numsec}

In this section we present numerical results based on the data of Rio de Janeiro EMS, as well as numerical results based on a mixture of this data with artificial data.
The methods were implemented in C++17 and Matlab, and the source code is available at \url{https://github.com/vguigues/Heuristics_Dynamic_Ambulance_Management}.
The experiments were performed on a computer with AMD Ryzen 5 2600 3.4 GHz CPU, 24GB of RAM, in a Ubuntu 22.04 OS.

All instances have $2$~types of ambulances: ALS and BLS, and 4 types of emergencies, indexed by 1, 2, 3, and 4.
Any ambulance can serve any emergency, but different emergency types have different priorities and different ambulance preference, as follows:
\begin{itemize}
\item
type 1 emergency:  high priority call that should preferably be served by an ALS ambulance;
\item
type 2 emergency: low priority call that should preferably be served by an ALS ambulance;
\item
type 3 emergency: high priority call that can be served equally well by any ambulance;
\item
type 4 emergency: low priority call that can be served equally well by any ambulance.
\end{itemize}
Correspondingly, the allocation cost is given by \eqref{costalloc} with quality of care coefficients $M_{ac}$ given in Table~\ref{tbl:compatibility_matrix}, as well as function {\tt penalization} of the waiting time given by \eqref{penfunc}  with $\theta_{c} = 4$ for emergency types 1 and 3 (high priority) and $\theta_{c} = 1$ for emergency types 2 and 4 (low priority).

\begin{table}[ht!]
\centering
\begin{tabular}{|c|c|c|c|c|}
\hline
   & 1:  High, ALS pref.  & 2: Low, ALS pref. & 3: High, no pref.  & 4: Low, no pref.\\ \hline
ALS & 0 & 0 & 1500 & 1500 \\ \hline
BLS & 6000 & 6000 & 0 & 0 \\ \hline
\end{tabular}
\caption{Quality of care coefficients $M_{ac}$.
The rows correspond to the ambulance types and the columns correspond to the emergency types.}
\label{tbl:compatibility_matrix}
\end{table}

\subsection{Experiments with real data}

We applied our methods to the management of Rio de Janeiro EMS, which serves several million people.
The study was done in collaboration with Rio de Janeiro EMS, and we used the history of emergencies for the period January 2016--February 2018, as well as the locations of ambulance stations and hospitals, and the set of ambulances.
The data of emergencies is available at \url{http://samu.vincentgyg.com/} (confidential data is omitted, and the available data contain the location of each emergency, the time each emergency call was received, and the type of each emergency).

\begin{figure}
\centering
\includegraphics[scale=0.6]{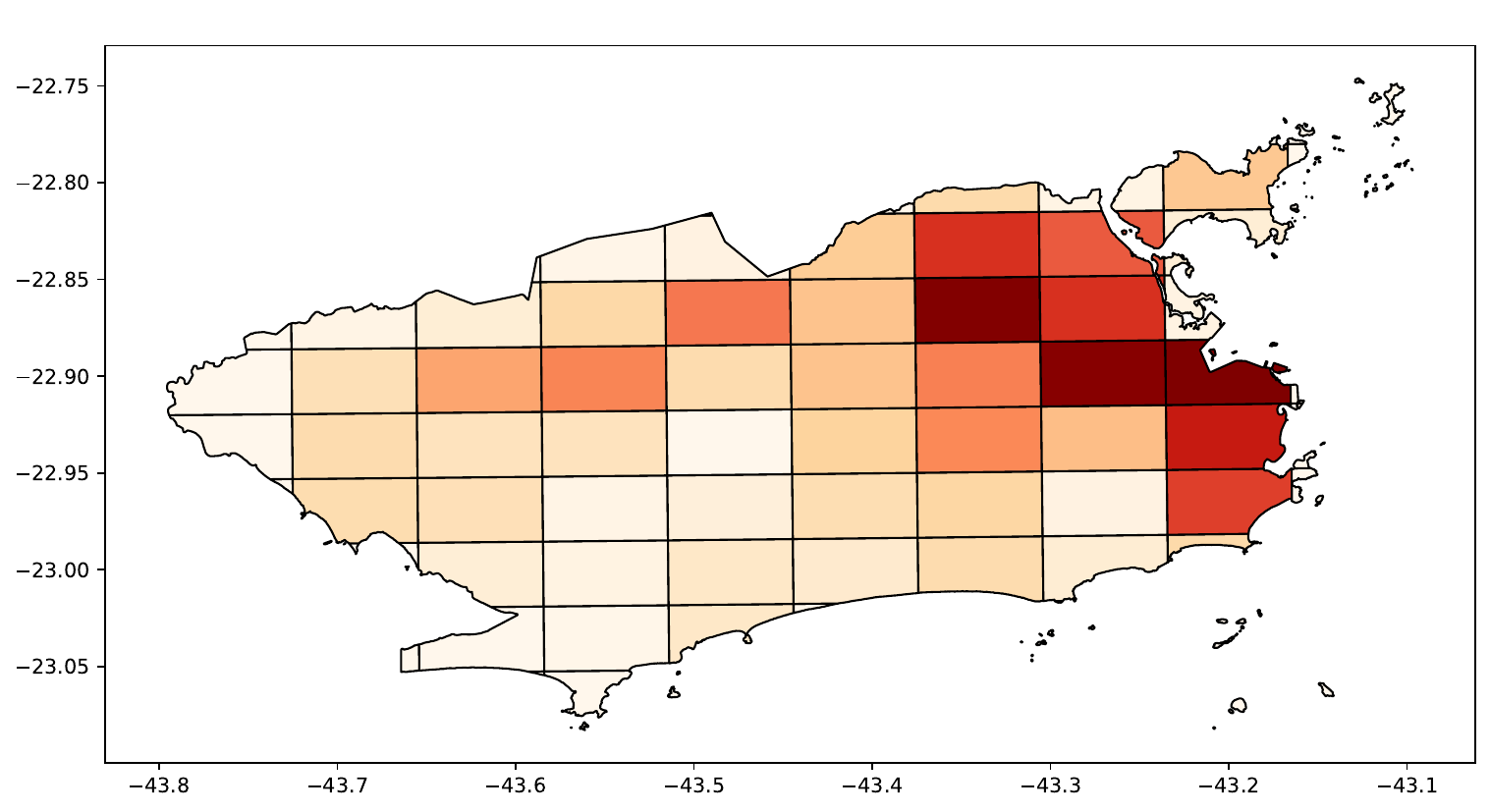}
\caption{Discretization of the service region of  Rio de Janeiro SAMU into $10 \times 10$ rectangles, and a heatmap of the intensity of emergency call rates aggregated by emergency type and time period for the period January 2016--February 2018.\label{figuresheat}}
\end{figure}
\begin{figure}
\centering
\resizebox{\textwidth}{!}{
\begin{tabular}{c}
\includegraphics[scale=0.025]{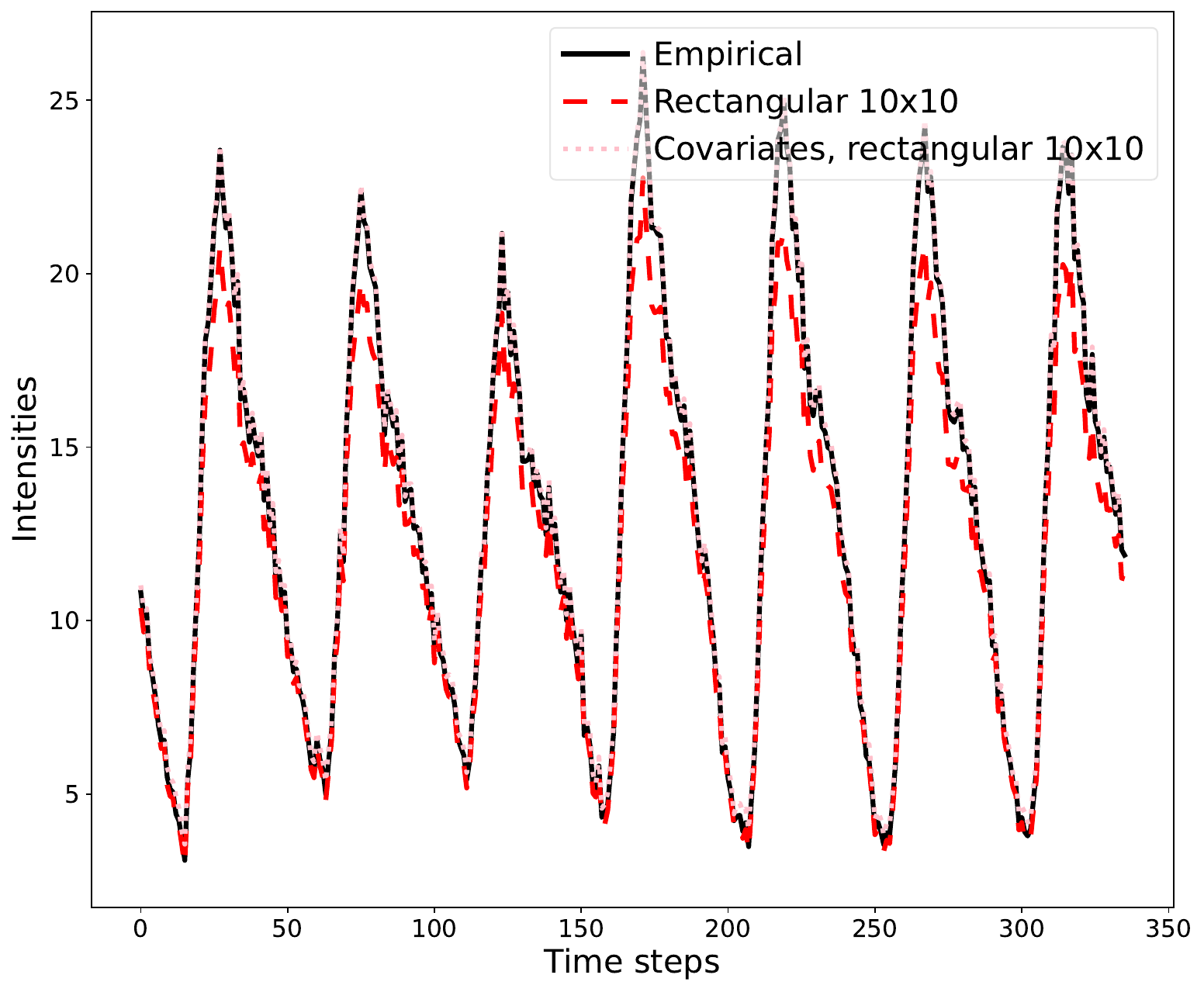}
\end{tabular}}
\caption{Aggregated estimates of the intensities over the week for 3 estimators: (i)~Empirical: the empirical intensities, (ii)~Rectangular $10 \times 10$: regularized model with the $10 \times 10$ rectangular space discretization, and (iii)~Covariates, rectangular $10 \times 10$: model with covariates with the $10 \times 10$ rectangular space discretization.
\label{fig:allstepsrect}}
\end{figure}

A spatio-temporal Poisson process was fitted to the emergency data.
For the calibration of the process, the city was discretized into $10 \times 10$ rectangular zones (see \citealt{laspatedpaper} and \citealt{laspatedmanual} for details of this procedure).
The fitted intensities were periodic with a period of one week, and the weekly period was discretized into time intervals of 30 minutes each.
Thus, the estimated Poisson intensities $\lambda_{c,i,t}$ are indexed by emergency type~$c$, zone~$i$, and time interval~$t$.
A heatmap of emergency call rates $\lambda_{c,i,t}$ aggregated by emergency type~$c$ and time period~$t$ is shown in Figure~\ref{figuresheat}.
Also, emergency call rates $\lambda_{c,i,t}$ aggregated by emergency type~$c$ and zone~$i$ for 3 estimators are shown in Figure~\ref{fig:allstepsrect}, see \citealt{laspatedpaper} and \citealt{laspatedmanual} for details on the different estimators.
For the numerical results in this paper, we use the regularized model from \cite{laspatedpaper} with weights obtained using cross validation.

We compare the $4$ ambulance dispatch heuristics described above (BM from Section~\ref{bestmyopic}, NM from Section~\ref{nonmyopic}, and GHP1 and GHP2 from Section~\ref{sec:prior}) with \textcolor{blue}{$7$} heuristics from the literature (\citealt{ande:07,lees:12,mayo:13,band:14,lees:14,jagt:15,lees:17}), as well as the Closest-Available (CA) heuristic.
Each of these heuristics is used both as a rollout policy~$H$, as described in Section~\ref{smodel}, as well as to make dispatch decisions directly without rollout.
For each of these heuristics, when an ambulance has to be sent to a cleaning station, it is sent to the closest cleaning station.
The simulation of the heuristics is done on a set of 100 scenarios and for each rollout heuristic, we use, each time a decision has to be made, a scenario tree of 100 future scenarios.

\par {\textbf{Deterministic travel times.}}
We first consider deterministic travel times,
assuming ambulances travel at mean speed 60 km/h.

Tables~\ref{costallocdetnoroll}
and~\ref{table_response_times_deterministic} show the allocation costs and response times for each of the heuristics applied directly without rollout.
Tables~\ref{costallocdetroll} and~\ref{table_response_times_deterministic_rollout} show the allocation costs and response times for each of the heuristics applied as rollout policies~$H$, as explained in Section~\ref{smodel}.

In these tables, the columns headed ``Mean'' show the mean performance metric per emergency, the columns headed ``Q0.9'' show the empirical $0.9$~quantile of the performance metric, and the columns headed ``Max'' show the maximum value of the performance metric per emergency over the $100$ sample paths.

In these results, all heuristics used the closest station rule (CSR) when sending an ambulance to an ambulance station.
The results are based on $100 $ independent sample paths for the peak period of Friday 6pm to 8pm.
Results are reported for settings with $10$, $12$, $14$, $18$, $22$, $24$, $26$ and $30$ ambulances.

\begin{figure}
    \centering
    \includegraphics[width=0.8\linewidth]{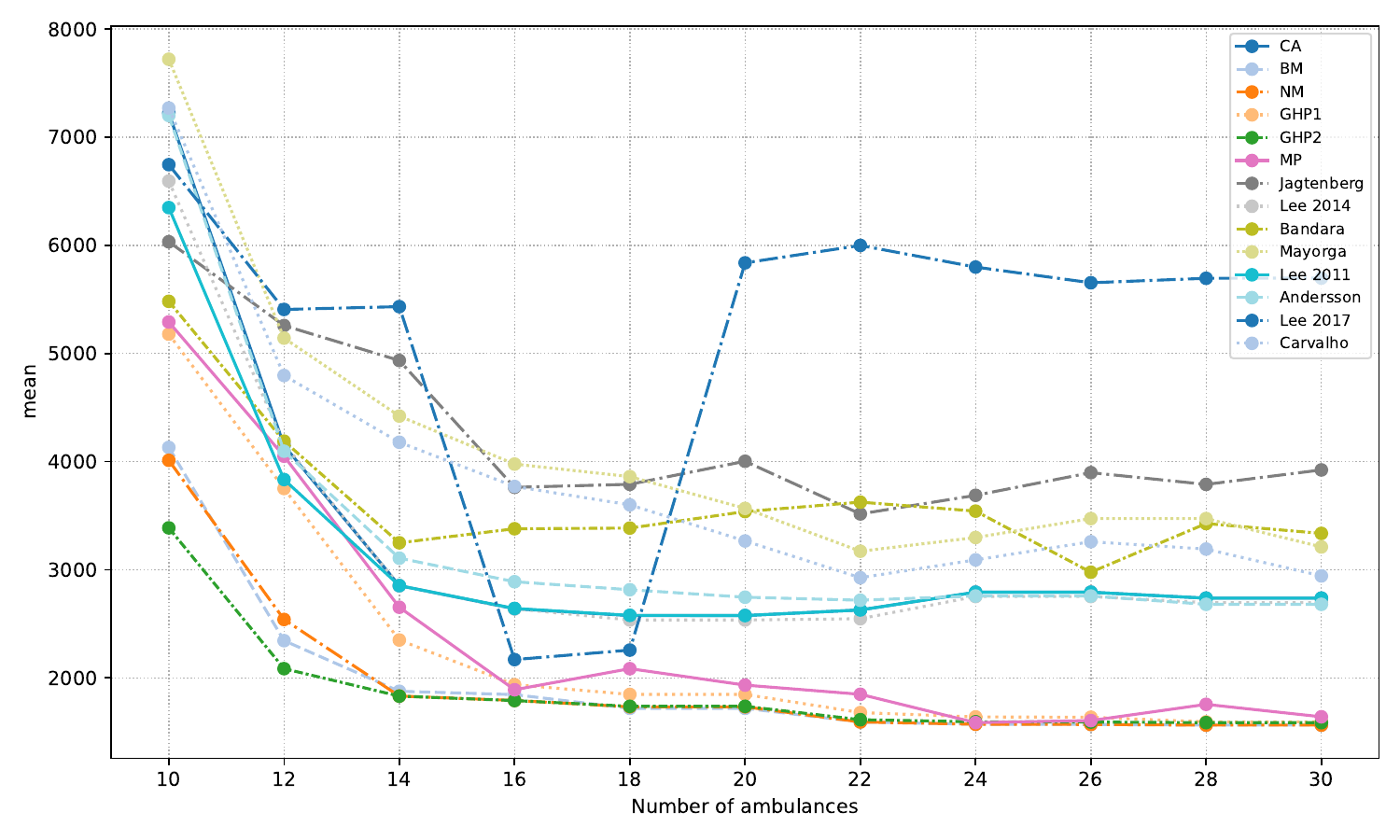}
    \caption{Mean allocation costs for each heuristic.}
    \label{fig:mean_alloc_costs_deterministic}
\end{figure}

\begin{figure}
    \centering
    \includegraphics[width=0.8\linewidth]{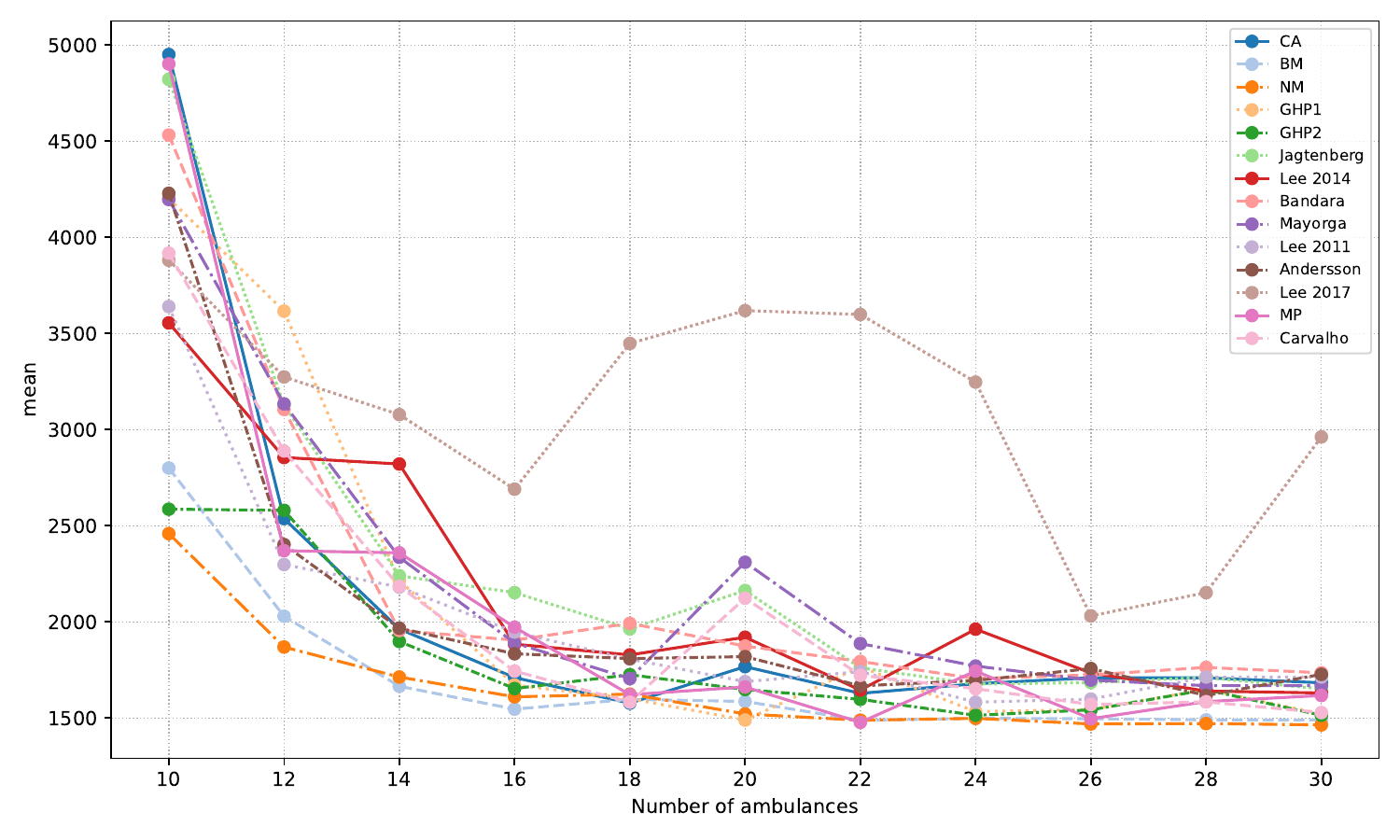}
    \caption{Mean allocation costs for each rollout heuristic.}
    \label{fig:mean_alloc_costs_deterministic_rollout}
\end{figure}

To ease the comparisons of the methods
in terms of allocation costs and response times, we plot in Figure \ref{fig:mean_alloc_costs_deterministic}
(resp. Figure \ref{fig:mean_alloc_costs_deterministic_rollout})
the mean allocation for each heuristic
without rollout (resp. with rollout)
as a function of the number of ambulances.

\begin{figure}
    \centering
    \includegraphics[width=0.8\linewidth]{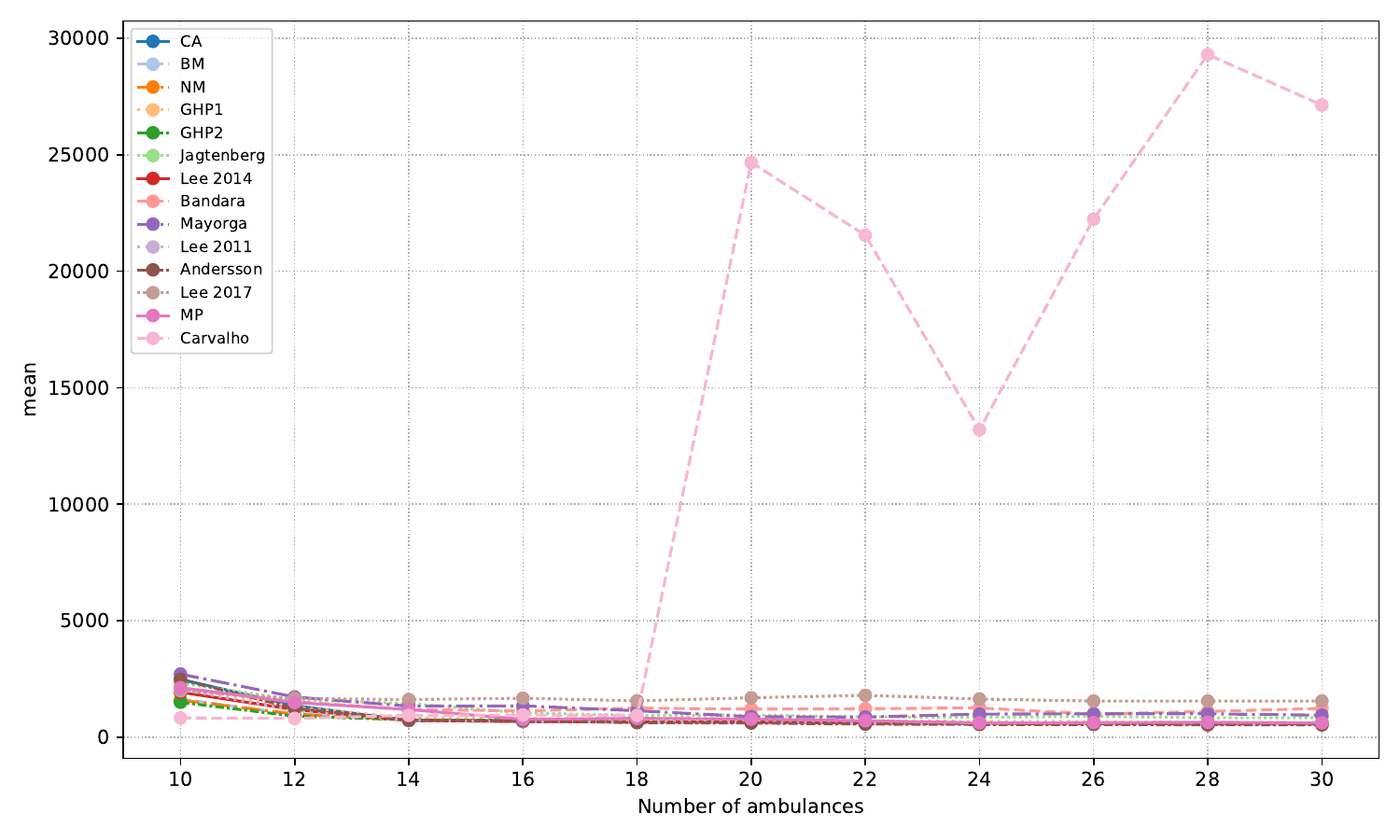}
    \caption{Mean response times for each heuristic.}
\label{fig:mean_response_times_deterministic}
\end{figure}

\begin{figure}
    \centering
    \includegraphics[width=0.8\linewidth]{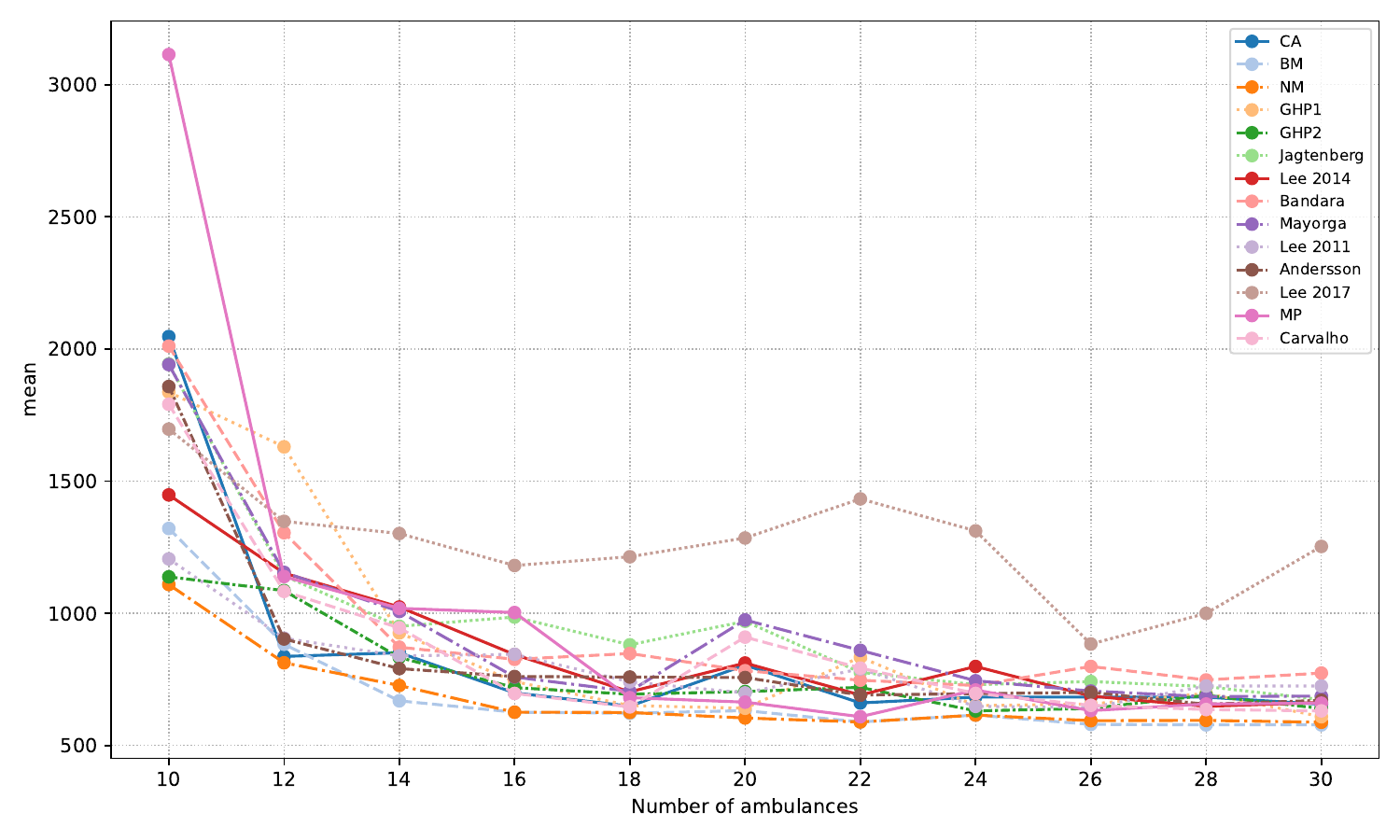}
    \caption{Mean response times for each rollout heuristic.}    \label{fig:mean_response_times_deterministic_rollout}
\end{figure}

We also plot in Figure \ref{fig:mean_response_times_deterministic}
(resp. Figure \ref{fig:mean_response_times_deterministic_rollout})
the mean response time for each heuristic
without rollout (resp. with rollout)
as a function of the number of ambulances.

\begin{figure}[htbp]
    \centering
    \includegraphics[width=0.45\textwidth]{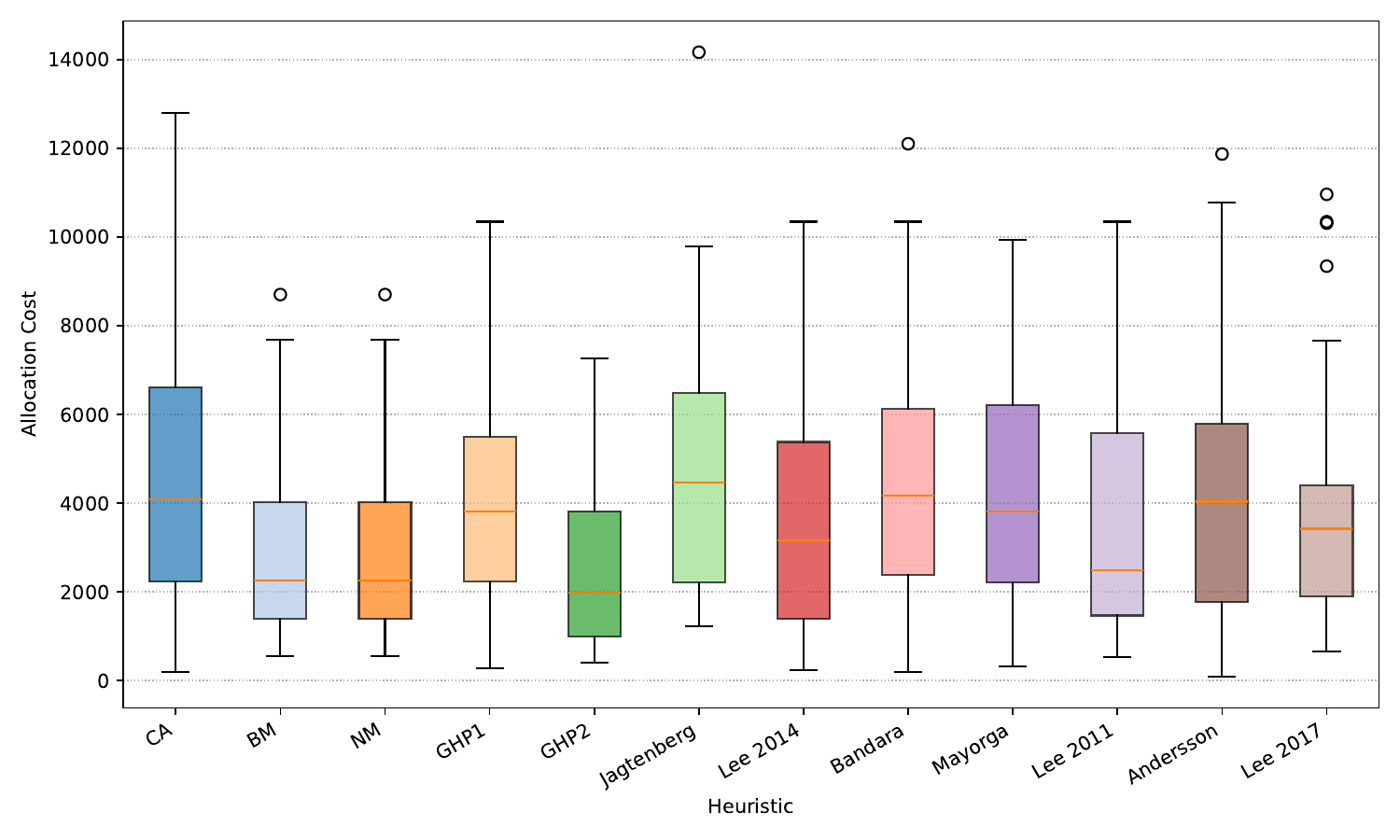}
    \includegraphics[width=0.45\textwidth]{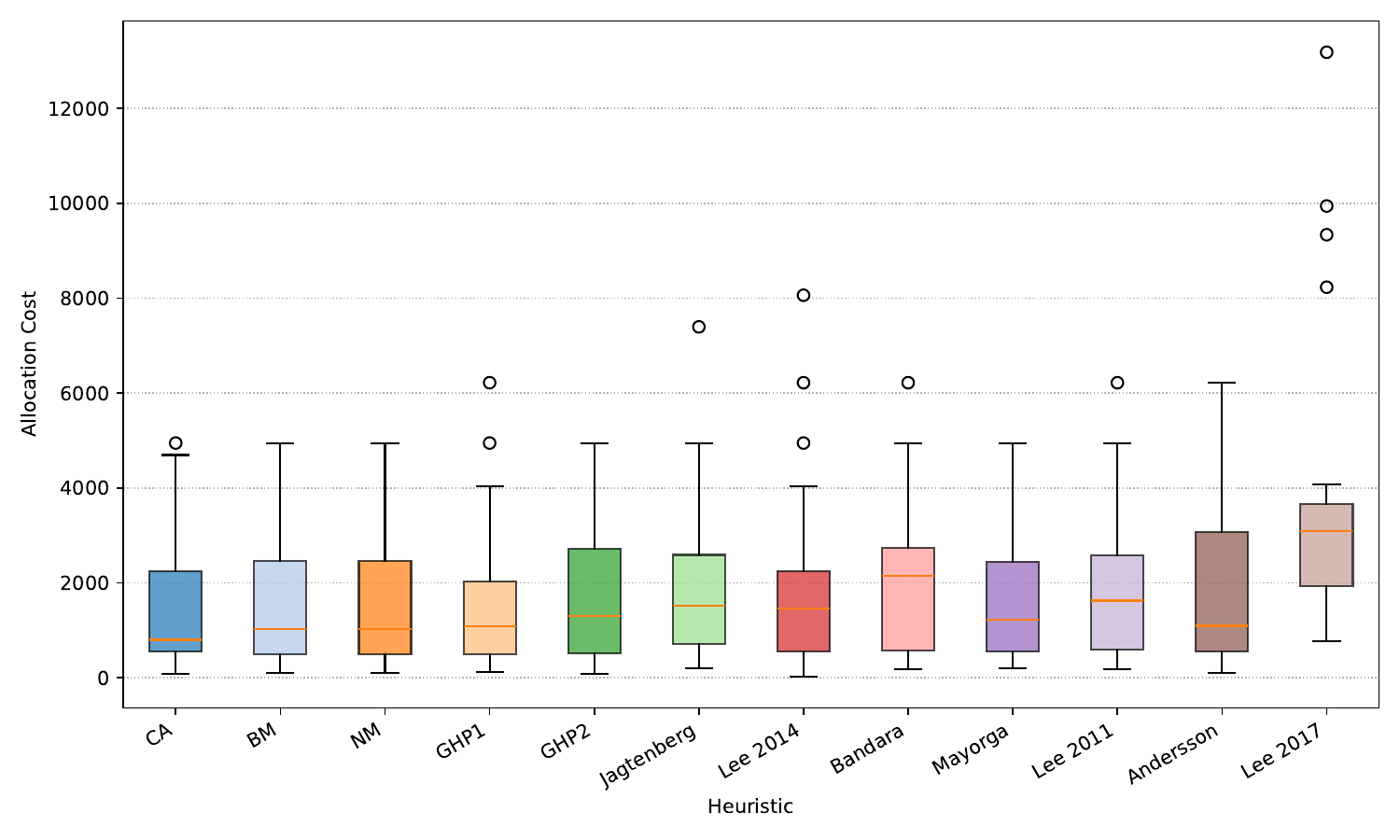}

    \vspace{0.3cm}

    \includegraphics[width=0.45\textwidth]{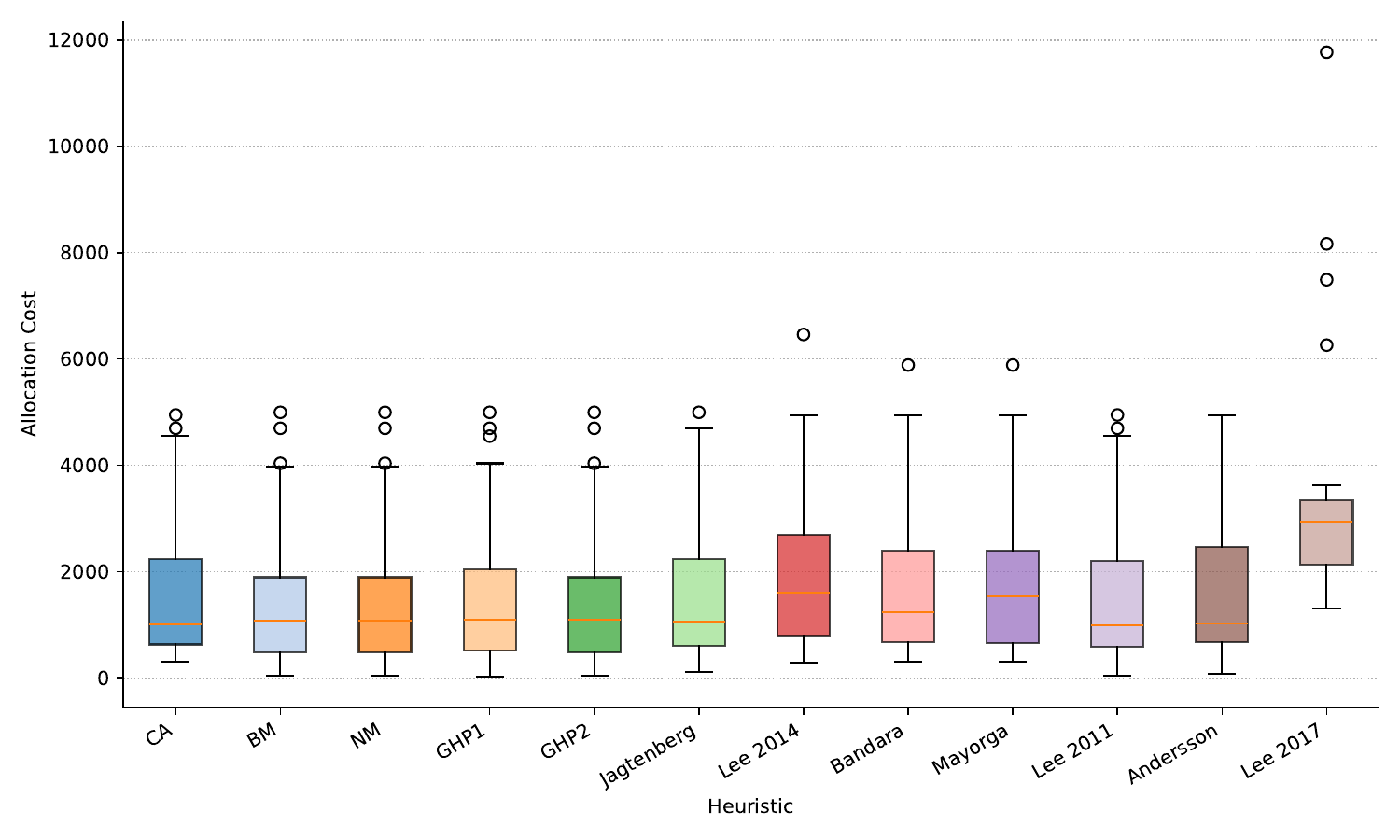}
    \includegraphics[width=0.45\textwidth]{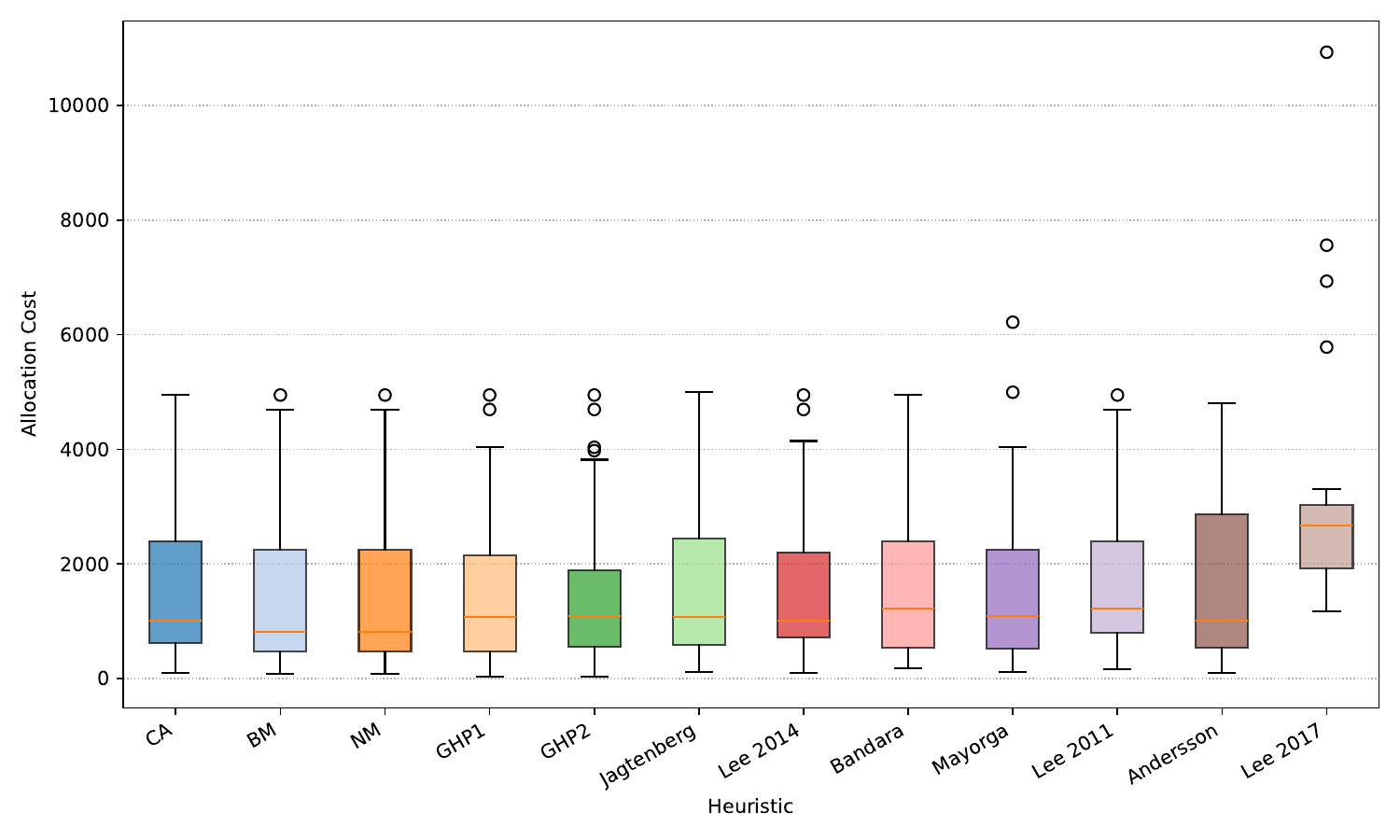}

    \caption{Box plots for allocation costs of the heuristics considering — from top left to bottom right — 10, 18, 24 and 30 ambulances.}
    \label{boxplots_allocation_costs}
\end{figure}

\begin{figure}[htbp]
    \centering
    \includegraphics[width=0.45\textwidth]{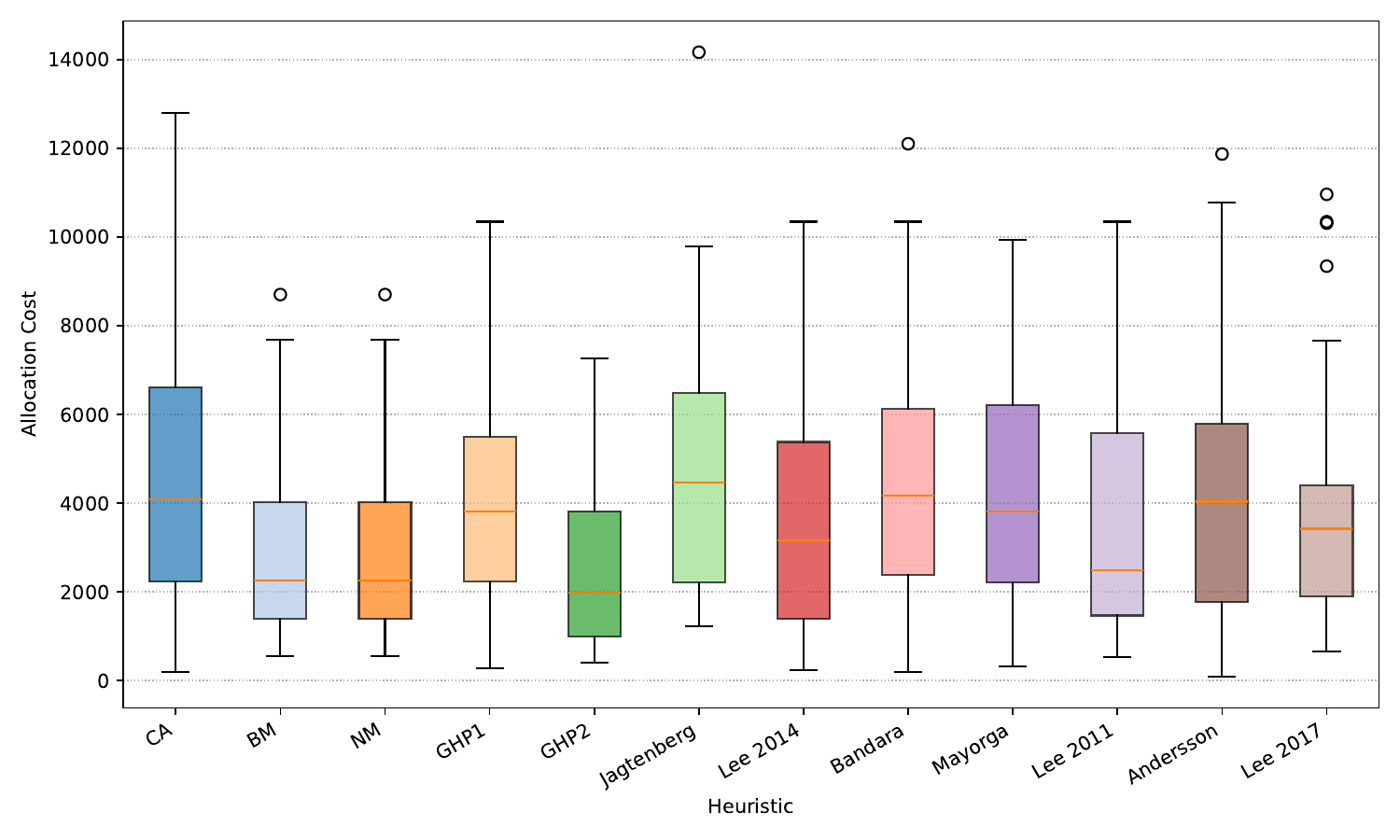}
    \includegraphics[width=0.45\textwidth]{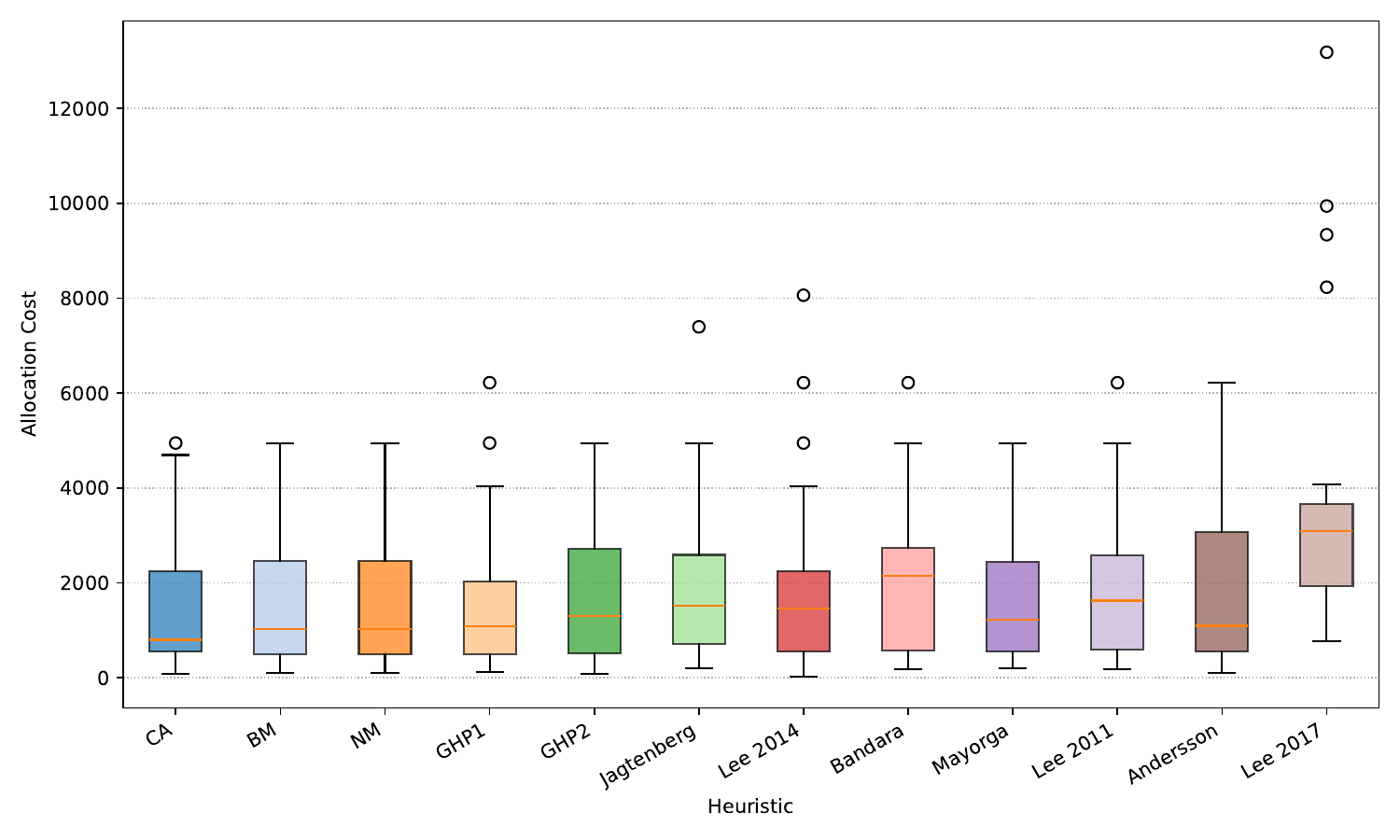}

    \vspace{0.3cm}

    \includegraphics[width=0.45\textwidth]{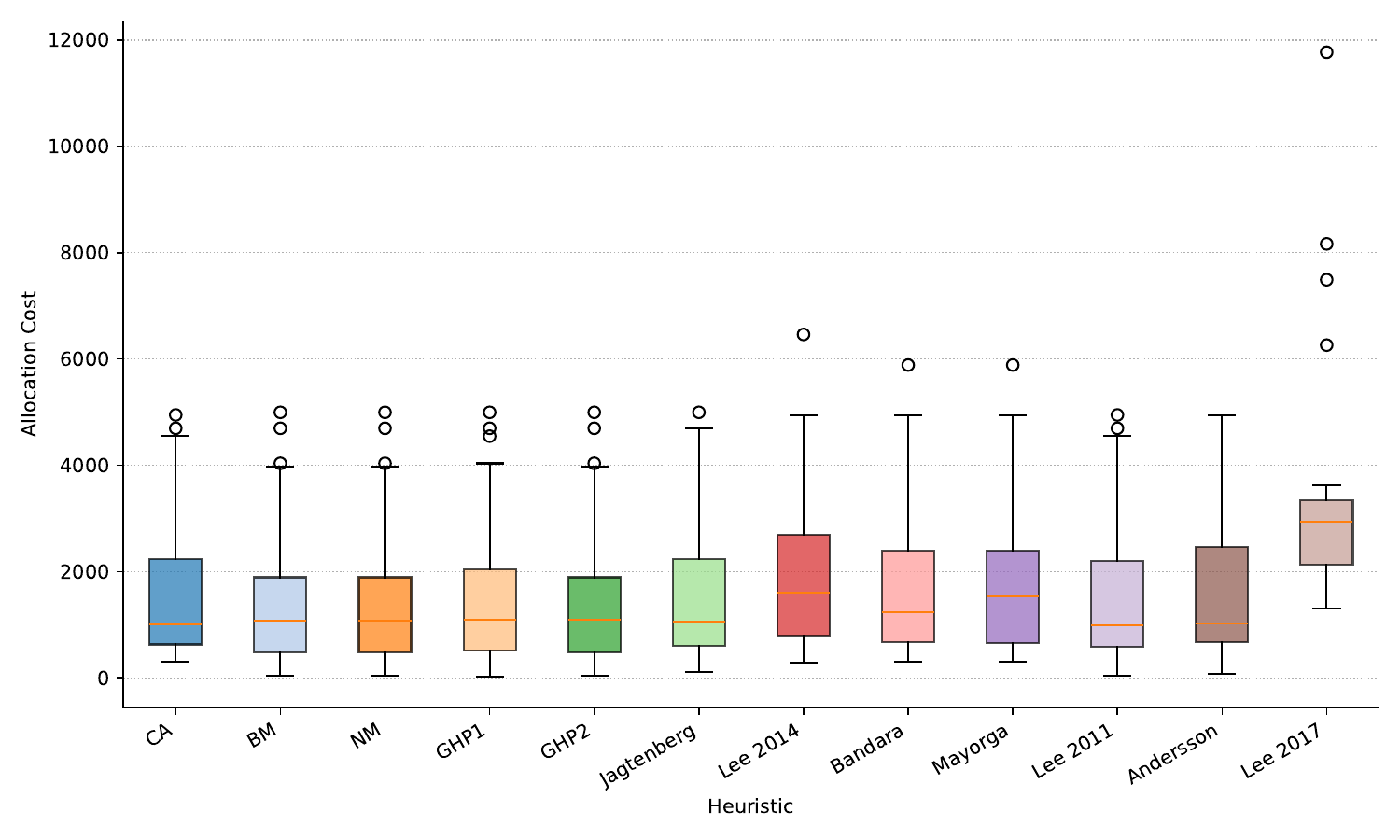}
    \includegraphics[width=0.45\textwidth]{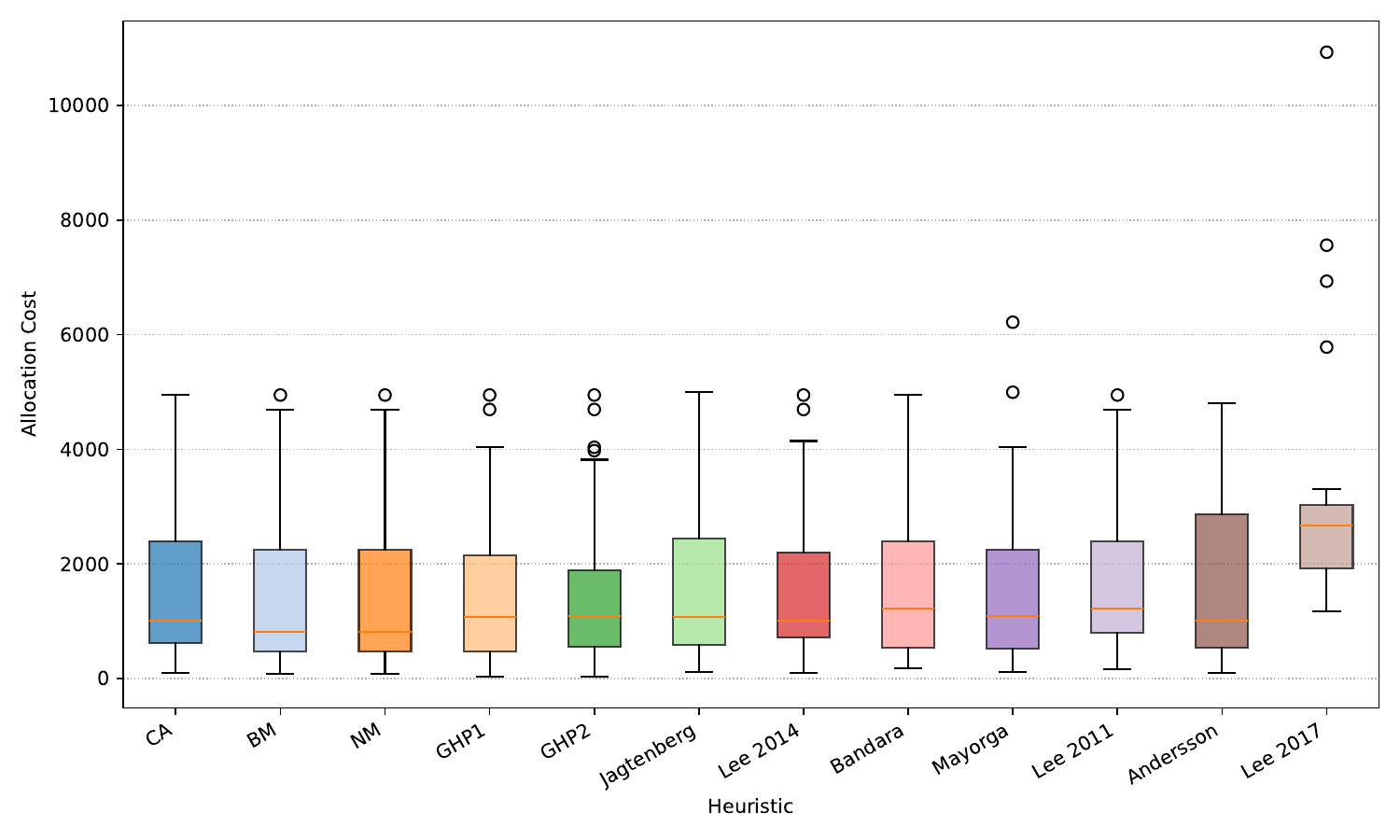}

    \caption{Box plots for allocation costs of the rollout heuristics considering — from top left to bottom right — 10, 18, 24 and 30 ambulances.}
    \label{boxplots_allocation_costs_rollout}
\end{figure}

Figure~\ref{boxplots_allocation_costs} displays the box plots  for the allocation costs of the 12 heuristics considering 10, 18, 24, and 30 ambulances. Figure~\ref{boxplots_allocation_costs_rollout} displays box plots of the allocation costs for the rollout version of the 12 heuristics.

\begin{figure}
    \centering
    \includegraphics[width=0.75\linewidth]{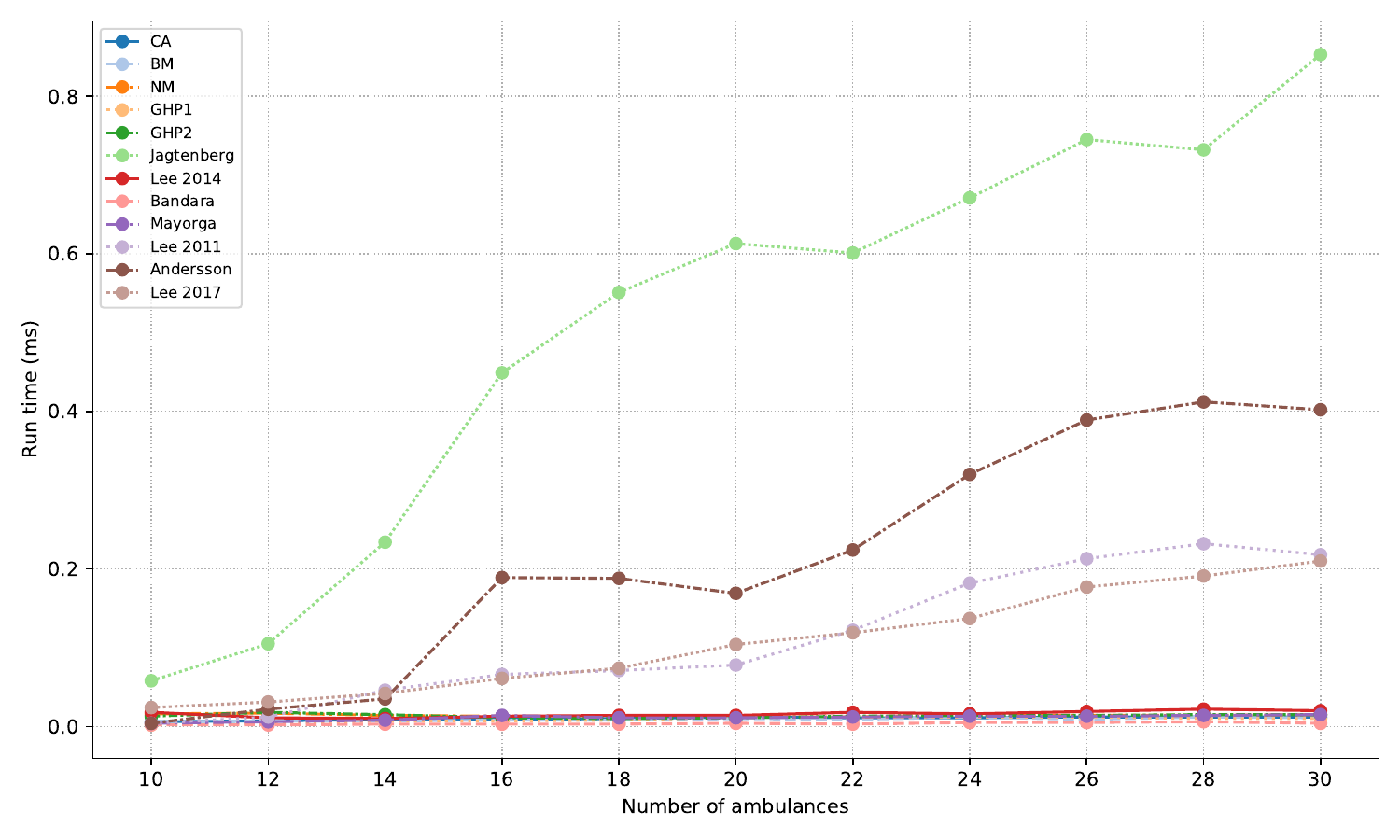}
    \caption{Mean run times for the heuristics when considering ambulance selection.}
    \label{fig:run_times_selection}
\end{figure}

\begin{figure}
    \centering
    \includegraphics[width=0.75\linewidth]{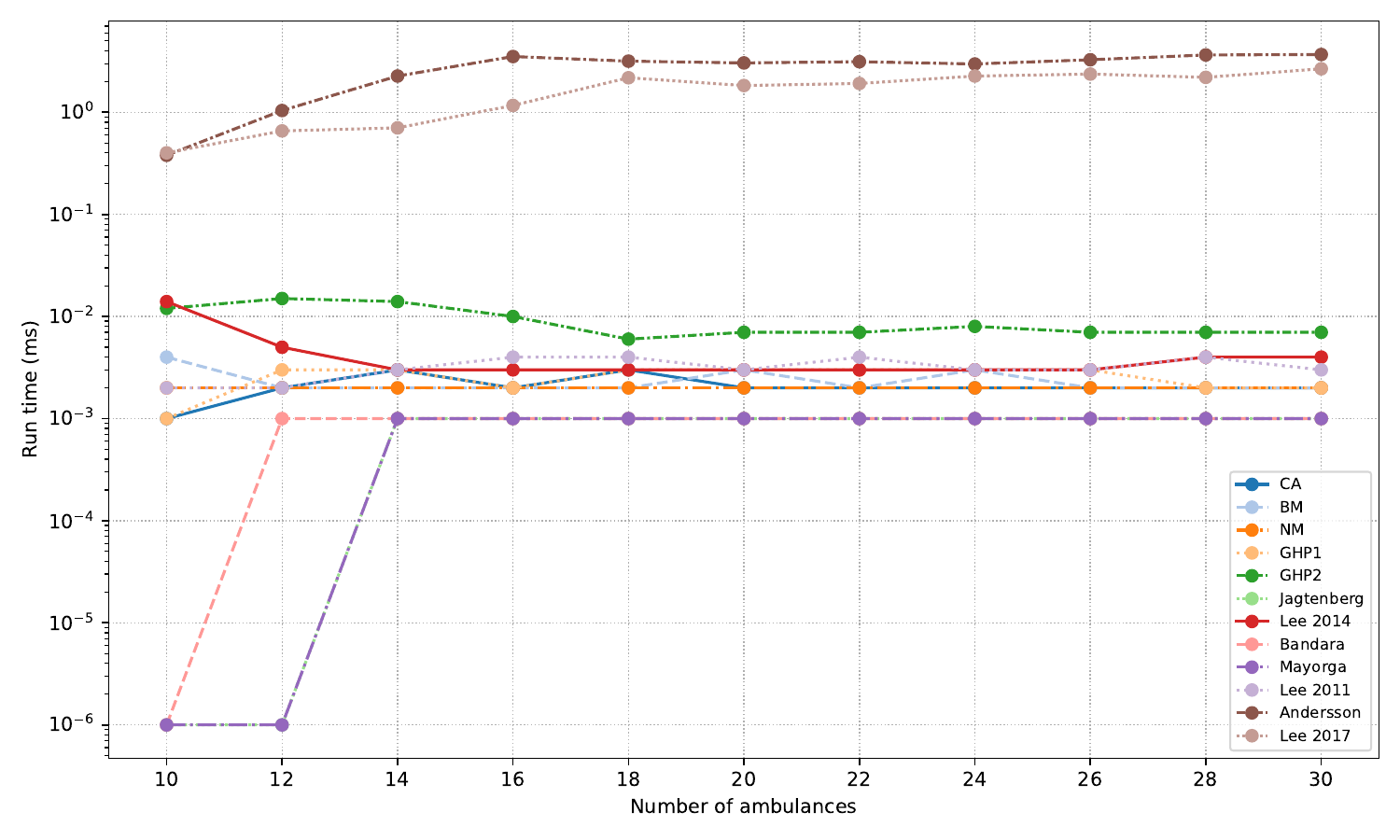}
    \caption{Mean run times for the heuristics when considering ambulance reassignment.}
    \label{fig:run_times_reassignment}
\end{figure}

Figures~\ref{fig:run_times_selection} and \ref{fig:run_times_reassignment} display the mean run times for each heuristic, in ms.
It can be seen that all decisions for all heuristics are computed extremely quickly, making them applicable in practice (since for this application, decisions have to be computed very quickly).

In all instances, the heuristics proposed in this paper outperform the other heuristics  from the literature in terms of allocation costs, and also in terms of response times. The heuristics of \citet{ande:07} and \citet{lees:11} are also competitive with a large number of ambulances.
Furthermore, for each heuristic, performance is better when using the rollout approach than when using the heuristics directly to make decisions.
The improved performance of the rollout approach is obtained at the cost of more computational effort than using the heuristics directly. 
Also, as expected, the performance of all heuristics improves as the number of ambulances increases. 
Such results can be used by an EMS to choose an ambulance fleet size and 
dispatch policy.

\if{
Recall, however, that trajectories follow geodesics.
It would be straightforward to use our code replacing these trajectories
by trajectories along the streets of the city, as long as we have
a map of these streets. In this context, the only functions that
need to be specialized to use our code are the functions travelTime and
positionBetweenOriginDestination}\fi

\begin{table}[ht!]
\scalebox{0.7}{
\centering
\begin{tabular}{@{}l|ccc|ccc|ccc|ccc@{}}
\toprule
Number of Ambulances & \multicolumn{3}{c|}{10} & \multicolumn{3}{c|}{12} & \multicolumn{3}{c|}{14} & \multicolumn{3}{c}{18} \\
\midrule
Heuristic & Mean & Q0.9 & Max & Mean & Q0.9 & Max & Mean & Q0.9 & Max & Mean & Q0.9 & Max \\
\midrule
CA & 7224 & 14628 & 18972 & 4149 & 7417 & 20721 & 2853 & 6524 & 7663 & 2576 & 4796 & 7627  \\
BM & 4131 & 8523 & 13389 & 2344 & 4796 & 9548 & 1876 & 4203 & 5532 & 1718 & 3919 & 4947  \\
NM & 4012 & 8145 & 10383 & 2540 & 5536 & 8828 & 1832 & 4030 & 4808 & 1732 & 4030 & 4808  \\
GHP1 & 5179 & 8852 & 19758 & 3751 & 7936 & 20721 & 2350 & 4796 & 7663 & 1847 & 4197 & 4947  \\
GHP2 & 3387 & 6630 & 7753 & 2085 & 4521 & 5930 & 1831 & 4030 & 4808 & 1737 & 4030 & 4808  \\
\cite{jagt:17a} & 6033 & 13154 & 22547 & 5258 & 10864 & 22547 & 4935 & 9225 & 22547 & 3790 & 8589 & 18710  \\
\cite{lees:14} & 6594 & 15603 & 20975 & 4165 & 7417 & 20721 & 2853 & 6524 & 7663 & 2535 & 4796 & 7627  \\
\cite{band:14} & 5482 & 10698 & 19758 & 4189 & 7967 & 17655 & 3249 & 6136 & 9843 & 3385 & 7008 & 9843  \\
\cite{mayo:13} & 7721 & 16096 & 25422 & 5142 & 8805 & 19661 & 4421 & 7995 & 19661 & 3861 & 7097 & 10414  \\
\cite{lees:11} & 6349 & 12682 & 23215 & 3834 & 7417 & 20721 & 2853 & 6524 & 7663 & 2578 & 4796 & 7627  \\
\cite{ande:07} & 7201 & 14554 & 18972 & 4096 & 7773 & 19645 & 3107 & 6687 & 7861 & 2815 & 5549 & 8446  \\
\cite{lees:17} & 6746 & 13376 & 21435 & 5407 & 10560 & 21020 & 5433 & 11640 & 19879 & 5474 & 11229 & 20930  \\
\bottomrule
\end{tabular}
}
\scalebox{0.7}{
\begin{tabular}{@{}l|ccc|ccc|ccc|ccc@{}}
\toprule
Number of Ambulances & \multicolumn{3}{c|}{22} & \multicolumn{3}{c|}{24} & \multicolumn{3}{c|}{26} & \multicolumn{3}{c}{30} \\
\midrule
Heuristic & Mean & Q0.9 & Max & Mean & Q0.9 & Max & Mean & Q0.9 & Max & Mean & Q0.9 & Max \\
\midrule
CA & 2628 & 6114 & 10414 & 2792 & 6572 & 10414 & 2792 & 6572 & 10414 & 2737 & 6572 & 10414  \\
BM & 1592 & 3906 & 4808 & 1571 & 3906 & 4808 & 1569 & 3906 & 4808 & 1563 & 3906 & 4808  \\
NM & 1592 & 3906 & 4808 & 1571 & 3906 & 4808 & 1569 & 3906 & 4808 & 1563 & 3906 & 4808  \\
GHP1 & 1679 & 4186 & 4947 & 1639 & 4186 & 4947 & 1636 & 4186 & 4947 & 1593 & 4186 & 4947  \\
GHP2 & 1613 & 3906 & 4808 & 1594 & 3906 & 4808 & 1591 & 3906 & 4808 & 1586 & 3906 & 4808  \\
\cite{jagt:17a} & 3516 & 8773 & 12442 & 3688 & 8773 & 12442 & 3898 & 8773 & 12442 & 3923 & 8773 & 12442  \\
\cite{lees:14} & 2547 & 6114 & 10414 & 2752 & 6572 & 10414 & 2752 & 6572 & 10414 & 2696 & 6572 & 10414  \\
\cite{band:14} & 3625 & 7016 & 9843 & 3542 & 6835 & 10585 & 2977 & 6776 & 9843 & 3336 & 6632 & 9843  \\
\cite{mayo:13} & 3171 & 6835 & 9843 & 3298 & 6835 & 9843 & 3473 & 7148 & 9843 & 3211 & 6835 & 9843  \\
\cite{lees:11} & 2628 & 6114 & 10414 & 2792 & 6572 & 10414 & 2792 & 6572 & 10414 & 2737 & 6572 & 10414  \\
\cite{ande:07} & 2718 & 6923 & 9928 & 2757 & 6572 & 9928 & 2757 & 6572 & 9928 & 2680 & 6572 & 9928  \\
\cite{lees:17} & 5999 & 10148 & 21020 & 5799 & 10148 & 21020 & 5654 & 10148 & 21020 & 5697 & 10148 & 21020  \\
\bottomrule
\end{tabular}
}
\caption{Allocation costs by heuristic and fleet size}\label{costallocdetnoroll}
\end{table}

\begin{table}[ht!]
\scalebox{0.7}{
\centering
\begin{tabular}{@{}l|ccc|ccc|ccc|ccc@{}}
\toprule
Number of Ambulances & \multicolumn{3}{c|}{10} & \multicolumn{3}{c|}{12} & \multicolumn{3}{c|}{14} & \multicolumn{3}{c}{18} \\
\midrule
Heuristic & Mean & Q0.9 & Max & Mean & Q0.9 & Max & Mean & Q0.9 & Max & Mean & Q0.9 & Max \\
\midrule
CA & 4951 & 9335 & 12799 & 2537 & 5455 & 9365 & 1963 & 4003 & 6218 & 1577 & 3906 & 4947  \\
BM & 2800 & 4957 & 8704 & 2029 & 4279 & 6629 & 1665 & 3904 & 4947 & 1597 & 3882 & 4947  \\
NM & 2459 & 5248 & 10347 & 1869 & 3991 & 4947 & 1713 & 3882 & 4996 & 1624 & 3882 & 4947  \\
GHP1 & 4203 & 8128 & 10347 & 3616 & 6910 & 9932 & 2215 & 4606 & 6975 & 1601 & 3882 & 6218  \\
GHP2 & 2586 & 5496 & 7265 & 2580 & 4690 & 7759 & 1898 & 4296 & 5512 & 1725 & 3882 & 4947  \\
\cite{jagt:17a} & 4822 & 8987 & 14168 & 3128 & 6669 & 9312 & 2240 & 4530 & 7396 & 1963 & 4240 & 7396  \\
\cite{lees:14} & 3555 & 6265 & 10347 & 2855 & 6004 & 8414 & 2821 & 6359 & 11872 & 1828 & 4000 & 8062  \\
\cite{band:14} & 4532 & 8350 & 12103 & 3105 & 6749 & 10380 & 1955 & 3993 & 4996 & 1992 & 3906 & 6218  \\
\cite{mayo:13} & 4197 & 7726 & 9933 & 3134 & 6388 & 10287 & 2336 & 4740 & 7578 & 1705 & 3882 & 4947  \\
\cite{lees:11} & 3640 & 7793 & 10347 & 2298 & 4796 & 7457 & 2181 & 4206 & 7624 & 1806 & 3882 & 6218  \\
\cite{ande:07} & 4229 & 6912 & 11872 & 2402 & 4796 & 7265 & 1967 & 4063 & 6218 & 1809 & 3978 & 6218  \\
\cite{lees:17} & 3881 & 8340 & 10964 & 3274 & 7073 & 9437 & 3078 & 6992 & 10347 & 3448 & 7074 & 13185  \\
\bottomrule
\end{tabular}
}
\scalebox{0.7}{
\begin{tabular}{@{}l|ccc|ccc|ccc|ccc@{}}
\toprule
Number of Ambulances & \multicolumn{3}{c|}{22} & \multicolumn{3}{c|}{24} & \multicolumn{3}{c|}{26} & \multicolumn{3}{c}{30} \\
\midrule
Heuristic & Mean & Q0.9 & Max & Mean & Q0.9 & Max & Mean & Q0.9 & Max & Mean & Q0.9 & Max \\
\midrule
CA & 1628 & 3906 & 4947 & 1678 & 3906 & 4947 & 1709 & 3906 & 4947 & 1683 & 3906 & 4947  \\
BM & 1488 & 3882 & 4947 & 1498 & 3882 & 4996 & 1496 & 3882 & 4947 & 1490 & 3882 & 4947  \\
NM & 1488 & 3882 & 4947 & 1498 & 3882 & 4996 & 1469 & 3882 & 4996 & 1464 & 3882 & 4947  \\
GHP1 & 1766 & 3882 & 4947 & 1533 & 3906 & 4996 & 1544 & 3906 & 4996 & 1521 & 3882 & 4947  \\
GHP2 & 1597 & 3882 & 4947 & 1514 & 3882 & 4996 & 1542 & 3906 & 4996 & 1514 & 3882 & 4947  \\
\cite{jagt:17a} & 1761 & 3906 & 4996 & 1676 & 3906 & 4996 & 1684 & 3906 & 4996 & 1651 & 3882 & 4996  \\
\cite{lees:14} & 1645 & 3906 & 4808 & 1962 & 4341 & 6461 & 1734 & 3919 & 4808 & 1630 & 4050 & 4947  \\
\cite{band:14} & 1794 & 3906 & 4947 & 1705 & 3882 & 5885 & 1724 & 4000 & 4947 & 1735 & 4000 & 4947  \\
\cite{mayo:13} & 1887 & 4000 & 5885 & 1770 & 3882 & 5885 & 1697 & 3906 & 6218 & 1668 & 3882 & 6218  \\
\cite{lees:11} & 1741 & 3906 & 4808 & 1582 & 3906 & 4947 & 1598 & 3906 & 4947 & 1713 & 3882 & 4947  \\
\cite{ande:07} & 1667 & 3882 & 4996 & 1695 & 4000 & 4947 & 1757 & 4222 & 4808 & 1727 & 4429 & 4808  \\
\cite{lees:17} & 3599 & 5408 & 12611 & 3247 & 4972 & 11771 & 2031 & 4218 & 10347 & 2962 & 4581 & 10930  \\
\bottomrule
\end{tabular}
}
\caption{Allocation costs by  rollout heuristic and fleet size}\label{costallocdetroll}
\end{table}

\begin{table}[ht!]
\scalebox{0.7}{
\centering
\begin{tabular}{@{}l|ccc|ccc|ccc|ccc@{}}
\toprule
Number of Ambulances & \multicolumn{3}{c|}{10} & \multicolumn{3}{c|}{12} & \multicolumn{3}{c|}{14} & \multicolumn{3}{c}{18} \\
\midrule
Heuristic & Mean & Q0.9 & Max & Mean & Q0.9 & Max & Mean & Q0.9 & Max & Mean & Q0.9 & Max \\
\midrule
CA & 2480 & 3968 & 5532 & 1350 & 2829 & 4805 & 719 & 1199 & 2087 & 649 & 1136 & 1559  \\
BM & 1641 & 2438 & 3615 & 970 & 1858 & 2387 & 745 & 1273 & 1997 & 671 & 1138 & 1611  \\
NM & 1590 & 2471 & 3188 & 991 & 1787 & 2204 & 729 & 1162 & 1639 & 654 & 1060 & 1202 \\
GHP1 & 2101 & 3889 & 4927 & 1294 & 2910 & 4805 & 820 & 1255 & 4081 & 792 & 1225 & 3981  \\
GHP2 & 1497 & 2611 & 3939 & 901 & 1670 & 2423 & 744 & 1210 & 1596 & 659 & 1060 & 1221  \\
\cite{jagt:17a} & 2318 & 4588 & 5262 & 1597 & 3057 & 5262 & 1444 & 2867 & 5262 & 933 & 1810 & 4302  \\
\cite{lees:14} & 1930 & 4411 & 6266 & 1183 & 2726 & 4805 & 719 & 1199 & 2087 & 649 & 1136 & 1559  \\
\cite{band:14} & 2115 & 3703 & 5036 & 1523 & 2819 & 4039 & 1183 & 2883 & 4232 & 1247 & 2643 & 4056  \\
\cite{mayo:13} & 2706 & 4706 & 6355 & 1723 & 3095 & 5139 & 1328 & 3120 & 4915 & 1132 & 2372 & 3857  \\
\cite{lees:11} & 2063 & 4560 & 7109 & 1137 & 2726 & 4805 & 719 & 1199 & 2087 & 650 & 1136 & 1559  \\
\cite{ande:07} & 2470 & 3968 & 5532 & 1240 & 3010 & 4536 & 725 & 1199 & 2087 & 622 & 1063 & 1237  \\
\cite{lees:17} & 1983 & 3946 & 7511 & 1659 & 3092 & 5043 & 1614 & 3686 & 4970 & 1557 & 3027 & 4858  \\
\bottomrule
\end{tabular}
}
\scalebox{0.7}{
\begin{tabular}{@{}l|ccc|ccc|ccc|ccc@{}}
\toprule
Number of Ambulances & \multicolumn{3}{c|}{22} & \multicolumn{3}{c|}{24} & \multicolumn{3}{c|}{26} & \multicolumn{3}{c}{30} \\
\midrule
Heuristic & Mean & Q0.9 & Max & Mean & Q0.9 & Max & Mean & Q0.9 & Max & Mean & Q0.9 & Max \\
\midrule
CA & 604 & 1047 & 1559 & 584 & 1047 & 1559 & 584 & 1047 & 1559 & 569 & 980 & 1559  \\
BM & 607 & 1051 & 1279 & 587 & 1051 & 1279 & 585 & 1051 & 1279 & 583 & 1051 & 1279  \\
NM & 607 & 1051 & 1279 & 587 & 1051 & 1279 & 585 & 1051 & 1279 & 583 & 1051 & 1279  \\
GHP1 & 705 & 1199 & 2654 & 672 & 1199 & 2246 & 669 & 1199 & 2246 & 627 & 1136 & 2087  \\
GHP2 & 626 & 1051 & 1672 & 607 & 1051 & 1672 & 605 & 1051 & 1672 & 603 & 1051 & 1672  \\
\cite{jagt:17a} & 871 & 1647 & 2439 & 849 & 1647 & 2439 & 893 & 1795 & 2982 & 833 & 1647 & 2525  \\
\cite{lees:14} & 604 & 1047 & 1559 & 584 & 1047 & 1559 & 584 & 1047 & 1559 & 569 & 980 & 1559  \\
\cite{band:14} & 1220 & 3092 & 3817 & 1253 & 3244 & 3931 & 986 & 2372 & 3981 & 1233 & 3047 & 4922  \\
\cite{mayo:13} & 861 & 1601 & 3351 & 988 & 1738 & 3938 & 1001 & 1738 & 3938 & 941 & 1641 & 3938  \\
\cite{lees:11} & 604 & 1047 & 1559 & 584 & 1047 & 1559 & 584 & 1047 & 1559 & 569 & 980 & 1559  \\
\cite{ande:07} & 558 & 946 & 1237 & 543 & 946 & 1237 & 543 & 946 & 1237 & 525 & 925 & 1237  \\
\cite{lees:17} & 1794 & 3820 & 4880 & 1636 & 3471 & 4880 & 1547 & 2952 & 4880 & 1550 & 2932 & 4880  \\
\bottomrule
\end{tabular}
}
\caption{Response times by heuristic and fleet size}
\label{table_response_times_deterministic}
\end{table}

\begin{table}[ht!]
\scalebox{0.7}{
\centering
\begin{tabular}{@{}l|ccc|ccc|ccc|ccc@{}}
\toprule
Number of Ambulances & \multicolumn{3}{c|}{10} & \multicolumn{3}{c|}{12} & \multicolumn{3}{c|}{14} & \multicolumn{3}{c}{18} \\
\midrule
Heuristic & Mean & Q0.9 & Max & Mean & Q0.9 & Max & Mean & Q0.9 & Max & Mean & Q0.9 & Max \\
\midrule
CA & 2047 & 3697 & 4853 & 837 & 1364 & 3190 & 851 & 1610 & 3981 & 649 & 1195 & 1693  \\
BM & 1321 & 2370 & 3737 & 883 & 1685 & 2147 & 669 & 1182 & 1581 & 622 & 1075 & 1588  \\
NM & 1109 & 2198 & 5073 & 813 & 1379 & 3305 & 727 & 1310 & 1783 & 625 & 1075 & 1588 \\
GHP1 & 1837 & 3063 & 5493 & 1630 & 3025 & 4229 & 926 & 1755 & 3981 & 651 & 1302 & 1693  \\
GHP2 & 1138 & 2141 & 3130 & 1086 & 2130 & 3066 & 833 & 1358 & 2563 & 694 & 1262 & 1695  \\
\cite{jagt:17a} & 1943 & 3450 & 4750 & 1142 & 2022 & 4890 & 951 & 1884 & 3981 & 881 & 1777 & 2654  \\
\cite{lees:14} & 1448 & 2588 & 4610 & 1152 & 2140 & 3515 & 1024 & 1702 & 3275 & 703 & 1324 & 2251  \\
\cite{band:14} & 2011 & 3731 & 5418 & 1305 & 2467 & 4132 & 872 & 1273 & 3981 & 848 & 1610 & 2788  \\
\cite{mayo:13} & 1941 & 3549 & 5615 & 1155 & 2380 & 3104 & 1007 & 2061 & 3981 & 705 & 1221 & 1588  \\
\cite{lees:11} & 1206 & 2063 & 3213 & 906 & 1578 & 3515 & 839 & 1426 & 3148 & 741 & 1364 & 2251  \\
\cite{ande:07} & 1858 & 3880 & 4832 & 903 & 1541 & 3515 & 791 & 1441 & 2311 & 759 & 1302 & 3981  \\
\cite{lees:17} & 1697 & 2786 & 4275 & 1348 & 2302 & 4693 & 1302 & 2382 & 3259 & 1214 & 1932 & 3789  \\
\bottomrule
\end{tabular}
}
\scalebox{0.7}{
\begin{tabular}{@{}l|ccc|ccc|ccc|ccc@{}}
\toprule
Number of Ambulances & \multicolumn{3}{c|}{22} & \multicolumn{3}{c|}{24} & \multicolumn{3}{c|}{26} & \multicolumn{3}{c}{30} \\
\midrule
Heuristic & Mean & Q0.9 & Max & Mean & Q0.9 & Max & Mean & Q0.9 & Max & Mean & Q0.9 & Max \\
\midrule
CA & 661 & 1125 & 1237 & 683 & 1063 & 1255 & 683 & 1043 & 1237 & 656 & 1043 & 1237  \\
BM & 589 & 1116 & 1588 & 615 & 1116 & 1588 & 580 & 1000 & 1298 & 578 & 1000 & 1298  \\
NM & 589 & 1116 & 1588 & 615 & 1116 & 1588 & 594 & 1075 & 1588 & 588 & 1075 & 1588 \\
GHP1 & 834 & 1613 & 3274 & 651 & 1310 & 1693 & 657 & 1204 & 1712 & 609 & 1116 & 1588  \\
GHP2 & 720 & 1302 & 3651 & 631 & 1204 & 1588 & 640 & 1286 & 1629 & 640 & 1199 & 1648  \\
\cite{jagt:17a} & 777 & 1258 & 3060 & 730 & 1204 & 2087 & 742 & 1204 & 2616 & 674 & 1116 & 1693  \\
\cite{lees:14} & 691 & 1152 & 1969 & 799 & 1384 & 1969 & 687 & 1064 & 1512 & 663 & 1014 & 1237  \\
\cite{band:14} & 747 & 1224 & 2525 & 727 & 1266 & 2138 & 799 & 1271 & 2389 & 774 & 1389 & 2525  \\
\cite{mayo:13} & 860 & 1495 & 4817 & 745 & 1331 & 2678 & 706 & 1258 & 1588 & 687 & 1258 & 1588  \\
\cite{lees:11} & 790 & 1370 & 1969 & 648 & 1060 & 1366 & 644 & 1112 & 1617 & 725 & 1203 & 1587  \\
\cite{ande:07} & 691 & 1191 & 1969 & 696 & 1082 & 1724 & 701 & 1152 & 1617 & 671 & 1130 & 1617  \\
\cite{lees:17} & 1433 & 2367 & 3899 & 1312 & 2176 & 3913 & 884 & 2010 & 3648 & 1253 & 1969 & 3947  \\
\bottomrule
\end{tabular}
}
\caption{Response times by rollout heuristic and fleet size}
\label{table_response_times_deterministic_rollout}
\end{table}

\begin{figure}
    \centering
    \includegraphics[width=0.6\linewidth]{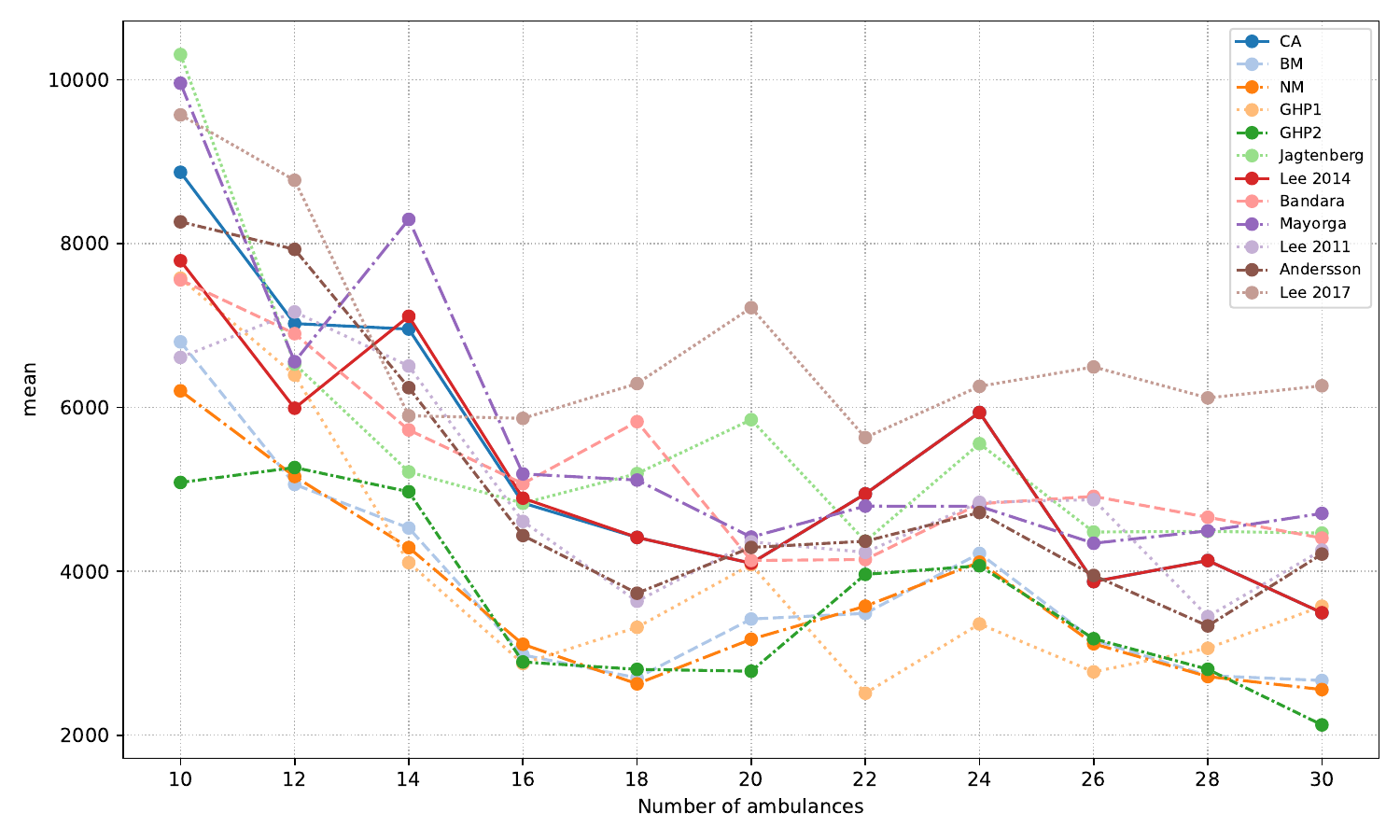}
    \caption{Mean allocation costs using random travel times.}
    \label{fig:random_allocation_costs_results}
\end{figure}

\begin{figure}
    \centering
    \includegraphics[width=0.6\linewidth]{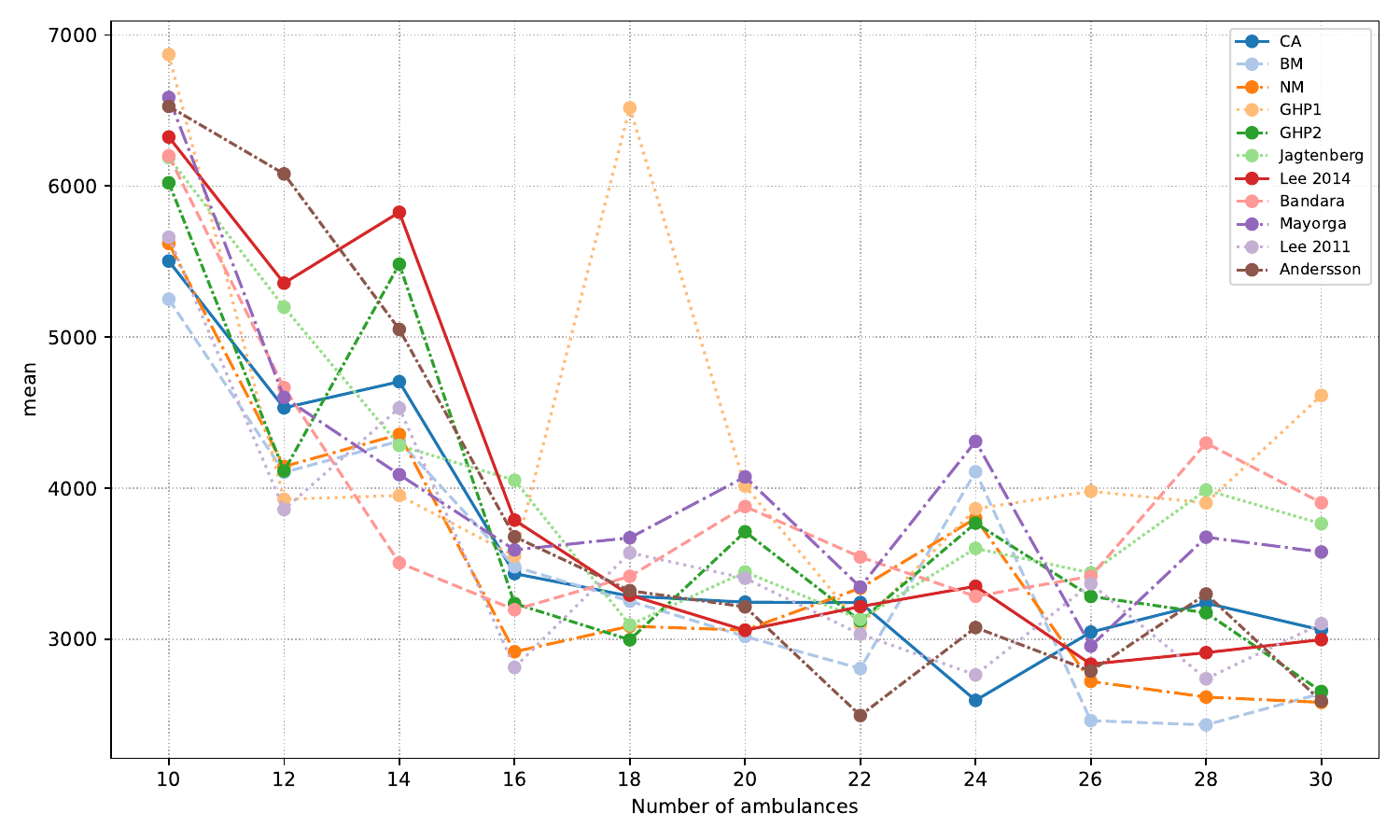}
    \caption{Mean allocation costs using random travel times with the rollout heuristics.}
    \label{fig:random_allocation_costs_results_rollout}
\end{figure}

\begin{figure}
    \centering
    \includegraphics[width=0.6\linewidth]{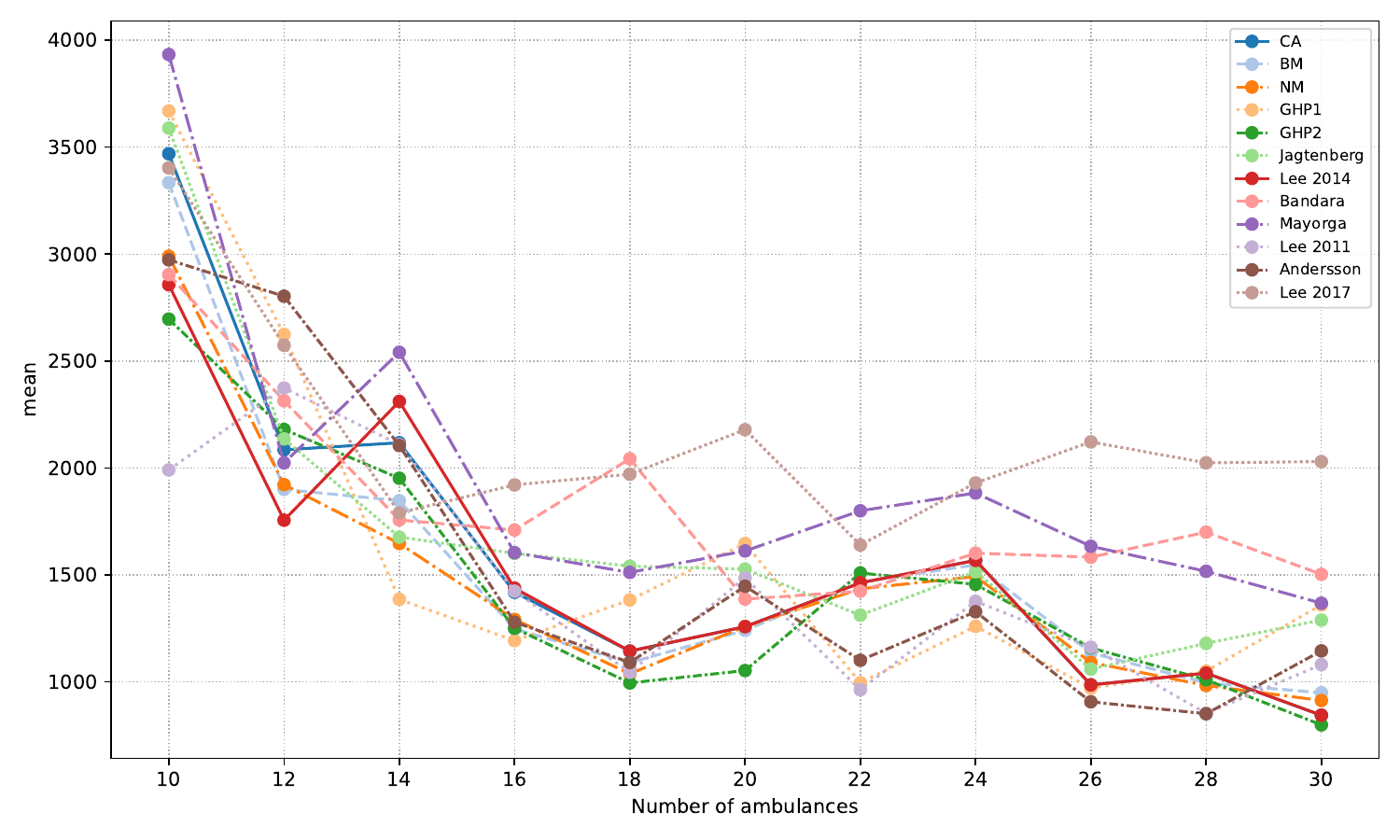}
    \caption{Mean response times using random travel times.}
    \label{fig:random_response_times_results}
\end{figure}

\begin{figure}
    \centering
    \includegraphics[width=0.6\linewidth]{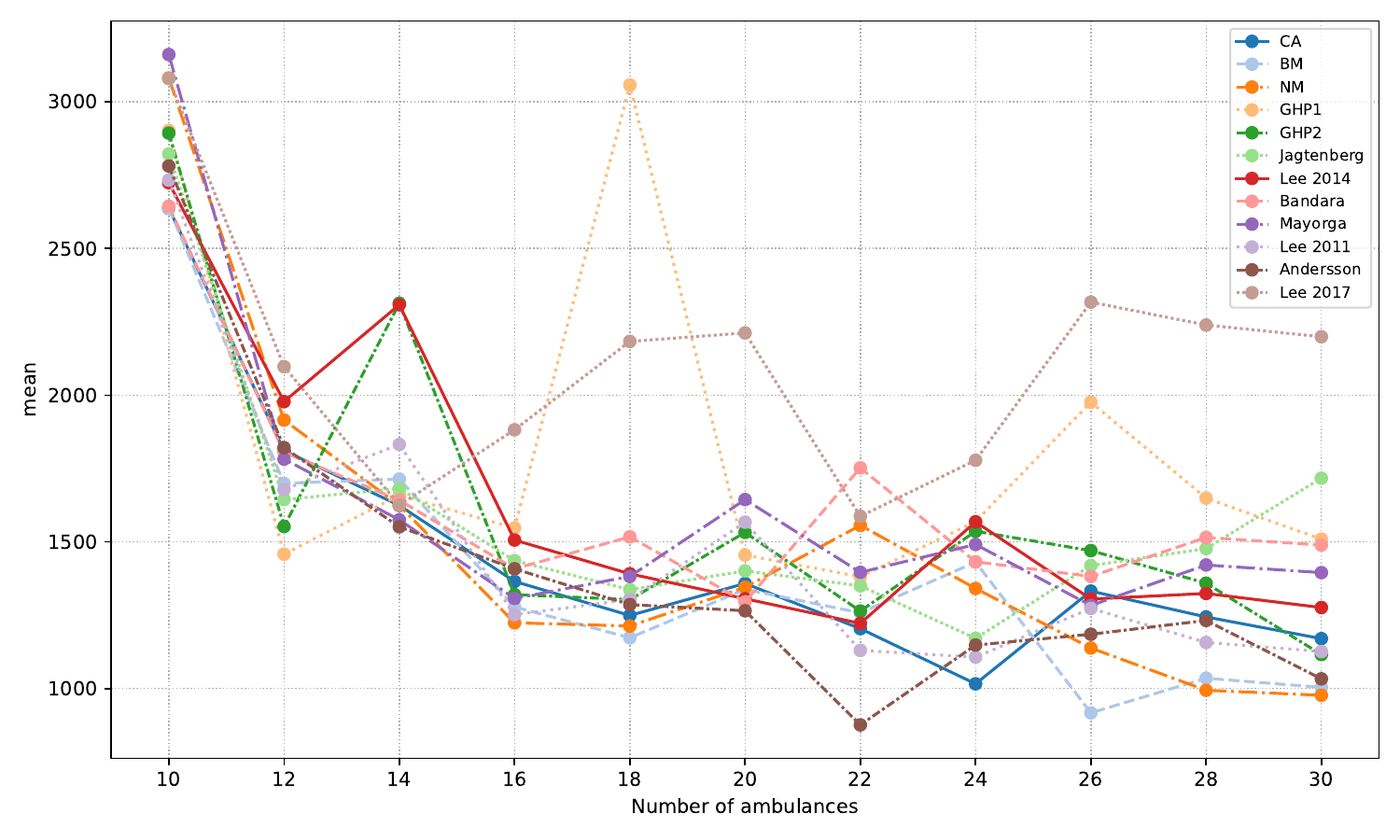}
    \caption{Mean response times using random travel times with the rollout heuristics.}
    \label{fig:random_response_times_results_rollout}
\end{figure}

\textbf{Numerical results with random travel times.} 
Travel time uncertainty is modeled using a lognormal distribution, following the approach proposed in \cite{west:16}.
Let $d$ denote the free flow travel time between two locations.
The random travel time $T$ is given by
\[
T \ \sim \ \text{LogNormal}(\mu, \sigma^2(d)),
\]
with
\[
\mu \ = \ \log(d + \tau_{0}),
\qquad
\sigma^2(d) \ = \ M e^{-\lambda d} + \delta.
\]
The intercept $\tau_{0}$ accounts for fixed delays unrelated to distance, while the variance model is obtained from the empirical behavior of travel times: shorter trips exhibit higher relative variability than longer trips.
Parameters $M$ and $\lambda$ control the magnitude and decay rate of variability, and $\delta$ represents a baseline variance floor.
For the numerical results in this paper, we use $\tau_{0} = 45$ seconds, $M = 1.05$, $\delta = 0.09$ and $\lambda = 0.001$.
Figures~\ref{fig:random_allocation_costs_results}--\ref{fig:random_response_times_results_rollout} report results with random travel times.
Figure~\ref{fig:random_allocation_costs_results} shows the mean allocation costs for each of the base heuristics and Figure~\ref{fig:random_allocation_costs_results_rollout} shows the mean allocation costs when using the rollout version of each heuristic.
Figures~\ref{fig:random_response_times_results} and \ref{fig:random_response_times_results_rollout} show mean response times for the base heuristics and mean response times for the rollout versions respectively.

\par {\textbf{Optimized choice of ambulance station.}} Next, we explore the effect of the choice of ambulance station, by comparing performance under the closest station rule (CSR) and under the best station rule (BSR) (described in Section~\ref{sec:newbase}), when combined with dispatch heuristics BM, NM, GHP1, and GHP2.
For this experiment, we consider the same peak period of Friday 6pm to 8pm as in the previous experiment, and selected $4$~zones in the west part of the city and $4$~zones in the south part of the city (these $8$~zones do not intersect).
In the first half-hour, we increased the call rates for the south zones while setting the call rates for the west zones to zero.
In the last half-hour, we did the opposite, increasing call rates for the west zones and setting the call rates in the south zones to zero.
For BSR, we used $\Delta = 90$ minutes.
Tables~\ref{us_pen_cbr_bbr_50} and \ref{us_pen_cbr_bbr_50_two_stage} report the 0.9-quantiles of the allocation costs, respectively without and with rollout.
In this experiment, CBR and BSR provided similar performance.

\begin{table}[ht!]
\centering
\scalebox{0.85}{
\begin{tabular}{@{}c|ll|ll|ll|ll@{}}
\toprule
Heuristic   & \multicolumn{2}{c|}{BM}                            & \multicolumn{2}{c|}{GHP1}                          & \multicolumn{2}{c|}{GHP2}                          & \multicolumn{2}{c}{NM}                            \\ \midrule
Number of Ambulances  & \multicolumn{1}{c}{CSR} & \multicolumn{1}{c|}{BSR} & \multicolumn{1}{c}{CSR} & \multicolumn{1}{c|}{BSR} & \multicolumn{1}{c}{CSR} & \multicolumn{1}{c|}{BSR} & \multicolumn{1}{c}{CSR} & \multicolumn{1}{c}{BSR} \\ \midrule
10 & 10405                   & 10381                    & 9527                    & 9476                     & 9952                    & 9687                     & 9300                    & 9210                    \\
12 & 6096                    & 6002                     & 6002                    & 6002                     & 6002                    & 5973                     & 5986                    & 5833                    \\
14 & 4701                    & 4649                     & 4725                    & 4698                     & 4745                    & 4737                     & 4654                    & 4649                    \\
16 & 3875                    & 3931                     & 3894                    & 3931                     & 3931                    & 3908                     & 3720                    & 3864                    \\
18 & 2396                    & 2270                     & 2466                    & 2380                     & 2489                    & 2712                     & 2342                    & 2714                    \\
20 & 2288                    & 2205                     & 2325                    & 2274                     & 2297                    & 2620                     & 2234                    & 2510                    \\ \bottomrule
\end{tabular}}
\caption{0.9-quantiles of allocation costs under dispatch heuristics BM, NM, GHP1, and GHP2, using CBR and BBR for the choice of ambulance station, for several values for the number of ambulances.}
\label{us_pen_cbr_bbr_50}
\end{table}

\begin{table}[ht!]
\centering
\scalebox{0.85}{
\begin{tabular}{@{}c|cc|cc|cc|cc@{}}
\toprule
Heuristic   & \multicolumn{2}{c|}{BM} & \multicolumn{2}{c|}{GHP1} & \multicolumn{2}{c|}{GHP2} & \multicolumn{2}{c}{NM} \\ \midrule
Number of Ambulances  & CSR        & BSR        & CSR         & BSR         & CSR         & BSR         & CSR        & BSR       \\ \midrule
10 & 7855       & 7707       & 7026        & 6968        & 7533        & 8194        & 8037       & 8015      \\
12 & 4316       & 4459       & 4654        & 4828        & 4819        & 5050        & 4597       & 4579      \\
14 & 3576       & 3332       & 3813        & 3332        & 3332        & 3824        & 3472       & 3439      \\
16 & 2664       & 2664       & 2735        & 2767        & 2680        & 2712        & 2664       & 2763      \\
18 & 1956       & 2198       & 2052        & 1904        & 2051        & 1975        & 2096       & 2002      \\
20 & 1820       & 1944       & 1831        & 1956        & 1756        & 1956        & 1719       & 1953      \\ \bottomrule
\end{tabular}}
\caption{0.9-quantiles of allocation costs under rollout with dispatch heuristics BM, NM, GHP1, and GHP2, using CBR and BBR for the choice of ambulance station, for several values for the number of ambulances.}
\label{us_pen_cbr_bbr_50_two_stage}
\end{table}

\subsection{Example with real and artificial data}

In this section, we consider a simple instance designed to make the rule for sending an idle ambulance to a station important.
In this instance, there are $4$ identical ambulances, $4$~stations, and all emergencies are of the same type.
The $4$~stations are (current) stations in the north, south, west and east parts of the city of Rio de Janeiro.
The simulation begins with one ambulance at each station.
The service time of each emergency is 45 minutes.
Calls arrive according to the following pattern: during the first 2 hours, say in $[0,2]$, there are $4$~emergencies in the south region with arrival times and locations independent uniformly distributed in time and space, and no emergencies in the remaining $3$~regions.
During the next $3$~hours, in $[2,5]$, no emergencies arrive anywhere.
This pattern is repeated $3$ times, to obtain $4$~emergencies in the south region during each of the time windows $[0,2]$, $[5,7]$, and $[10,12]$, and no emergencies during the time windows $[2,5]$, $[7,10]$, and $[12,15]$.
Next the pattern is repeated replacing the south region with the north region, with $4$~emergencies in the north region during each of the $3$ time windows $[15,17]$, $[20,22]$, and $[25,27]$, and no emergencies during time windows $[17,20]$, $[22,25]$, and $[27,30]$ anywhere in the city.
For BSR, $q(b,c,t_{1},t_{2})$ was the maximum number of emergencies during $[t_{1},t_{2}]$.
Table~\ref{table:bbr_hbr} reports the response times under heuristics CA, BM, NM, GHP1, and GHP2, with the home station rule (HBR) and best station rule (BSR) for sending idle ambulances to stations.
In this example, the response times with BSR are approximately half the response times with HBR. This example shows that we can build instances where BSR  performs much better than HSR and should in general be privileged to account for future high call loads in some regions.

\begin{table}[ht!]
\centering
\begin{tabular}{@{}ccc@{}}
\toprule
     & HSR  & BSR \\ \midrule
CA   & 1505 & 734 \\
BM   & 1468 & 734 \\
GHP1 & 1468 & 734 \\
GHP2 & 1468 & 734 \\
NM   & 1468 & 709 \\ \bottomrule
\end{tabular}
\caption{Response times (in seconds) under ambulance dispatch heuristics CA, BM, NM, GHP1, GHP2, combined with either the Home Station Rule (HSR) or the Best Station Rule (BSR) for sending idle ambulances to stations.}
\label{table:bbr_hbr}
\end{table}

\section{Conclusion and future work}

In this paper, we proposed four heuristics for ambulance selection and one
new heuristic for ambulance reassignment
for the management of an ambulance
fleet under uncertainty. These  heuristics
were combined with a rolling horizon approach
where second stage costs are computed with these
heuristics.
The resulting policies 
performed better than a large set
of recent heuristics from the literature for the corresponding problem.
The heuristics provide both quick response times for the calls
and decisions that are computed quickly. They are also nonanticipative,
as it should be, meaning that an ambulance can only be sent to a call
that has already arrived (and not to call that would arrive in the future, which would be possible in a deterministic setting where the whole set of calls is know over the entire planning horizon).

As a future work, we intend to extend our heuristics to the situation where
ambulance capacities are taken into account.

\addcontentsline{toc}{section}{References}
\bibliographystyle{plainnat}
\bibliography{New_Ambulance_Heuristics}

\begin{thebibliography}{62}
\providecommand{\natexlab}[1]{#1}
\providecommand{\url}[1]{\texttt{#1}}
\expandafter\ifx\csname urlstyle\endcsname\relax
  \providecommand{\doi}[1]{doi: #1}\else
  \providecommand{\doi}{doi: \begingroup \urlstyle{rm}\Url}\fi

\bibitem[Alanis et~al.(2013)Alanis, Ingolfsson, and Kolfal]{alan:13}
R.~Alanis, A.~Ingolfsson, and B.~Kolfal.
\newblock A {Markov} chain model for an {EMS} system with repositioning.
\newblock \emph{Production and Operations Management}, 22\penalty0
  (1):\penalty0 216--231, 2013.

\bibitem[Andersson and V\"{a}rbrand(2007)]{ande:07}
T.~Andersson and P.~V\"{a}rbrand.
\newblock Decision support tools for ambulance dispatch and relocation.
\newblock \emph{Journal of the Operational Research Society}, 58\penalty0
  (2):\penalty0 195--201, 2007.

\bibitem[Bandara et~al.(2012)Bandara, Mayorga, and McLay]{band:12}
D.~Bandara, M.~E. Mayorga, and L.~A. McLay.
\newblock Optimal dispatching strategies for emergency vehicles to increase
  patient survivability.
\newblock \emph{International Journal of Operational Research}, 15\penalty0
  (2):\penalty0 195--214, 2012.

\bibitem[Bandara et~al.(2014)Bandara, Mayorga, and McLay]{band:14}
D.~Bandara, M.~E. Mayorga, and L.~A. McLay.
\newblock Priority dispatching strategies for {EMS} systems.
\newblock \emph{Journal of the Operational Research Society}, 65:\penalty0
  572--587, 2014.

\bibitem[B\'{e}langer et~al.(2016)B\'{e}langer, Kergosien, Ruiz, and
  Soriano]{bela:16}
V.~B\'{e}langer, Y.~Kergosien, A.~Ruiz, and P.~Soriano.
\newblock An empirical comparison of relocation strategies in real-time
  ambulance fleet management.
\newblock \emph{Computers and Industrial Engineering}, 94:\penalty0 216--229,
  2016.

\bibitem[Berlin and Liebman(1974)]{berl:74}
G.~N. Berlin and J.~C. Liebman.
\newblock Mathematical analysis of emergency ambulance location.
\newblock \emph{Socio-Economic Planning Sciences}, 8\penalty0 (6):\penalty0
  323--328, 1974.

\bibitem[Brotcorne et~al.(2003)Brotcorne, Laporte, and Semet]{brot:03}
L.~Brotcorne, G.~Laporte, and F.~Semet.
\newblock Ambulance location and relocation models.
\newblock \emph{European Journal of Operational Research}, 147\penalty0
  (3):\penalty0 451--463, 2003.

\bibitem[Burwell et~al.(1993)Burwell, Jarvis, and McKnew]{burw:93}
T.~H. Burwell, J.~P. Jarvis, and M.~A. McKnew.
\newblock Modeling co-located servers and dispatch ties in the hypercube model.
\newblock \emph{Computers and Operations Research}, 20\penalty0 (2):\penalty0
  113--119, 1993.

\bibitem[Church and ReVelle(1974)]{chur:74}
R.~L. Church and C.~S. ReVelle.
\newblock The maximal covering location problem.
\newblock \emph{Papers of the Regional Science Association}, 32:\penalty0
  101--118, 1974.

\bibitem[Daskin(1983)]{dask:83}
M.~S. Daskin.
\newblock A maximum expected covering location model: Formulation, properties
  and heuristic solution.
\newblock \emph{Transportation Science}, 7\penalty0 (1):\penalty0 48--70, 1983.

\bibitem[Daskin and Stern(1981)]{dask:81}
M.~S. Daskin and E.~H. Stern.
\newblock A hierarchical objective set covering model for {Emergency Medical
  Service} vehicle deployment.
\newblock \emph{Transportation Science}, 15\penalty0 (2):\penalty0 137--152,
  1981.

\bibitem[Degel et~al.(2015)Degel, Wiesche, Rachuba, and Werners]{dege:15}
D.~Degel, L.~Wiesche, S.~Rachuba, and B.~Werners.
\newblock Time-dependent ambulance allocation considering data-driven
  empirically required coverage.
\newblock \emph{Health Care Management Science}, 18:\penalty0 444--458, 2015.

\bibitem[Erkut et~al.(2009)Erkut, Ingolfsson, Sim, and Erdo\u{g}an]{erku:09}
E.~Erkut, A.~Ingolfsson, T.~Sim, and G.~Erdo\u{g}an.
\newblock Computational comparison of five maximal covering models for locating
  ambulances.
\newblock \emph{Geographical Analysis}, 41:\penalty0 43--65, 2009.

\bibitem[Fitzsimmons(1973)]{fitz:73}
J.~A. Fitzsimmons.
\newblock A methodology for emergency ambulance deployment.
\newblock \emph{Management Science}, 19\penalty0 (6):\penalty0 627--636, 1973.

\bibitem[Gendreau et~al.(1997)Gendreau, Laporte, and Semet]{gend:97}
M.~Gendreau, G.~Laporte, and F.~Semet.
\newblock Solving an ambulance location model by {Tabu Search}.
\newblock \emph{Location Science}, 5\penalty0 (2):\penalty0 75--58, 1997.

\bibitem[Gendreau et~al.(2001)Gendreau, Laporte, and Semet]{gend:01}
M.~Gendreau, G.~Laporte, and F.~Semet.
\newblock A dynamic model and parallel tabu search heuristic for real-time
  ambulance relocation.
\newblock \emph{Parallel Computing}, 27\penalty0 (12):\penalty0 1641--1653,
  2001.

\bibitem[Gendreau et~al.(2006)Gendreau, Laporte, and Semet]{gend:06}
M.~Gendreau, G.~Laporte, and F.~Semet.
\newblock The maximal expected coverage relocation problem for emergency
  vehicles.
\newblock \emph{Journal of the Operational Research Society}, 57\penalty0
  (1):\penalty0 22--28, 2006.

\bibitem[Goldberg and Paz(1991)]{gold:91a}
J.~Goldberg and L.~Paz.
\newblock Locating emergency vehicle bases when service time depends on call
  location.
\newblock \emph{Transportation Science}, 25\penalty0 (4):\penalty0 264--280,
  1991.

\bibitem[Goldberg and Szidarovszky(1991)]{gold:91b}
J.~Goldberg and F.~Szidarovszky.
\newblock Methods for solving nonlinear equations used in evaluating emergency
  vehicle busy probabilities.
\newblock \emph{Operations Research}, 39\penalty0 (6):\penalty0 903--916, 1991.

\bibitem[Goldberg et~al.(1990)Goldberg, Dietrich, Chen, Mitwasi, Valenzuela,
  and Criss]{gold:90b}
J.~Goldberg, R.~Dietrich, J.~M. Chen, M.~G. Mitwasi, T.~Valenzuela, and
  E.~Criss.
\newblock Validating and applying a model for locating emergency medical
  vehicles in {Tucson, AZ}.
\newblock \emph{European Journal of Operational Research}, 49\penalty0
  (3):\penalty0 308--324, 1990.

\bibitem[Guigues et~al.(2022)Guigues, Kleywegt, and Nascimento]{guiklevhn2022}
V.~Guigues, A.~J. Kleywegt, and V.~H. Nascimento.
\newblock Operation of an ambulance fleet under uncertainty.
\newblock arXiv:2203.16371v2 [math.OC], 2022.

\bibitem[Guigues et~al.(2023{\natexlab{a}})Guigues, Kleywegt, Amorim, Krauss,
  and Nascimento]{laspatedmanual}
V.~Guigues, A.~J. Kleywegt, G.~Amorim, A.~M. Krauss, and V.~H. Nascimento.
\newblock {LASPATED: A Library for the Analysis of SPAtio-TEmporal Discrete
  Data (User Manual)}.
\newblock arXiv:2407.13889 [stat.CO], 2023{\natexlab{a}}.

\bibitem[Guigues et~al.(2023{\natexlab{b}})Guigues, Kleywegt, Amorim, Krauss,
  and Nascimento]{laspatedpaper}
V.~Guigues, A.~J. Kleywegt, G.~Amorim, A.~M. Krauss, and V.~H. Nascimento.
\newblock {LASPATED: A Library for the Analysis of SPAtio-TEmporal Discrete
  Data}.
\newblock arXiv:2401.04156v2 [stat.ME], 2023{\natexlab{b}}.

\bibitem[Guigues et~al.(2024)Guigues, Kleywegt, Nascimento, Salles, and
  Viana]{websiteambrouting24}
V.~Guigues, A.~Kleywegt, V.~H. Nascimento, Victor Salles, and Thais Viana.
\newblock Management and visualization tools for emergency medical services.
\newblock arXiv, 2024.

\bibitem[Henderson and Mason(1999)]{hend:99}
S.~G. Henderson and A.~J. Mason.
\newblock Estimating ambulance requirements in {Auckland, New Zealand}.
\newblock In \emph{Proceedings of the 1999 Winter Simulation Conference},
  volume~2, pages 1670--1674, 1999.

\bibitem[Henderson and Mason(2004)]{hend:04}
S.~G. Henderson and A.~J. Mason.
\newblock Ambulance service planning: Simulation and data visualisation.
\newblock In M.~Brandeau, F.~Sainfort, and W.~Pierskalla, editors,
  \emph{Operations Research and Health Care: A Handbook of Methods and
  Applications, International Series in Operations Research and Management
  Science 70}, chapter~4, pages 77--102. Kluwer, Dordecht, 2004.

\bibitem[Hill et~al.(1984)Hill, Hill, and Jacobs]{hill:84}
E.~D. Hill, J.~L. Hill, and L.~M. Jacobs.
\newblock Planning for emergency ambulance service systems.
\newblock \emph{The Journal of Emergency Medicine}, 1:\penalty0 331--338, 1984.

\bibitem[Hogan and Revelle(1986)]{hoga:86}
K.~Hogan and C.~S. Revelle.
\newblock Concepts and applications of backup coverage.
\newblock \emph{Management Science}, 32\penalty0 (11):\penalty0 1434--1444,
  1986.

\bibitem[Ingolfsson et~al.(2008)Ingolfsson, Budge, and Erkut]{ingo:08}
A.~Ingolfsson, S.~Budge, and E.~Erkut.
\newblock Optimal ambulance location with random delays and travel times.
\newblock \emph{Health Care Management Science}, 11:\penalty0 262--274, 2008.

\bibitem[Jagtenberg et~al.(2015)Jagtenberg, Bhulai, and Van~der Mei]{jagt:15}
C.~J. Jagtenberg, S.~Bhulai, and R.~D. Van~der Mei.
\newblock An efficient heuristic for real-time ambulance redeployment.
\newblock \emph{Operations Research for Health Care}, 12\penalty0 (4):\penalty0
  27--35, 2015.

\bibitem[Jagtenberg et~al.(2017{\natexlab{a}})Jagtenberg, Bhulai, and Van~der
  Mei]{jagt:17a}
C.~J. Jagtenberg, S.~Bhulai, and R.~D. Van~der Mei.
\newblock Dynamic ambulance dispatching: Is the closest-idle policy always
  optimal?
\newblock \emph{Operations Research for Health Care}, 20\penalty0 (4):\penalty0
  517--531, 2017{\natexlab{a}}.

\bibitem[Jagtenberg et~al.(2017{\natexlab{b}})Jagtenberg, Van~den Berg, and
  Van~der Mei]{jagt:17b}
C.~J. Jagtenberg, P.~L. Van~den Berg, and R.~D. Van~der Mei.
\newblock Benchmarking online dispatch algorithms for {Emergency Medical
  Services}.
\newblock \emph{European Journal of Operational Research}, 258\penalty0
  (2):\penalty0 715--725, 2017{\natexlab{b}}.

\bibitem[Jarvis(1985)]{jarv:85}
J.~P. Jarvis.
\newblock Approximating the equilibrium behavior of multi-server loss systems.
\newblock \emph{Management Science}, 31\penalty0 (2):\penalty0 235--239, 1985.

\bibitem[Knight et~al.(2012)Knight, Harper, and Smith]{knig:12}
V.~A. Knight, P.~R. Harper, and L.~Smith.
\newblock Ambulance allocation for maximal survival with heterogeneous outcome
  measures.
\newblock \emph{Omega}, 40:\penalty0 918--926, 2012.

\bibitem[Larson(1974)]{lars:74}
R.~C. Larson.
\newblock A hypercube queuing model for facility location and redistricting in
  urban emergency services.
\newblock \emph{Computers and Operations Research}, 1:\penalty0 67--95, 1974.

\bibitem[Larson(1975)]{lars:75}
R.~C. Larson.
\newblock Approximating the performance of urban emergency service systems.
\newblock \emph{Operations Research}, 23\penalty0 (5):\penalty0 845--868, 1975.

\bibitem[Lee(2011)]{lees:11}
S.~Lee.
\newblock The role of preparedness in ambulance dispatching.
\newblock \emph{Journal of the Operational Research Society}, 62\penalty0
  (10):\penalty0 1888--1897, 2011.

\bibitem[Lee(2012)]{lees:12}
S.~Lee.
\newblock The role of centrality in ambulance dispatching.
\newblock \emph{Decision Support Systems}, 54\penalty0 (1):\penalty0 282--291,
  2012.

\bibitem[Lee(2014)]{lees:14}
S.~Lee.
\newblock Role of parallelism in ambulance dispatching.
\newblock \emph{IEEE Transactions on Systems, Man, and Cybernetics: Systems},
  44\penalty0 (8):\penalty0 1113--1122, 2014.

\bibitem[Lee(2017)]{lees:17}
S.~Lee.
\newblock A new preparedness policy for {EMS} logistics.
\newblock \emph{Health Care Management Science}, 20:\penalty0 105--114, 2017.

\bibitem[Li and Saydam(2016)]{lisay:16}
X.~Li and C.~Saydam.
\newblock Balancing ambulance crew workloads via a tiered dispatch policy.
\newblock \emph{Pesquisa Operacional}, 36\penalty0 (3):\penalty0 399--419,
  2016.

\bibitem[Mason(2013)]{maso:13}
A.~J. Mason.
\newblock Simulation and real-time optimised relocation for improving ambulance
  operations.
\newblock In B.~T. Denton, editor, \emph{Handbook of Healthcare Operations
  Management: Methods and Applications, International Series in Operations
  Research and Management Science 184}, chapter~11, pages 289--317. Springer,
  New York, 2013.

\bibitem[Maxwell et~al.(2009)Maxwell, Henderson, and Topaloglu]{maxw:09}
M.~S. Maxwell, S.~G. Henderson, and H.~Topaloglu.
\newblock Ambulance redeployment: An approximate dynamic programming approach.
\newblock In M.~D. Rossetti, R.~R. Hill, B.~Johansson, A.~Dunkin, and R.~G.
  Ingalls, editors, \emph{Proceedings of the 2009 Winter Simulation
  Conference}, pages 1850--1860, 2009.

\bibitem[Maxwell et~al.(2010)Maxwell, Restrepo, Henderson, and
  Topaloglu]{maxw:10}
M.~S. Maxwell, M.~Restrepo, S.~G. Henderson, and H.~Topaloglu.
\newblock Approximate dynamic programming for ambulance redeployment.
\newblock \emph{INFORMS Journal on Computing}, 22\penalty0 (2):\penalty0
  266--281, 2010.

\bibitem[Maxwell et~al.(2013)Maxwell, Henderson, and Topaloglu]{maxw:13}
M.~S. Maxwell, S.~G. Henderson, and H.~Topaloglu.
\newblock Tuning approximate dynamic programming policies for ambulance
  redeployment via direct search.
\newblock \emph{Stochastic Systems}, 3\penalty0 (2):\penalty0 322--361, 2013.

\bibitem[Maxwell et~al.(2014)Maxwell, Ni, Tong, Henderson, Topaloglu, and
  Hunter]{maxw:14}
M.~S. Maxwell, E.~C. Ni, C.~Tong, S.~G. Henderson, H.~Topaloglu, and S.~R.
  Hunter.
\newblock A bound on the performance of an optimal ambulance redeployment
  policy.
\newblock \emph{Operations Research}, 62\penalty0 (5):\penalty0 1014--1027,
  2014.

\bibitem[Mayorga et~al.(2013)Mayorga, Bandara, and McLay]{mayo:13}
M.~E. Mayorga, D.~Bandara, and L.~A. McLay.
\newblock Districting and dispatching policies for emergency medical service
  systems to improve patient survival.
\newblock \emph{IIE Transactions on Healthcare Systems Engineering}, 3\penalty0
  (1):\penalty0 39--56, 2013.

\bibitem[Nair and Miller-Hooks(2009)]{nair:09}
R.~Nair and E.~Miller-Hooks.
\newblock Evaluation of relocation strategies for emergency medical service
  vehicles.
\newblock \emph{Transportation Research Record}, 2137:\penalty0 63--73, 2009.

\bibitem[Repede and Bernardo(1994)]{repe:94}
J.~F. Repede and J.~J. Bernardo.
\newblock Developing and validating a decision support system for locating
  emergency medical vehicles in {Louisville, Kentucky}.
\newblock \emph{European Journal of Operational Research}, 75:\penalty0
  567--581, 1994.

\bibitem[Restrepo et~al.(2009)Restrepo, Henderson, and Topaloglu]{rest:09}
M.~Restrepo, S.~G. Henderson, and H.~Topaloglu.
\newblock Erlang loss models for the static deployment of ambulances.
\newblock \emph{Health Care Management Science}, 12\penalty0 (67):\penalty0
  67--79, 2009.

\bibitem[ReVelle and Hogan(1989)]{reve:89}
C.~ReVelle and K.~Hogan.
\newblock The maximum availability location problem.
\newblock \emph{Transportation Science}, 23\penalty0 (3):\penalty0 192--200,
  1989.

\bibitem[Schilling et~al.(1979)Schilling, Elzinga, Cohon, Church, and
  ReVelle]{schi:79}
D.~A. Schilling, D.~J. Elzinga, J.~Cohon, R.~L. Church, and C.~S. ReVelle.
\newblock The {TEAM/FLEET} models for simultaneous facility and equipment
  siting.
\newblock \emph{Transportation Science}, 13\penalty0 (2):\penalty0 163--175,
  1979.

\bibitem[Schmid(2012)]{schm:12}
V.~Schmid.
\newblock Solving the dynamic ambulance relocation and dispatching problem
  using approximate dynamic programming.
\newblock \emph{European Journal of Operational Research}, 219:\penalty0
  611--621, 2012.

\bibitem[Schmid and Doerner(2010)]{schm:10}
V.~Schmid and K.~F. Doerner.
\newblock Ambulance location and relocation problems with time-dependent travel
  times.
\newblock \emph{European Journal of Operational Research}, 207:\penalty0
  1293--1303, 2010.

\bibitem[Sorensen and Church(2010)]{sore:10}
P.~Sorensen and R.~Church.
\newblock Integrating expected coverage and local reliability for emergency
  medical services location problems.
\newblock \emph{Socio-Economic Planning Sciences}, 44\penalty0 (1):\penalty0
  8--18, 2010.

\bibitem[Swoveland et~al.(1973{\natexlab{a}})Swoveland, Uyeno, Vertinsky, and
  Vickson]{swov:73a}
C.~Swoveland, D.~Uyeno, I.~Vertinsky, and R.~Vickson.
\newblock Ambulance location: A probabilistic enumeration approach.
\newblock \emph{Management Science}, 20\penalty0 (4):\penalty0 686--698,
  1973{\natexlab{a}}.

\bibitem[Swoveland et~al.(1973{\natexlab{b}})Swoveland, Uyeno, Vertinsky, and
  Vickson]{swov:73b}
C.~Swoveland, D.~Uyeno, I.~Vertinsky, and R.~Vickson.
\newblock A simulation-based methodology for optimization of ambulance service
  policies.
\newblock \emph{Socio-Economic Planning Sciences}, 7\penalty0 (6):\penalty0
  697--703, 1973{\natexlab{b}}.

\bibitem[Toregas et~al.(1971)Toregas, Swain, ReVelle, and Bergman]{tore:71}
C.~Toregas, R.~Swain, C.~ReVelle, and L.~Bergman.
\newblock The location of emergency service facilities.
\newblock \emph{Operations Research}, 19\penalty0 (6):\penalty0 1363--1373,
  1971.

\bibitem[Van~Barneveld et~al.(2016)Van~Barneveld, Bhulai, and Van~der
  Mei]{vanbarn:16}
T.~C. Van~Barneveld, S.~Bhulai, and R.~D. Van~der Mei.
\newblock The effect of ambulance relocations on the performance of ambulance
  service providers.
\newblock \emph{European Journal of Operational Research}, 252:\penalty0
  257--269, 2016.

\bibitem[Van~Barneveld et~al.(2017)Van~Barneveld, Bhulai, and Van~der
  Mei]{vanbarn:17}
T.~C. Van~Barneveld, S.~Bhulai, and R.~D. Van~der Mei.
\newblock A dynamic ambulance management model for rural areas: Computing
  redeployment actions for relevant performance measures.
\newblock \emph{Health Care Management Science}, 20:\penalty0 165--186, 2017.

\bibitem[Volz(1971)]{volz:71}
R.~A. Volz.
\newblock Optimum ambulance location in semi-rural areas.
\newblock \emph{Transportation Science}, 5\penalty0 (2):\penalty0 193--203,
  1971.

\bibitem[Westgate et~al.(2016)Westgate, Woodard, Matteson, and
  Henderson]{west:16}
Bradford~S. Westgate, Dawn~B. Woodard, David~S. Matteson, and Shane~G.
  Henderson.
\newblock Large-network travel time distribution estimation for ambulances.
\newblock \emph{European Journal of Operational Research}, 252\penalty0
  (1):\penalty0 322--333, 2016.
\newblock ISSN 0377-2217.
\newblock \doi{https://doi.org/10.1016/j.ejor.2016.01.004}.
\newblock URL
  \url{https://www.sciencedirect.com/science/article/pii/S0377221716000102}.

\end{thebibliography}

\end{document}